\documentclass[a4paper,oneside,11pt]{article} 

\usepackage{mystyle}

\title{On heavy-tail phenomena in some large deviations problems}  
\author{Fanny Augeri\footnote{Institut de Mathématiques de Toulouse, France, E-mail: \href{mailto:faugeri@math.univ-toulouse.fr}{faugeri@math.univ-toulouse.fr} } } 
\date{\today} 

\begin{document}

\maketitle
\begin{abstract}

In this paper, we revisit the proof of the large deviations principle of Wiener chaoses partially given by Borell \cite{Borell}, and then by Ledoux \cite{LedouxWiener} in its full form. We show that some heavy-tail phenomena observed in large deviations can be explained by the same mechanism as for the Wiener chaoses, meaning that the deviations are created, in a sense, by translations. More precisely, we prove a general large deviations principle for a certain class of functionals $f_n : \RR^n \to \mathcal{X}$, where $\mathcal{X}$ is some metric space, under the $n$-fold probability measure $\nu_{\alpha}^n$, where $\nu_{\alpha} =Y_{\alpha}^{-1}e^{-|x|^{\alpha}}dx$, $\alpha \in (0,2]$, for which the large deviations are due to translations. We retrieve, as an application, the large deviations principles known for the Wigner matrices without Gaussian tails in \cite{Bordenave}, \cite{LDPei}, \cite{LDPtr} of the empirical spectral measure, the largest eigenvalue, and traces of polynomials. We also apply our large deviations result to the last-passage time, which yields a large deviations principle when the weights follow the law $Z_{\alpha}^{-1} e^{-x^{\alpha}} \Car_{x\geq 0}dx$, with $\alpha \in (0,1)$.

\end{abstract}
\section{Introduction}
In \cite{LedouxWiener}, Ledoux proposed a large deviations principle for the Wiener chaoses based on the approach Borell gave in \cite{Borell} for estimating their tail distribution. The main feature which stands out of the proof is that the large deviations of Wiener chaoses are due to translations by elements of the Cameron-Martin space. The lower bound consists in an application of the Cameron-Martin formula, whereas the upper bound relies on the Gaussian isoperimetric inequality. 

More precisely, let $(E,\mathcal{H}, \mu)$ be an abstract Wiener space, where $E$ is a separable Banach space, $\mu$ is a Gaussian measure on $E$, and $\mathcal{H}$ the reproducing kernel (see \cite{Lifshits} or \cite[chapter 4]{Ledouxflour} for proper definitions). Let also $\Psi$ be a homogenous Wiener chaos of degree $d$ taking values in some Banach space $B$, that is, a random variable in the subspace spanned in $L^2(\mu; B)$ by Hermite polynomials of degree $d$. From \cite{LedouxWiener}, we know that $t^{-d} \Psi$ follows a large deviations principle with speed $t^2$ and good rate function $I_{\Psi}$ defined by,
\begin{equation} \label{ratefuncWiener} \forall x \in B, \ I_{\Psi}(x) = \inf\Big \{ \frac{1}{2} |h|^2 : x = \Psi^{(d)}(h),  h \in \mathcal{H} \Big\},\end{equation}
where $| \ |$ denotes the norm of the reproducing kernel $\mathcal{H}$, and 
\begin{equation} \label{detequivWiener} \forall h \in \mathcal{H}, \ \Psi^{(d)}(h) = \int \Psi(x+h) d\mu(x).\end{equation}

We believe Borell and Ledoux's approach to be extremely fruitful, and can shed a new light on heavy-tail phenomena appearing in the large deviations of certain models, where the large deviations are created also, in a sense, by translations. We already used this approach in a previous work \cite{LDPtr} to deal with the question of the large deviations of traces of powers of Gaussian Wigner matrices. Indeed, this problem can be reformulated as understanding the large deviations of Gaussian chaoses defined on spaces with growing dimension. Although this problem cannot be solved directly by using the large deviations principle of Wiener chaoses, the same outline of proof was carried out in this case, and yields a rate function having a similar structure as \eqref{ratefuncWiener}.

We would like here to push further this approach in a more general setting, and give some elements showing that heavy-tail phenomena in the large deviations of certain models can be understood using the paradigm of the Wiener chaoses. To this end, we propose a general large deviations result for a certain class of functionals $f_n : \RR^n \to \mathcal{X}$, where $\mathcal{X}$ is some metric space, under the $n$-fold probability measure $\nu_{\alpha}^n$, where $\nu_{\alpha} = Y_{\alpha}^{-1} e^{-|x|^{\alpha}} dx$, with $\alpha \in (0,2]$, for which the large deviations are governed by translations.

As an application of this result, we will retrieve the large deviations principles of different spectral functionals of the so-called Wigner matrices without Gaussian tails. Introduced in  \cite{Bordenave} by Bordenave and Caputo,  the model of Wigner matrices without Gaussian tails designates Wigner matrices whose entries have tail distributions behaving like $e^{-c t^{\alpha}}$, with $c>0$, and $\alpha \in (0,2)$. This model gives rise to a heavy-tail phenomenon which enables one to derive full large deviations principles for the spectral measure \cite{Bordenave} (see \cite{Groux} in the Wishart matrix case), the largest eigenvalue \cite{LDPei}, and the traces of powers \cite{LDPtr}.

In the more restricted setting where we assume that the entries have a density with respect to Lebesgue measure which is proportional to $e^{-c |x|^{\alpha}}$, with $c>0$, and $\alpha\in (0,2)$, the large deviations principles of these spectral functionals will fall in a unified way from our general large deviation result.  

Another application of this result will consist in a large deviations principle for the last-passage time when the weights are independent and have a density on $\RR^+$ proportional to $e^{-x^{\alpha}}$ for $\alpha \in (0,1)$.

\section{Main results}
Let us present the main results of this paper.
For $\alpha>0$, we denote by $\nu_{\alpha}$ the probability measure on $\RR$ with density $Y_{\alpha}^{-1} e^{-|x|^{\alpha}}$ with respect to Lebesgue measure, and $\nu_{\alpha}^n$ its $n$-fold product measure on $\RR^n$. Similarly, we define $\mu_{\alpha}$ the probability measure on $\RR^+$ with density $Z_{\alpha}^{-1}e^{-x^{\alpha}}$. We will denote for any $h\in \RR^n$,
$$|| h ||_{\ell^{\alpha}} = \big(\sum_{i=1}^n |x_i|^{\alpha}\big )^{1/\alpha}.$$
We recall that a sequence of random variables $(Z_n)_{n\in \NN}$ taking value in some topological space $\mathcal{X}$ equipped with the Borel $\sigma$-field $\mathcal{B}$, follows a large deviations principle (LDP) with speed $\upsilon(n)$, and rate function $J : \mathcal{X} \to [0, +\infty]$, if $J$ is lower semicontinuous and $\upsilon(n)$ increases to infinity and for all $B\in \mathcal{B}$,
$$- \inf_{B^{\circ}}J \leq \liminf_{n\to +\infty} \frac{1}{\upsilon(n)} \log \PP\left(Z_n \in B\right) \leq \limsup_{n\to +\infty} \frac{1}{\upsilon(n)}\log\PP\left(Z_n \in B\right) \leq -\inf_{\overline{B}  } J,$$
where $B^{\circ}$ denotes the interior of $B$ and $\overline{B}$ the closure of $B$. We recall that $J$ is lower semicontinuous if its $t$-level sets $\{ x \in \mathcal{X} : J(x) \leq t \}$ are closed, for any $t\in [0,+\infty)$. Furthermore, if all the level sets are compact, then we say that $J$ is a good rate function.

The purpose of the general large deviations result we will present, is to identify a class of functionals $f_n : \RR^n \to \mathcal{X}$, where $\mathcal{X}$ is some metric space, for which the large deviations are created by translations. Let us describe first informally the assumptions we will make. Let $X_n$ follow the law $\nu_{\alpha}^n$. We will assume that $f_n(X_n)$ admits a kind of deterministic equivalent under additive deformations, given by a certain function $F_n$, that is,
\begin{equation} \label{equidetinfomel}f_n(X_n+v(n)^{1/\alpha}h_n) \simeq F_n(h_n),\end{equation}
in probability, for any sequence $h_n\in \RR^n$, $\sup_n||h_n||_{\ell^{\alpha}}<+\infty$, where $v(n)$ will eventually be the speed of deviations. It is convenient to think of $F_n(h_n)$ as a deterministic equivalent of $f_n(X_n+v(n)^{1/\alpha}h_n)$, where we took the large $n$ limit on the variable $X_n$. Under this assumption, we will show that a large deviations lower bound for $f_n(X_n)$ at speed $v(n)$, holds with rate function,
$$J_{\alpha} =  \sup_{\delta >0} \limsup_{\underset{n\in N}{n\to +\infty}} I_{n,\delta},$$
where
$$ \forall x \in \mathcal{X}, \ I_{n,\delta }(x) = \inf \{ ||h||_{\ell^{\alpha}}^{\alpha} : d(F_n(h),x)<\delta, h\in \RR^n \}.$$
This rate function $J_{\alpha}$ can be interpreted by saying that to make a deviation around some $F_n(h_n)$, $X_n$ needs to make a translation by $v(n)^{1/\alpha}h_n$, which one pays at the exponential scale $v(n)$ by $|| h_n ||_{\ell^{\alpha}}^{\alpha}$. 

For the upper bound, we will further assume that for any $r>0$, the deterministic equivalent \eqref{equidetinfomel} holds uniformly in $||h_n||_{\ell^{\alpha}}\leq r$. The upper bound will rely on  sharp large deviation inequalities for $\nu_{\alpha}^n$, where we will need, excepted in the Gaussian 
case, to neglect the Euclidean enlargements appearing naturally. We thus make the assumption that $f_n$ has a small, in expectation, local Lipschitz constant with respect to $|| \ ||_{\ell^2}$ when $\alpha<2$. Finally, under some compactness property of $F_n$, we will prove that a large deviations upper bound holds for $f_n(X_n)$ with speed $v(n)$ and rate function,
 $$I_{\alpha} =  \sup_{\delta >0} \inf_{n\in N} I_{n,\delta}.$$
Thus, if we moreover assume that the upper bound rate function $I_{\alpha}$ matches the lower rate function, we will get a full large deviations principle with speed $v(n)$.
More precisely, we will prove the following result.

\begin{The}\label{theoremgene}
Let $(\mathcal{X},d)$ be a metric space. Let $\alpha \in (0,2]$ and $N\subset \NN$ an infinite subset. Let $X_n$ be a random variable with law $\nu_{\alpha}^{ n}$. Let $f_n, F_n : \RR^n \to \mathcal{X}$ be measurable functions. Let $(v(n))_{n\in N}$ be a sequence going to $+\infty$. Define for $\delta>0$ and $n\in N$, the function
$$ \forall x \in \mathcal{X}, \ I_{n,\delta }(x) = \inf \{ ||h||_{\ell^{\alpha}}^{\alpha} : d(F_n(h),x)<\delta, h\in \RR^n \}.$$
We set 
\begin{equation} \label{deftaux}\forall x \in \mathcal{X}, \ I_{\alpha }(x)  = \sup_{\delta >0} \inf_{n\in N} I_{n,\delta}(x).\end{equation}
We assume:  \\
(i).(Uniform deterministic equivalent). For any $r>0$,
$$ \sup_{ h_n \in r B_{\ell^{\alpha}}} d \big(f_n(X_n+v(n)^{1/\alpha} h_n), F_n(h_n)\big) \underset{\underset{n\in N}{n\to +\infty}}{\longrightarrow} 0,$$
in probability.\\
(ii).(Control of the Lipschitz constant). If $\alpha <2$, then for any $\delta>0$ and $r>0$, there is a sequence $t_{\delta}(n)$ such that, $$ \EE \sup_{||h||_{\ell^2} \leq t_{\delta}(n)} \mathcal{L}_n(h ) \leq \delta,$$
with 
\begin{equation} \label{defLn}\mathcal{L}_n(h) = \sup_{X_n + rv(n)^{1/\alpha} B_{\ell^{\alpha} }}d\big( f_n(x+h) , f_n(x)\big),\end{equation}
satisfying, 
$$(\log n)^{\alpha/2} = o(\log \frac{t_{\delta}(n)^2}{v(n)}) \text{ if }\alpha \neq 1, \text{ or 
 } v(n) = o(t_{\delta}(n)^2)\text{ if }\alpha=1.$$
(iii).(Compactness). For any $r>0$, $\cup_{n\in N} F_n( rB_{\ell^{\alpha}})$ is relatively compact.\\
(iv).(Upper bound = lower bound). For any $x\in \mathcal{X}$,
\begin{equation} \label{condtaux} I_{\alpha }(x)  = \sup_{\delta >0} \limsup_{\underset{n\in N}{n\to +\infty}} I_{n,\delta}(x).\end{equation}
Then $(f_n(X_n))_{n\in N}$ satisfies a LDP with speed $v(n)$ and good rate function $I_{\alpha}$.
\end{The}

%

Let us make some remarks on the assumptions of this theorem.

\begin{Rems}\label{remtheogene}
(a). We will prove that under the assumption that for any sequence $h_n\in \RR^n$, $n\in N$, such that $\sup_n || h_n ||_{\ell^{\alpha}} < +\infty$,
\begin{equation}\label{equidet} d(f_n(X_n+v(n)^{1/\alpha} h_n), F_n(h_n)) \underset{ \underset{n\in N}{n\to+\infty}}{\longrightarrow} 0,\end{equation}
in probability, the lower bound of the LDP holds with the rate function \eqref{condtaux}. 

(b). The assumption $(i)$ that the approximation \eqref{equidet} holds uniformly in $h_n\in rB_{\ell^{\alpha}}$ is crucial for deriving the upper bound of the LDP with rate function \eqref{deftaux}, and is one of the most constraining assumptions of Theorem \ref{theoremgene}. In the applications we develop when $\alpha <2$, this is proven by some concentration inequality and chaining arguments, which can be carried out successfully due to the ``sparsity'' of the ball $B_{\ell^{\alpha}}$.

(c).\label{remlipconst} The formulation of assumption $(ii)$ on the Lipschitz constant of $f_n$ is specially designed to include polynomial functionals $f_n$, as the trace of a polynomial of random matrices. In other words, it says that the ``local'' Lipschitz constant of $f_n$, is small enough  uniformly on the set $X_n + rv(n)^{1/\alpha} B_{\ell^{\alpha}}$. Note that when $f_n$ is $L_2(n)$-Lipschitz with respect to $||\ ||_{\ell^2}$, a sufficient condition for assumption $(ii)$ to be fulfilled is \begin{equation} \label{condlipconst} (\log n)^{\alpha/2} = o\big( \log \frac{1}{L_2(n)^2v(n)}\big) \text{ if } \alpha \in (1,2), \text{ and } v(n) = o\big( \frac{1}{L_2(n)^2} \big)  \text{ if } \alpha =1.\end{equation}
This assumption ensures that the deviations of $f_n(X_n)$ are explained by a heavy-tail phenomenon. For example, it fails to hold for empirical means under $\nu_{\alpha}^n$ when $\alpha \in [1,2)$.

(d). The compactness assumption of $(iii)$ is made to ensure that $I_{\alpha}$ is a good rate function. As one can observe in the proof, without it, the upper bound of the LDP holds only for compact sets.

(e). \label{remlsi} The rate function $I_{\alpha}$ can be simplified in certain cases. Define the function $\tilde{I}_{\alpha}$ by,
$$\forall x \in \mathcal{X}, \ \tilde{I}_{\alpha}(x) = \inf_n\{||h||_{\ell^{\alpha}}^{\alpha} : x = F_n(h), \ h \in \RR^n \}.$$
One can see that,
$$I_{\alpha} = \sup_{\delta>0} \inf_{B(x,\delta)} \tilde{I}_{\alpha}.$$
Thus if $\tilde{I}_{\alpha}$ is lower semi-continuous, then $I_{\alpha} = \tilde{I}_{\alpha}$.
\end{Rems}
The proof is in line with the ideas and the framework developed by Borell and Ledoux in \cite{Borell}, \cite{Borell2} and \cite{LedouxWiener}, \cite{Ledouxflour}, for the large deviations for Wiener chaoses. To make a parallel with their approach, one can observe that the first step in their proof is to show some deterministic equivalent for the Wiener chaoses when deformed in a direction of the reproducing kernel, that is, by \cite[chapter 5 (5.7)]{Ledouxflour}, for any $h\in \mathcal{H}$,
 \begin{equation} \label{equivWiener}||t^{-d}\Psi(x+th) -\Psi^{(d)}(h)||\underset{t \to +\infty}{ \longrightarrow} 0,\end{equation}
in probability, and even uniformly in $h \in \mathcal{O}$ the unit ball of $\mathcal{H}$, along a discretization of $\Psi$ by \cite[chapter 5 (5.9)]{Ledouxflour}, where $\Psi^{(d)}$ defined in \eqref{detequivWiener}. Similarly, we make the assumption $(i)$ that a uniform deterministic equivalent holds for the functionals $f_n$. 

 For the lower bound, we replace the use of the Cameron-Martin formula, used in the context of abstract Wiener space, with a lower bound estimate of the probability of translated events, that is,
\begin{equation} \label{lowerboundnu} \liminf_{\underset{n \in  N}{n \to +\infty}} \frac{1}{v(n)} \log \nu_{\alpha}^n( E + v(n)^{1/\alpha} h_n)\geq - \limsup_{\underset{n\in N}{n \to +\infty}} c_{\alpha}(h_n),\end{equation}
for a given sequence $h_n\in \RR^n$, subsets $E$ such that $\liminf_n \nu_{\alpha}^n(E)>0$, and where $c_{\alpha}$ is some weight function. In the Gaussian case $\alpha=2$, the translation formula of the Gaussian measure gives this estimate with $c_{\alpha}(h) = || h ||_{\ell^2}^2$. When $\alpha <2$, one can mimic the Gaussian case to get such an estimate \eqref{lowerboundnu} with $c_{\alpha}(h) = ||h||_{\ell^{\alpha}}^{\alpha}$, whereas when $\alpha>2$, we believe that there is a competition between the speed and the dimension which is not workable in the applications. 

Whereas the Gaussian isoperimetric inequality is used in the proof of the upper bound of the deviations of Wiener chaoses, ours will rely on sharp large deviation inequalities for $\nu_{\alpha}^n$ with respect to the weight function $c_{\alpha}$, that is
\begin{equation}\label{upperbound} \limsup_{\underset{n \in N}{n \to +\infty}} \frac{1}{v(n)} \log \nu_{\alpha}^n( x \notin E + \{ c_{\alpha} \leq r v(n) \}) \leq -r,\end{equation}
for some ``large enough'' subsets $E$. We will show that we can take $c_{\alpha} = ||h||_{\ell^{\alpha}}^{\alpha}$, which together with  \eqref{lowerboundnu} will allow us to make the upper and lower bound match. In the Gaussian case, this is due to the Gaussian isoperimetric inequality, whereas when $\alpha <2$, we will have to call for sharp inf-convolution inequalities for $\nu_{\alpha}^n$. This is in particular where assumption $(ii)$ plays its role since it enables us, when $\alpha <2$, to neglect the Euclidean balls which come naturally in the deviation inequality of $\nu_{\alpha}^n$, and consider subsets $E$ which are indeed large enough.

These two estimates \eqref{lowerboundnu} and \eqref{upperbound} are behind the limitation in Theorem \ref{theoremgene} to the probability measures $\nu_{\alpha}^n$ for $\alpha \in (0,2]$. For example, if one replaces the measure $\nu_{\alpha}$ by the probability measure on $\RR_+$ with density $ Z_{\alpha}^{-1} e^{-x^{\alpha}}$, one can show that \eqref{lowerboundnu} holds provided $h_n$ has all its coordinates non-negative (and $ n = o(v(n))$ if $\alpha>1$). But then, we will have to prove \eqref{upperbound} with $c_{\alpha}(h) = || h ||_{\ell^{\alpha}}^{\alpha}$ if the coordinates of $h$ are non-negative, and $+\infty$ otherwise, which we do not know how to obtain for the subsets $E$ we are dealing with in the proof.

This said, we can give a version of Theorem \ref{theoremgene} for the probability measure $\mu_{\alpha}$, with density $Z_{\alpha}^{-1}e^{-x^{\alpha}} \Car_{x\geq 0}$, which will be sufficient to prove a LDP result for the last-passage time.

\begin{The}\label{theoremgenesup}
 Let $\alpha \in (0,1]$ and $N\subset \NN$ an infinite subset. Let $X_n$ be a random variable distributed according to $\mu_{\alpha}^{ n}$. Let $f_n, F_n : \RR^n \to \RR$ be measurable functions. Let $(v_n)_{n\in N}$ be a sequence going to $+\infty$. Define $I_{\alpha}$ as in \eqref{deftaux}, and for $\delta>0$ and $n\in N$,
$$I_{n,\delta}^+(x) = \inf\big \{ ||h||_{\ell^{\alpha}}^{\alpha} : d(F_n(h),x)<\delta,  h\in \RR_+^n\big \}.$$
 Assume $(i)-(ii)-(iii)$ from Theorem \ref{theoremgene}, and,\\
(iv)'. For any $x\in \mathcal{X}$,
$$I_{\alpha}(x) = \sup_{\delta>0} \limsup_{\underset{n\in N}{n\to+\infty}}I_{n,\delta}^+(x) .$$   
Then $(f_n(X_n))_{n\in N}$ satisfies a LDP with speed $v(n)$ and good rate function $I_{\alpha}$.
\end{The}

\begin{Rem}
We only state this result for $\alpha \in (0,1]$ because for $\alpha >1$, we know how to get the lower bound \eqref{lowerboundnu} for a sequence $h_n\in \RR_+^n$ only under the additional assumption on the speed that $n= o(v(n))$. But this condition and the requirement $(ii)$ cannot be met simultaneously in the applications we will present.

\end{Rem}

\subsection{Applications to Wigner matrices}

\label{appliWignerIntro}
We present now the applications of Theorem \ref{theoremgene} to Wigner matrices. We denote by $\mathcal{H}_n^{(\beta)}$ the set of Hermitian matrices when $\beta=2 $, and symmetric matrices when $\beta=1$, of size $n$.
We define $\mathcal{S}_{\alpha}$ the class of Wigner matrices whose law is of density $Z_{W_{\alpha}}^{-1}e^{-W_{\alpha}}$ with respect to the Lebesgue measure $\ell_n^{(\beta)}$ on $\mathcal{H}_n^{(\beta)}$, where
\begin{equation} \label{defW} \forall A \in \mathcal{H}_n^{(\beta)}, \ W_{\alpha}(A) = b\sum_{i}|A_{i,i}|^{\alpha} + \sum_{i<j} \Big(a_1 |\Re A_{i,j}|^{\alpha} + a_2 |\Im A_{i,j}|^{\alpha}\Big),\end{equation}
for some $b,a_1, a_2 \in (0,+\infty)$, and where $Z_{W_{\alpha}}$ is the normalizing constant. 

We will denote by $\mu_A$ the empirical spectral measure of  a matrix $A\in \mathcal{H}_n^{(\beta)}$, that is,
$$ \mu_A = \frac{1}{n} \sum_{i=1}^n \delta_{\lambda_i},$$
where $\lambda_1,...,\lambda_n$ are the eigenvalues of $A$, and we will denote by $\lambda_A$ the largest eigenvalue of $A$.

We will say that $X$ is a Wigner matrix if $X$ is a random Hermitian matrix with independent coefficients (up to the symmetry) such that $(X_{i,i})_{1\leq i \leq n}$ are identically distributed and $(X_{i,j})_{i<j}$ are identically distributed. If $\EE|X_{1,2}- \EE X_{1,2}|^2 = 1$, then by Wigner's theorem (see \cite[Theorem 2.1.1, Exercice 2.1.16]{Guionnet}, \cite[Theorem 2.5]{Silverstein}), almost surely,
$$\mu_{X/\sqrt{n}} \underset{n\to\infty}{\leadsto} \mu_{sc},$$
where $\leadsto$ denotes the weak convergence, and $\mu_{sc}$ is the semi-circular law defined by,
$$ \mu_{sc} = \frac{1}{2\pi} \sqrt{4-x^2} \Car_{|x| \leq 2} dx.$$
If we assume furthermore that $\EE {X_{1,1}}^2< +\infty$ and $\EE X_{1,2}^4<+\infty$, then we know by \cite{Bai}, \cite[Theorem 5.1]{Silverstein},
$$\lambda_{X/\sqrt{n}} \underset{n\to +\infty}{\longrightarrow} 2,$$
in probability.

As a consequence of Theorem \ref{theoremgene}, we have the following large deviations principles, originally proven in  \cite{Bordenave}, in the case of the empirical spectral measure and in \cite{LDPei} for the largest eigenvalue.
\begin{The}\label{LDPmsp} Let $\alpha \in (0,2)$. Assume $X $ is in the class $\mathcal{S}_{\alpha}$ such that $\EE|X_{1,2}|^2 =1$.
$(\mu_{X/\sqrt{n}})_{n\in \NN}$ follows a LDP with respect to the weak topology with speed $n^{1+\alpha/2}$, and good rate function $I_{\alpha}$, defined for any probability measure $\mu$ on $\RR$ by,
$$ I_{\alpha}(\mu) =  \sup_{\delta > 0} \inf_{n\in\NN} \{ W_{\alpha}(A)  : A \in \mathcal{H}_n^{(\beta)}, \ d(\mu ,\mu_{sc} \boxplus \mu_{n^{1/\alpha} A} )<\delta\},$$
where $d$ is a distance compatible with the weak topology, $\boxplus$ stands for the free convolution (see \cite[section 2.3.3]{Guionnet} for a definition), and $\mu_{sc}$ is the semi-circular law.

 \end{The}
\begin{Rem}
In \cite{Bordenave}, the rate function $I_{\alpha}$ is computed explicitly for measures $\mu_{sc} \boxplus \nu$, where $\nu$ is a symmetric probability measure, for which we have 
$$I_{\alpha}(\mu_{sc}\boxplus \nu) = \min\big(b,\frac{a}{2} \big) \int |x|^{\alpha} d\nu(x).$$

\end{Rem}
\begin{The}\label{LDPvp}Let $\alpha \in (0,2)$. Assume that $X$ is in the class $\mathcal{S}_{\alpha}$ such that $\EE|X_{1,2}|^2 =1$.
$(\lambda_{X/\sqrt{n}})_{n\in \NN}$ follows a LDP with speed $n^{\alpha/2}$ and good rate function $J_{\alpha}$, defined for any $x\in \RR$ by,
$$ J_{\alpha}(x) = \begin{cases}
cg_{\mu_{sc}}(x)^{-\alpha}& \text{ if } x>2,\\
0 & \text{ if } x=2,\\
 +\infty& \text{ if } x<2,
\end{cases}$$ 
with 
$$c = \inf\big \{ W_{\alpha}(A)  : A \in \cup_{n \in \NN} \mathcal{H}_n^{(\beta)}, \ \lambda_A =1 \big\},$$
and where $g_{\mu_{sc}}$ denotes the Stieltjes transform of $\mu_{sc}$, that is,
$$ \forall z \in \CC\setminus (-2,2), \ g_{\mu_{sc}}(z) = \int \frac{d\mu_{sc}(x)}{z-x}.$$
 \end{The}

\begin{Rem}
The constant $c$ can be computed explicitly, we refer the reader to \cite[section 8]{LDPei} for more details.
\end{Rem}

If $\textbf{X}=(X_1,...,X_p)$ is a collection of independent centered Wigner matrices such that $\EE M_{1,2}^2=1$ for any $M \in \{X_1,...,X_n\}$, and with entries having finite moments of order $d$, then for any non-commutative polynomial $P \in \CC\langle \textbf{X} \rangle$ of total degree $d$, we know by \cite[Theorem 5.4.2]{Guionnet}, 
$$ \tau_n [ P(\textbf{X}/\sqrt{n})] \underset{n\to +\infty}{\longrightarrow} \tau[P(\textbf{s})],$$
in probability, where $\tau_n = \frac{1}{n} \tr$ and $\textbf{s} =(s_1,...,s_p)$ is a free family of $p$ semi-circular variables in a non-commutative probability space $(\mathcal{A},\tau)$ (see \cite[section 5.3]{Guionnet} for a definition). 

Concerning the large deviations of such normalized traces of polynomials in independent matrices in the class $\mathcal{S}_{\alpha}$, with $\alpha \in (0,2]$ we have the following result. 
\begin{The}\label{LDPpoly}
Let $\alpha \in (0,2]$ and $p,d\in \NN$, $d> \alpha$. Assume $\textbf{X} = (X_1,...,X_p)$ is a collection of independent Wigner matrices in the class $\mathcal{S}_{\alpha}$, such that for $M \in \{ X_1,...,X_p\}$, $\EE|M_{1,2}|^2 =1$. We assume that $X_i$ is distributed according to $Z_{W_{\alpha}}^{-1}e^{-W_{\alpha,i}} d\ell_n^{(\beta)}$, where $W_{\alpha,i}$ is of the form \eqref{defW}. Let $P \in \CC\langle \textbf{X} \rangle$ be a non-commutative polynomial of total degree $d$. We denote by $\tau_n$ the state $\frac{1}{n} \tr$ on $\mathcal{H}_n^{(\beta)}$. The sequence 
$$\tau_n [P(\textbf{X}/\sqrt{n})]$$
 satisfies a LDP with speed $n^{\alpha\big( \frac{1}{2}+\frac{1}{d} \big)}$ and good rate function $K_{\alpha}$, defined for all $x\in \RR$ by
$$ K_{\alpha}(x) = \begin{cases}
c_{1}\big(x-\tau(P(\textbf{s}))\big)^{\frac{\alpha}{d}} & \text{ if } x > \tau(P(\textbf{s})),\\
0 & \text{ if } x= \tau(P(\textbf{s})),\\
c_{-1}\big |x-\tau(P(\textbf{s})) \big|^{\frac{\alpha}{d}} & \text{ if } x < \tau(P(\textbf{s})),
\end{cases}$$
where for any $\sigma \in \{-1,1\}$,
$$c_{\sigma}=\inf\big\{ W_{\alpha}( \textbf{H}):  \textbf{H} \in \cup_{n\in \NN} (\mathcal{H}_n^{(\beta)})^p, \sigma = \tr P_d(\textbf{H}) \big\}\in [0,+\infty],$$
where $W_{\alpha}(\textbf{H}) = \sum_{i=1}^p W_{\alpha,i}(H_i)$ and $P_d$ is the homogeneous part of degree $d$ of $P$.
\end{The}

\begin{Rem}Unlike the previous results on deviations of the spectral measure and the largest eigenvalue, this one allows us to consider Gaussian matrices. As we will see in the proof, the mechanism of deviations of traces of polynomials is the same in both cases $\alpha \in (0,2)$, and $\alpha=2$. This is essentially due to the fact that still in the Gaussian case there is a heavy-tail phenomena which appears when the degree of the polynomial is strictly greater than $2$ since there is no exponential moments.

This large deviations principle is an extension, although in a more restricted setting, of the large deviations principle proven in \cite{LDPtr}, in the case where $p=1$ and $P = X^d$ for some $d\geq 3$, for Gaussian matrices and Wigner matrices without Gaussian tails. 
\end{Rem}

\subsection{Application to last-passage percolation}\label{LPP}
Let $d\in \NN$, $d\geq 2$. We denote by $\ZZ_+^d$ the subset of vectors of $\ZZ^d$ with non-negative coordinates. Let $(X_v)_{v\in \ZZ_+^d}$ be a collection of weights. We will call a \textit{directed path} a path in which at each step, one coordinate is increased by $1$. For $v_1,v_2 \in \ZZ^d_+$, we denote by $\Pi(v_1,v_2)$ the set of directed paths from $v_1$ to $v_2$. We will  identify a path with the set of its vertices. We define the \textit{last-passage time} $T_{v_1,v_2}(X)$, by
$$ T_{v_1,v_2}(X) = \sup_{\pi \in \Pi(v_1,v_2)} \sum_{v \in \pi} X_{v},$$
We know by a work of Martin \cite{Martin}, that if the weights $X_v$ are i.i.d random variables with common distribution function $F$ satisfying,
\begin{equation} \label{condfuncdistr} \int_0^{+\infty}(1-F(t))^{1/d} dt <+\infty,\end{equation}
then for any $v\in \RR^d_+$, 
\begin{equation} \label{defg} \frac{1}{n}\EE T_{0, \lfloor nv \rfloor }(X) \underset{n\to+\infty}{\longrightarrow} g(v),\end{equation}
where $g$ is a continuous function on $\RR^d_+$.

As an application of Theorem \ref{theoremgenesup}, we will get the following LDP for the last-passage time.

\begin{The}\label{LDPLPP}Let $\alpha \in (0,1)$. For any $n\in \NN$, we set $T(X) = T_{0,(n,...,n)}(X)$.
Let $(X_v)_{v\in \ZZ^d_+}$ be a family of i.i.d random variables distributed according to $\mu_{\alpha}$.
 The sequence $T(X)/n$ satisfies a LDP with speed $n^{\alpha}$ and good rate function $L_{\alpha}$, defined by 
$$L_{\alpha}(x) = \begin{cases} 
(x-g(1,...,1))^{\alpha} & \text{ if } x\geq g(1,...,1),\\
+\infty  & \text{ otherwise}.
\end{cases}$$

\end{The}

\subsection{Concentration inequalities}\label{concineintro}
In order to prove that assumption $(i)$ holds in the context of Wigner matrices in the class $\mathcal{S}_{\alpha}$ when $\alpha \in (0,2)$ for the largest eigenvalue and the empirical spectral measure, we will prove some concentration inequalities for Wigner matrices which we would like to present as they can be of independent interest. 

To derive such concentration inequalities for functions of the spectrum of random matrices, we will follow the classical argument which consists in considering our functionals as functions of the entries, and taking advantage of the concentration property of the law of the underlying random matrix. This approach is made possible in the setting where the spectrum is a smooth function of the entries, which will be our case as we will work with Hermitian matrices. 

For Wigner matrices with bounded entries, or satisfying a Log-Sobolev inequality, or also for certain unitarily or orthogonally invariant models, concentration inequalities for Lipschitz (convex) linear statistics of the eigenvalues and for the largest eigenvalue, have been extensively studied by Guionnet-Zeitouni \cite{GZconc}, Guionnet \cite[Part II]{GuionnetFlour},  and Ledoux \cite[Chapter 8 \S 8.5]{Ledouxmono} (see also \cite[sections 2.3, 4.4]{Guionnet}). 

More precisely, we will provide concentration inequalities for the linear statistics, the spectral measure and the largest eigenvalue of random Hermitian matrices satisfying a certain concentration property which will be indexed by some $\alpha \in (0,2]$.
As we will see, this concentration property will capture the gradation of speeds of large deviations for the spectral functionals we are interested in, as it has been observed in Theorems \ref{LDPmsp} and \ref{LDPvp}. 

We now present the concentration property with which we will be working. 

\begin{Def}\label{defCalpha}
Let $\alpha\in (0,2]$. We will say in the following that a Wigner matrix $X$ satisfies the \textit{concentration property $\mathcal{C}_{\alpha}$}, if there is a constant $\kappa>0$, such that for any Borel subset $A$ of $\mathcal{H}_n^{(\beta)}$, such that $\PP(X\in A )\geq 1/2$, and any $t>0$, 
\begin{equation} \label{concprop}\PP(X \notin A + \kappa \sqrt{t}B_{\ell^2} + \kappa t^{1/\alpha} B_{\ell^{\alpha}} ) \leq 2e^{-t},\end{equation}
if $\alpha\in [1,2]$, and 
\begin{equation} \label{calpha0}\PP\big( X\notin A + \kappa (\log n)^{\frac{1}{\alpha}-1}\big( \sqrt{r} B_{\ell^2} + r B_{\ell^1} \big) +\kappa  r^{\frac{1}{\alpha}} B_{\ell^{\alpha}} \Big) \leq 4 e^{-r},\end{equation}
if $\alpha \in (0,1)$,
where for any $p>0$, 
$$B_{\ell^p}= \big\{ Y \in \mathcal{H}_n^{(\beta)} : ||Y||_{\ell^p} \leq 1 \big\},$$
with
$$ \forall Y \in \mathcal{H}_n^{(\beta)}, \ ||Y||_{\ell^p}^p = \sum_{i,j} |Y_{i,j}|^p.$$
\end{Def}

When $\alpha \in [1,2]$, the motivation for defining this concentration property $\mathcal{C}_{\alpha}$ comes from Talagrand's famous two-levels deviation inequality \cite{Talagrand} for the measure $\nu_{\alpha}^n$, which says that there is a constant $L>0$ such that for any $n\in \NN$, any Borel subset $A$ of $\RR^n$ with $\nu_{\alpha}^n(A)>0$, and $r>0$,
\begin{equation} \label{Taldev}  \nu_{\alpha}^n(x \notin A+ \sqrt{r} B_{\ell^2} + r^{\frac{1}{\alpha}} B_{\ell^{\alpha}}) \leq \frac{e^{-Lr}}{\nu_{\alpha}^n(A)},\end{equation}
and similarly for $\mu_{\alpha}$.

In particular, the Wigner matrices in the  class $\mathcal{S}_{\alpha}$ for $\alpha \in [1,2]$ satisfy the concentration property $\mathcal{C}_{\alpha}$ with some $\kappa$ depending on the parameters $b,a_1,a_2$ of the law of $X$ (see \eqref{defW}). More generally, we know by the results of Bobkov-Ledoux \cite[Corollary 3.2]{BobLed}, and Gozlan \cite[Proposition 1.2]{Gozlan} that if $X$ is a Wigner matrix with entries satisfying a certain Poincaré-type inequality, where the underlying metric on $\RR^m$, $m=1,2$, is the following,
\begin{equation} \label{defdalpha}\forall x,y \in \RR^m,  \ d_{\omega_{\alpha}}(x,y) = \Big( \sum_{i=1}^m |\omega_{\alpha}(x_i)-\omega_{\alpha}(y_i)|^2\Big)^{1/2},\end{equation}
where $\omega_{\alpha}(t) = \sg(t) \max( |t|,|t|^{\alpha})$, $\sg(t)$ standing for the sign of $t$, then $X$ satisfies the concentration property $\mathcal{C}_{\alpha}$ with some constant $\kappa$ depending on the spectral gap. We will get into more details in section \ref{Chapconc} about this functional inequality, and  present some workable criterion available for a Wigner matrix to satisfy $\mathcal{C}_{\alpha}$ when $\alpha \in [1,2]$.


When $\alpha \in (0,1)$, the concentration property of the law of Wigner matrices in the class $\mathcal{S}_{\alpha}$ differs significantly from the case where $\alpha \in [1,2]$. 
We know by Talagrand \cite[Proposition 5.1]{Talexpo} that as $\nu_{\alpha}$ does not have exponential tails, $\nu_{\alpha}^n$ cannot satisfy a dimension-free concentration inequality. Transporting $\nu_1^n$ onto $\nu_{\alpha}^n$, we will prove the following deviation inequality.

\begin{Pro}\label{devnualpha0}
Let $n\in \NN$, $n\geq 2$. There is a constant $c>0$ depending on $\alpha$, such that for any $r>0$, $A$ Borel subset of $\RR^n$, and $C>0$ such that $\nu_{\alpha}^{n}(A) > 1/C$,
$$\nu_{\alpha}^{n} \Big( x \notin A + C(\log n)^{\frac{1}{\alpha}-1}\big( \sqrt{r} B_{\ell^2} + r B_{\ell^1} \big) + r^{\frac{1}{\alpha}} B_{\ell^{\alpha}} \Big) \leq \frac{e^{-c r}}{\nu_{\alpha}^{n}(A) - 1/C}.$$
\end{Pro}
We will discuss in remark \ref{remdev0} in section \ref{sectiondevineq0} the optimality of such a deviation inequality for $\nu_{\alpha}$. 
The above proposition justifies the definition of the concentration property $\mathcal{C}_{\alpha}$ in the case where $\alpha \in (0,1)$, as it implies that Wigner matrices in the class $\mathcal{S}_{\alpha}$ satisfy this property when $\alpha \in (0,1)$.

Regarding the linear statistics of Wigner matrices having concentration $\mathcal{C}_{\alpha}$, we will consider different families of function whether $\alpha \in (0,1)$ or $\alpha \in [1,2]$. To this end, we define $\mathcal{M}_s^{\alpha}$ the set of finite signed measures $\sigma$ such that its total variation $|\sigma|$ has a finite $\alpha^{\text{th}}$-moment. Following \cite[Chapter 2 \S 5.1]{Samko}, we define when $\alpha \in (0,1)$, the fractional integrals of order $\alpha+1$ of $\sigma \in \mathcal{M}_s^{\alpha}$, by
\begin{align}
\forall t \in \RR,\ & (\mathcal{I}^{\alpha+1}_+ \sigma) (t) =\frac{1}{\Gamma(\alpha+1)} \int_{-\infty}^t (t-x)^{\alpha} d\sigma(x),\nonumber\\
& (\mathcal{I}^{\alpha+1}_- \sigma)(t) = \frac{1}{\Gamma(\alpha+1)}\int_t^{+\infty} (x-t)^{\alpha}d\sigma(x).
\label{defintfrac}\end{align}
This definition interpolates for non-integer order the usual iterated integral (see \cite[Chapter 1 \S 2.3]{Samko} for more details). 
With these definitions, we will prove the following deviations inequalities.

\begin{Pro}
\label{conclinearstat}Let $\alpha \in (0,2]$.
Let $X$ be a Wigner matrix having concentration $\mathcal{C}_{\alpha}$ with some $\kappa>0$.   There is a constant $c_{\alpha}>0$ such that if $\alpha \in [1,2]$ and $f : \RR\to \RR$ is some $1$-Lipschitz function, then for any $t>0$,
$$ \PP\big( \mu_{X/\sqrt{n}}(f) - m_f >t \big) \leq 2 \exp\Big( -c_{\alpha} \min\Big(\frac{n^2t^2}{\kappa^2}, \frac{n^{1+\frac{\alpha}{2}}t^{\alpha}}{\kappa^{\alpha}} \Big)\Big),$$
if $\alpha \in (0,1)$, $f$ is $1$-Lipschitz and moreover $f = \mathcal{I}_{\pm}^{\alpha+1}(\sigma)$ for some $\sigma \in \mathcal{M}_s^{\alpha}$ such that $|\sigma|(\RR)\leq m$, then for any $t>0$,
$$ \PP\big( \mu_{X/\sqrt{n}}(f) - m_f >t \big) \leq 4 \exp\Big( -c_{\alpha} \min\Big(\frac{n^2 t^2}{\kappa^2(\log n)^{2(\frac{1}{\alpha}-1)} }, \frac{n^{\frac{3}{2}}t}{\kappa (\log n)^{\frac{1}{\alpha}-1}},\frac{n^{1+\frac{\alpha}{2}} t}{\kappa m}  \Big)\Big),$$
where $m_f$ denotes a median of $\mu_{X/\sqrt{n}}$. 
\end{Pro}

\begin{Rem}
The reason for considering the class of function $\mathcal{I}_{\pm}^{\alpha+1}(\mathcal{M}_s^{\alpha})$ in the case $\alpha\in (0,1)$, comes from the fact that we only understand the stability of the empirical spectral measure with respect to $|| \ ||_{\ell^{\alpha}}$, by using a certain distance $d_{\alpha}$ which controls this class of functions (see section \ref{spvarsection} for more details). 

Still in the case $\alpha<1$, note that one cannot expect the above concentration inequality to be true for all Lipschitz functions, since a change of large deviations speed may occur as the entries of $X$ do not have exponential tails. Indeed, for example if $X$ is in the class $\mathcal{S}_{\alpha}$, Theorem \ref{LDPpoly} tells us the speed of large deviations of $\frac{1}{n}\tr (X/\sqrt{n})$ is $n^{3\alpha/2}$.

\end{Rem}

\begin{Rem}
One can identify the image $\mathcal{I}_{\pm}^{\alpha +1}(\mathcal{M}_s^{\alpha})$, by a minor change of \cite[Theorem 6.3]{Samko}. To ease the notation, we will only describe $\mathcal{I}_{+}^{\alpha +1}(\mathcal{M}_s^{\alpha})$.
For any $\phi \in L^1(\RR)$, one can define the fractional integral of order $\alpha$ by,
$$\forall x \in \RR, \   \mathcal{I}_{+}^{1-\alpha}(\phi)(x) = \frac{1}{\Gamma(1-\alpha)} \int_{0}^{+\infty} t^{-\alpha} \phi(x- t) dt.$$
The function above is well-defined almost everywhere as $t^{-\alpha}\phi(x- t)$ is integrable on a neighborhood of $0$ for almost all $x$ by Fubini theorem. With this definition, the set $\mathcal{I}_{+}^{\alpha +1}(\mathcal{M}_s^{\alpha})$ consists of the functions $f$ such that there is some $\phi \in L^1(\RR)$ and $\sigma \in \mathcal{M}_s^{\alpha}$, such that
$$ \forall x \in \RR, \ f(x) = \int_{-\infty}^x \phi(t) dt,\text{ and } \mathcal{I}_{+}^{1-\alpha}(\phi)(x) = \sigma(-\infty,x].$$
\end{Rem}
\begin{Rem}\label{concsimpl0}

Note also that the exponential bound can be simplified in the case $\alpha \in (0,1)$ if $m\geq c_0$, where $c_0$ is a constant independent of $n$. One gets then, for any $t>0$,
$$ \PP\big( \mu_{X/\sqrt{n}}(f) - m_f >t \big) \leq 4 \exp\Big( -c_{\alpha} \min\Big(\frac{n^2 t^2}{\kappa^2(\log n)^{2(\frac{1}{\alpha}-1)} }, \frac{n^{1+\frac{\alpha}{2}} t}{\kappa m}  \Big)\Big).$$
\end{Rem}

In order to state our concentration inequality for the spectral measure, we will work with the following distance $d$ defined on the set of probability measures on $\RR$, denoted by $\mathcal{P}(\RR)$, in order to quantify the deviations:
\begin{equation} \label{defdStiel} \forall \mu, \nu \in \mathcal{P}(\RR), \ d(\mu, \nu) = \sup_{z\in \mathcal{K}} |g_{\mu}(z) - g_{\nu}(z) |,\end{equation}
where $\mathcal{K}$ is a compact subset of $\{ z \in \CC : \Im z \geq 2\}$ with an accumulation point, such that $\mathrm{diam}(\mathcal{K}) \leq 1$, and with $g_{\mu}$ the Stieltjes transform of $\mu$, that is,
$$\forall z \in \CC^+, \ g_{\mu}(z) = \int \frac{d\mu(t)}{z-t},$$
where $\CC^+ = \{ z \in \CC : \Im z >0\}$. This distance metrizes the weak topology on $\mathcal{P}(\RR)$  by \cite[Theorem 2.4.4]{Guionnet}.

We will prove the following concentration inequalities for the empirical spectral measure and the largest eigenvalues of Wigner matrices having concentration $\mathcal{C}_{\alpha}$.

\begin{Pro} \label{concspintro}
Let $\alpha \in (0,2]$. Let $X$ be a Wigner matrix satisfying $\mathcal{C}_{\alpha}$ with some $\kappa>0$. There exists a constant $c_{\alpha}>0$, depending on $\alpha$, such that for any $t>0$,
$$ \PP\left( d\big( \mu_{X/\sqrt{n}}, \EE \mu_{X/\sqrt{n}} \big) > t +\delta_n\right) \leq \frac{32}{t^2}\exp\big ( -c_{\alpha}k_{\alpha}(t) \big),$$
where $\delta_n = O\big(\kappa n^{-1} ( \log n)^{(1/\alpha-1)_+}\big)$, and where for $\alpha \in [1,2]$,
$$ k_{\alpha}(t) = \min\Big(\frac{n^2t^2}{\kappa^2}, \frac{n^{1+\frac{\alpha}{2}}t^{\alpha}}{\kappa^{\alpha}} \Big),$$
whereas for $\alpha \in (0,1)$
$$ k_{\alpha}(t) =  \min\Big(\frac{n^2 t^2}{\kappa^2(\log n)^{2(\frac{1}{\alpha}-1)} }, \frac{n^{1+\frac{\alpha}{2}} t}{\kappa}  \Big).$$
\end{Pro}

%
\begin{Pro}\label{concvpintro}
Let $\alpha \in (0,2]$. Let $X$ be a Wigner matrix satisfying $\mathcal{C}_{\alpha}$ for some $\kappa>0$. There is a constant $c_{\alpha}>0$, such that for any $t>0$,
$$ \PP\left( \big| \lambda_{X/\sqrt{n}}-\EE  \lambda_{X/\sqrt{n}}\big| > t +\eps_n\right) \leq 8\exp\big(-c_{\alpha} h_{\alpha}(t)  \big),$$
where 
\begin{equation}\label{h1vp} h_{\alpha}(t) = \min\Big( \frac{t^2n}{\kappa^2}, \frac{t^{\alpha}n^{\frac{\alpha}{2}}}{\kappa^{\alpha}} \Big),\end{equation}
if $\alpha \in [1,2]$, and 
\begin{equation}\label{h0vp} h_{\alpha}(t) = \min\Big( \frac{t^2n}{\kappa^2(\log n)^{2(\frac{1}{\alpha}-1) }},\frac{t\sqrt{n}}{\kappa (\log n)^{\frac{1}{\alpha}-1} }, \frac{t^{\alpha}n^{\frac{\alpha}{2}}}{\kappa^{\alpha} } \Big),\end{equation}
if $\alpha \in (0,1)$, and where $\eps_n = O( \kappa n^{-1/2} (\log n)^{(1/\alpha-1)_+})$, uniformly in $H \in \mathcal{H}_n^{(\beta)}$.
\end{Pro}

\subsection{Spectral variation inequalities}
We would like also to advertise for  some spectral variation inequalities, which are not particularly new, but which are maybe a little less known in the form we will propose. 
Indeed, to obtain the concentration inequality of Proposition \ref{concspintro}, we need to understand the stability of the spectrum of Hermitian matrices with respect to the distance $|| \ ||_{\ell^p}$ for $p\geq 1$ or $|| \ ||_{\ell^p}^p$ when $p<1$.

For $p\geq 1$, define the $L^p$-Wasserstein distance on the set of probability measures on $\RR$ with finite $p^{\text{th}}$-moment by,
$$\mathcal{W}_p(\mu,\nu) = \Big( \inf_{\pi} \int |x-y|^p  d\pi(x,y) \Big)^{1/p},$$
where the infimum is over all coupling $\pi$ between $\mu$ and $\nu$, two probability measures on $\RR$ with finite $p^{\text{th}}$-moment.

When $p\geq1$, we get as a mere consequence of Lidskii's theorem (see \cite[Theorem III.4.1]{Bhatia}) the following lemma.
\begin{Lem}\label{spvar1intro}
Let $p\in [1,2]$, and $A,B \in \mathcal{H}_n^{(\beta)}$.
$$\mathcal{W}_p(\mu_A,\mu_B) \leq \frac{1}{n^{1/p}} || A-B||_{\ell^p}.$$
As a consequence,
$$d(\mu_A,\mu_B) \leq \frac{1}{n^{1/p}} || A-B||_{\ell^p}.$$
\end{Lem}

Whereas for $p<1$, we obtain by Rofteld's inequality (see \cite[Theorem IV.2.14]{Bhatia} or  \cite{Thompson}) the following.

\begin{Lem}\label{spvar0intro}Let $p\in (0,1)$.
Let $A,B \in \mathcal{H}_n^{(\beta)}$. For any $t\in \RR$,
$$ \big| \sum_{i=1}^n (t-\lambda_i(A))_+^p - \sum_{i=1}^n (t-\lambda_i(B))_+^{p}\big| \leq \sum_{i=1}^n |\lambda_i(A-B)|^p,$$
where $\lambda_1(A),...,\lambda_n(A)$ denote the eigenvalues of $A$, and similarly for $B$.
Furthermore, there is a positive constant $C_p$, such that for any  $A,B \in \mathcal{H}_n^{(\beta)}$,
$$ d(\mu_{A},\mu_B) \leq \frac{C_p}{n} ||A-B||_{\ell^p}^p,$$
with 
$$C_p = \sqrt{\pi}(p+1) \frac{\Gamma\big( \frac{p+1}{2}\big)}{\Gamma\big( 1+\frac{p}{2}\big)}.$$
\end{Lem}


\subsection*{Acknowledgements}
I would like to thank my supervisor Charles Bordenave for his inspiring advice and the many fruitful conversations which helped me build the present paper. I am also grateful to Franck Barthe and Michel Ledoux for precious conversations and references, as well as Guillaume Aubrun for pointing me out the result of \cite[Proposition 3.2.2]{ETentropy}. I would like also to thank IMPA for its welcome, where this work was partially carried out.

\subsection{Organization of the paper}

In the section \ref{infconv}, we prove some inf-convolution inequalities for $\nu_{\alpha}^{ n}$. As the large deviations of our functional $f_n$ are governed by translates, we will need some sharp deviation inequalities with respect to the metric $||\ ||_{\ell^{\alpha}}$ (or $|| \ ||_{\ell^{\alpha}}^{\alpha}$ when  $\alpha<1$). We will provide a family of weights $W_{\alpha,\eps}$ which captures the asymptotics of the tail distribution of $\nu_{\alpha}^{n}$, that is, behaving like $|| x ||_{\ell^{\alpha}}^{\alpha}$ when $||x||_{\infty} \gg 1$. This will be done by transporting and tensoring the family of optimal weights known for the exponential law due to Talagrand \cite[Theorem 1.2]{Talagrand2}. 

In the section \ref{prooftheorem}, we give a proof of Theorems \ref{theoremgene} and \ref{theoremgenesup}. The upper bound relies on Proposition \ref{grandev} which gives a large deviations sharp upper bound for $\nu_{\alpha}^{n}$ with respect to the metric $|| \ ||_{\ell^{\alpha}}$ using the inf-convolution inequalities proven in section \ref{infconv}.  The lower bound is given by Proposition \ref{lowb} which estimates at the exponential scale $v(n)$ the probability, under $\nu_{\alpha}^n$, of an event translated by some element $v(n)^{1/\alpha} h_n$. 

The rest of the paper is devoted to applications to Wigner matrices and the last-passage time.

In the section \ref{Chapconc}, we prove the concentration inequalities of Propositions \ref{concspintro} and \ref{concvpintro} for the largest eigenvalue, linear statistics and empirical spectral measure of Wigner matrices satisfying the concentration property $\mathcal{C}_{\alpha}$ defined in \eqref{concprop} and \eqref{calpha0}.
To do so, we will prove and discuss the spectral variations inequalities in Lemmas \ref{spvar1intro} and \ref{spvar0intro} in section \ref{spvarsection}.

In section \ref{secdetermexpoWigner}, we show some uniform deterministic equivalents for the spectral measure, largest eigenvalue and traces of non-commutative polynomials of deformed Wigner matrices in the class $\mathcal{S}_{\alpha}$. To make the equivalents for the spectral measure and largest eigenvalue of hold uniformly for $\alpha<2$, we make use of the concentration inequalities we proved in section \ref{Chapconc}, and perform a classical chaining argument.

In section \ref{secdetermequivLPT}, we provide a deterministic equivalent for the last-passage time under additive deformations of the weights. The strategy to make our equivalent hold uniformly will be the same as for the case of the spectral measure and largest eigenvalue of Wigner matrices in the class $\mathcal{S}_{\alpha}$, meaning that it will rely on concentration and chaining arguments.

In section \ref{Wigner}, we apply Theorem \ref{theoremgene} in the setting of Wigner matrices in the class $\mathcal{S}_{\alpha}$, to the spectral measure, the largest eigenvalue (for $\alpha \in (0,2)$) and to traces of non-commutative polynomials (for $\alpha \in (0,2]$). Using of the uniform deterministic equivalents we proved in section \ref{secdetermexpoWigner}, we give a proof of Theorems \ref{LDPmsp}, \ref{LDPvp}, and \ref{LDPpoly}. 

Finally we prove in section \ref{LPPLDP},  the large deviations principle for the last-passage time of Theorem \ref{LDPLPP} by applying Theorem \ref{theoremgenesup} and using the uniform deterministic equivalent proved in section  \ref{secdetermequivLPT}.

\section{Inf-convolution inequalities for $\nu_{\alpha}^{ n}$} \label{infconv}
Let $\nu$ be a probability measure on $\RR^n$, and let $w$ be a measurable function on $\RR^n$ taking non-negative values. Following Maurey (see \cite{Maurey}), we will say that $(\nu, w)$ satisfies the $\tau$-property if for any non-negative measurable function $f$ on $\RR^n$,
\begin{equation}\label{IC} \big(\int e^{f\Box w } d\nu \big) \big(\int e^{-f} d\nu \big) \leq 1,\end{equation}
where $ \Box $  denotes the inf-convolution, that is,
$$ \forall x \in \RR^n, \ f\Box w (x) = \inf_{y \in \RR^n} \{ f(y) + w(x-y)\}.$$ \label{nomen:infconv} \nomenclature[]{$\Box $}{inf-convolution}{}{\pageref{nomen:infconv}}

The $\tau$-property is closely linked to transportation-cost inequalities. By the Kantorovitch duality (see \cite[Theorem 5.10]{Villani}), and the duality of the entropy (see \cite[Lemma 6.2.13]{Zeitouni}), it is known that under mild assumptions on $w$ that the following general inf-convolution inequality,
\begin{equation} \label{infconvgene} \int e^{f\Box w} d\nu \leq e^{\int f d\nu},\end{equation}
satisfied for any non-negative measurable function $f$ is equivalent to the following transportation-cost inequality: for any $\mu$ probability measure on $\RR^n$, 
\begin{equation} \label{transportcost} \mathcal{W}_w(\mu, \nu) \leq D( \mu||\nu),\end{equation}
where $D(\mu || \nu)$ is the relative entropy of $\mu$ with respect to $\nu$, and 
\begin{equation} \label{distKR} \mathcal{W}_w(\mu,\nu) = \inf\big\{  \int w(x-y) d\pi (x,y) :\pi \text{ has marginals } \mu \text{ and } \nu \big\}.\end{equation}
 In particular, under the assumption that $w$ is upper semi-continuous, Kantorovitch duality is valid by  \cite[Theorem 5.10]{Villani}, so that the equivalence above between \eqref{infconvgene} and \eqref{distKR} holds.

One can observe that if $(\nu,w)$ satisfies the $\tau$-property, then  by Jensen's inequality, it satisfies also the general inf-convolution inequality \eqref{infconvgene}, and therefore $\nu$ satisfies the transportation-cost inequality \eqref{transportcost} with cost function $w$. 

Conversely, according to \cite[Proposition 4.13]{Gozlan},  if $\nu$ satisfies the transportation-cost inequality \eqref{transportcost} with cost function $w$, then $(\nu, w\Box w)$ satisfies the $\tau$-property. If moreover $w$ is sub-additive, then one can see that $w\Box w =w$ and thus $(\nu, w)$ satisfies the $\tau$-property. Whereas if $w$ is convex, then $w \Box w =  2w(./2)$ so that $(\nu, 2w(./2))$ satisfies the $\tau$-property. This remark will be useful later when we will need to translate a transportation-cost inequality into a $\tau$-property.

More importantly for us, the $\tau$-property yields deviations bounds with respect to enlargements by the weight $w$. We know from \cite[Lemma 4]{Maurey}, that if $(\nu, w)$ satisfies the $\tau$-property, then for any Borel subset $A$ of $\RR^n$, and any $t>0$,
\begin{equation} \label{dev} \nu\big( x \notin A +\{w \leq r\}\big)\leq \frac{e^{-r}}{\nu(A)}.\end{equation}

We define another form of inf-convolution inequality, designed to enable us to get the best constants in our weight functions, (and also to deal with the measure $\nu_{\alpha}^{n}$ when $\alpha \in (0,1)$), which we will call the \textit{truncated $\tau$-property}. More precisely, we will say that a measure $\nu$ on $\RR^n$ with the weight function $w$, satisfies the \textit{$A_0$-truncated $\tau$-property}, where $A_0$ is a Borel subset of $\RR^n$, if \eqref{IC} is true for any non-negative measurable function $f$ such that $f = +\infty$ on $A_0^c$.

This $A_0$-truncated $\tau$-property yields a deviation inequality with respect to enlargement by the weight $w$ of the following form: for any Borel subset $A$ of $\RR^n$ such that $\nu(A)>0$, and any $r>0$,
 \begin{equation}\label{devtrunc} \nu\big( x \notin A + \{w \leq r\} \big) \leq \frac{e^{-r}}{\nu(A\cap A_0)}.\end{equation}

The goal of this section is to find, for the measure $\nu_{\alpha}^{n}$, when $\alpha \in (0,2)$, a family of weights $W_{\alpha, \eps}$ for which a truncated $\tau$-property is satisfied, and which captures the asymptotics of the tail distribution of $\nu_{\alpha}^{n}$. More precisely, we will prove the following proposition.

\begin{Pro}\label{tau}
Let $\alpha >0$. 
If $\alpha=1$, then for any $\eps<1/2$, $(\nu_1^n, W_{1,\eps})$ satisfies the $\tau$-property with 
$$\forall x \in \RR, \ W_{1,\eps}(x) = \sum_{i=1}^n w_{\eps}(x_i),$$ where
$$ w_{\eps}(t) = \begin{cases}
\frac{\eps e^{-1/\eps} t^2}{8} & \text{ if } |t| \leq 2/\eps^2,\\
(1-2\eps)|t| & \text{ if } |t|> 2/\eps^2.
\end{cases}$$
 If $\alpha \neq 1$, there are some constants $\kappa>0$ and $\eps_0 \in (0,1)$ such that for $\eps \in (0,\eps_0)$ and $m\geq 1$, 
 $(\nu_{\alpha}^{n}, W_{\alpha,\eps}^{(m)})$ satisfies the $mB_{\ell^{\infty}}$-truncated $\tau$-property, where
\begin{equation}\label{weight0}\forall x \in \RR^n, \  W_{\alpha,\eps}^{(m)}(x) = \sum_{i=1}^n w_{\alpha,\eps}^{(m)}(x_i),\end{equation}
with 
$$ w_{\alpha,\eps }^{(m)}(t) = \begin{cases}
\kappa^{-1} e^{-(\frac{m}{\eps})^{\alpha/2}} t^2 & \text{ if } |t| \leq m\eps^{-1},\\
(1- \kappa\eps^{(\alpha/2)\wedge 1} )|t|^{\alpha}  & \text{ if } |t| >  m \eps^{-1}.
\end{cases}$$

\end{Pro}

The rest of this section will be devoted to proving the above proposition. We will reduce the problem in a first phase to the one-dimensional case, and to an estimation of the monotone rearrangement of $\nu_{1}$ onto $\nu_{\alpha}$.

As the usual $\tau$-property (see \cite[Lemma 1]{Maurey}), the truncated version of the $\tau$-property tensorizes in the following way.

\begin{Lem} \label{tensortau}Let $\nu_i$ be a probability measure defined on some measurable space $\mathcal{X}_i$, $A_i$ be some measurable subset of $\mathcal{X}_i$ and $w_i : \mathcal{X}_i \to \RR_+$ be a measurable function, for $i=1,2$. 

If $(\nu_i,w_i)$ satisfies the $A_i$-truncated $\tau$-property for $i=1,2$, then $(\nu_1\otimes \nu_2, w)$ satisfies the $A_1\times A_2$-truncated $\tau$-property with
$$ \forall (x,y) \in \mathcal{X}_1\times \mathcal{X}_2, \ w(x,y)  = w_1(x) + w_2(y).$$

\end{Lem}
Since we are dealing with the product measure $\nu_{\alpha}^n$, we can focus on studying the $\tau$-property for the one-dimensional marginal $\nu_{\alpha}$. 

For the exponential measure, we have the following result due to Talagrand, which gives a family of optimal weights $c_{\lambda}$.
\begin{Pro}[{\cite[Theorem 1.2]{Talagrand}}]\label{Tal}
Let $\lambda \in (0,1)$. Define the weight function $c_{\lambda}$ for any $x\in \RR$ by,
$$ c_{\lambda}(x) =\big( \frac{1}{\lambda} -1 \big) (e^{-\lambda |x|} -1 +\lambda |x|). $$
For any $\lambda \in (0,1)$, $\nu_1$ satisfies a transportation-cost inequality \eqref{transportcost} with cost function $c_{\lambda}$.
\end{Pro}

Note that,  $c_{\lambda}(x)  \sim_{\pm \infty} (1-\lambda)|x|$. Thus, when $\lambda \ll1$, $c_{\lambda}$  captures the exact asymptotics of the tail distribution of the exponential law.

For technical reasons, we prefer to work with a different family of weights than the one defined in Proposition \ref{Tal}. In the following corollary, we reformulate Talagrand's result for the symmetric exponential measure $\nu_1$. 

\begin{Cor}\label{corexp}
Let $\delta >0$. We define the weight function $w_{\eps}$, for any $t \in \RR$, by
$$ w_{\delta}(t) = \begin{cases}
\frac{\delta e^{-1/\delta} t^2}{8} & \text{ if } |t| \leq 2/\delta^2,\\
(1-2\delta)|t| & \text{ if } |t|> 2/\delta^2.
\end{cases}$$
For any $\delta \in (0,1/2)$, $(\nu_1,w_{\delta})$ satisfies the $\tau$-property. As a consequence, $(\nu_1^n, W_{1,\delta})$ satisfies the $\tau$-property, with $W_{1,\delta}$ defined in Proposition \ref{tau}.
\end{Cor}
This reformulation reveals in particular the structure of the enlargements given by the weights $c_{\lambda}$ which consist in a mixture of $\ell^2$ and $\ell^1$-balls.

\begin{proof}As $c_{\lambda}$ is a convex function, we know by \cite[Proposition 4.13]{Gozlan} that $(\nu_1, 2c_{\lambda}(./2))$ satisfies the $\tau$-property.
To prove Corollary \ref{corexp}, it suffices to prove that $w_{\delta} \leq 2c_{\delta}(./2)$ for any $\delta \in (0,1/2)$. Since both functions are even, it is sufficient to prove the inequality on $\RR_+$.  Let $t>0$. By Taylor's formula
$$ e^{-\delta t}-1+\delta t= \delta^2e^{-\delta y} \frac{t^2}{2},$$
for some $y\in [0,t]$. If $t\leq 2/\delta^2$ and $\delta\leq 1/2$, we get
$$2c_{\delta}(t/2)  \geq \delta\big(1-\delta) e^{-1/\delta} \frac{t^2}{4}\geq w_{\delta}(t).$$
If $t\geq 1/\delta^2$, we have
$$c_{\delta}(t) \geq \big( \frac{1}{\delta}-1\big)(-1+\delta t)\geq (1-\delta)t-\frac{1}{\delta}\geq (1-2\delta)t.$$
Thus, $2c_{\delta}(t/2)\geq (1-2\delta)t$ for $t\geq 2/\delta^2$.

After tensorization (see \cite[Lemma 1]{Maurey}), we obtain that $(\nu_1^n, W_{1,\delta})$ satisfies the $\tau$-property with $W_{1,\delta}$ defined in Proposition \ref{tau}.

\end{proof}

For $\alpha \neq1$, the general strategy is to transport this $\tau$-property of the symmetric exponential law to obtain a $\tau$-property for $\nu_{\alpha}$. 
It extends in our setting of truncated $\tau$-property, a result of Maurey \cite[Lemma 2]{Maurey}.
\begin{Lem}\label{transp0}
Let $A$ be a Borel subset of $\RR^n$. Let $\mu$ be a probability measure on $\RR^n$ and let $\psi : \RR^n \to \RR^n$ be a bijective measurable map. Assume $(\mu, w)$ satisfies the $\tau$-property. Let $A$ be a Borel subset of $\RR^n$ and let $\tilde{w}$ be a weight function such that,
$$\forall x \in \RR^n, y \in A,\ \tilde{w}(x-y)\leq  w\big(\psi^{-1}(x) - \psi^{-1}(y)\big).$$
Then, $(\mu \circ \psi^{-1}, \tilde{w})$ satisfies the $A$-truncated $\tau$-property.
\end{Lem}

\begin{proof}
Let $f : \RR^n \to \RR$ be a measurable non-negative function being $+\infty$ on $A^c$. Applying the $\tau$-property of $(\mu,w)$ to $f\circ \psi$, we get
$$\big( \int e^{f\circ \psi \Box w} d\mu \big) \big( \int e^{-f \circ \psi} d\mu \big) \leq 1.$$ 
But, as $\psi$ is a bijection and $f=+\infty$ on $A^c$,
\begin{align*}  f \circ \psi \Box w(\psi^{-1}(x))& = \inf_{y\in \RR^n} \{ f(y) + w( \psi^{-1}(x) - \psi^{-1}(y)) \}\\
& = \inf_{y\in A} \{ f(y) + w( \psi^{-1}(x) - \psi^{-1}(y)) \}.
\end{align*}
From the assumption on $\tilde{w}$, we deduce
$$ f \circ \psi \Box w(\psi^{-1}(x)) \geq  \inf_{y\in A} \{ f(y) + w( x -y) \} = f \Box w(x).$$
Therefore,
$$\big( \int e^{f \Box w} d \mu \circ \psi^{-1} \big) \big( \int e^{-f} d\mu\circ \psi^{-1} \big)\leq 1.$$

\end{proof}

In particular, in the one-dimensional case, if $(\mu, w)$ satisfies the $\tau$-property and $w$ is even and non-decreasing on $\RR_+$, then $\mu \circ \psi^{-1}$ satisfies the $A$-truncated $\tau$-property with any even weight function $\tilde{w}$ such that
$$ \forall s \geq 0, \ \tilde{w}(s) \leq w\big( \Delta_A(s)\big),$$
where  $\Delta_A$ is defined for any $s\geq 0$ by,
\begin{equation*}  \Delta_A(s) =  \inf \big\{ |\psi^{-1}(x) - \psi^{-1}(y)|, |x-y| = s, \ x \in A \big\}.\end{equation*}
If $\mu$ and $\nu$ are two probability measures on $\RR$, we define the \textit{monotone rearrangement} $T$ of $\mu$ onto $\nu$ by,
$$\forall t \in \RR, \ \mu(-\infty,t] =  \nu(-\infty, T(t)].$$
This defines a unique non-decreasing map if the distribution function of $\nu$ is invertible, which sends $\mu$ to $\nu$.

Let $\psi$ be the monotone rearrangement of $\nu_1$ onto $\nu_{\alpha}$.
One can easily check that $\psi$ is an odd function, and that its restriction $\phi$ on $\RR^+$ satisfies,
$$ \forall x \geq 0, \  e^{-x} = \int_{\phi(x)}^{+\infty} e^{-u^{\alpha}} \frac{du}{Z_{\alpha}},$$
where $Z_{\alpha}$ is the normalizing constant of $\mu_{\alpha}$, so that $\phi$ is the monotone rearrangement of $\mu_1$ onto $\mu_{\alpha}$.
Thus, we are reduced to understand the behavior of the map $\phi$ and how it deforms the weights $c_{\eps}$ of Proposition \ref{Tal}.

\subsection{Behavior of the monotone rearrangement}\label{sectionBrenier}
When $\alpha \geq 1$, we have the following estimate on the monotone rearrangement due to Talagrand \cite{Talagrand}.
\begin{Lem}[{\cite[Lemma 2.5]{Talagrand}}]\label{Taltransp}Let $\alpha \geq 1$.
Let $\psi$ be the monotone rearrangement sending $\nu_1$ to $\nu_{\alpha}$. Denote by $\Delta$ the function defined for any $s\geq 0$ by,
\begin{equation}\label{defdelta} \Delta(s) = \inf_{|x-y| = s} |\psi^{-1}(x) - \psi^{-1}(y)|.\end{equation}
There is a constant $c>0$ depending on $\alpha$ such that for any $s\geq 0$,
$$\Delta(s) \geq c \max(s ,s^{\alpha}).$$
\end{Lem}

\begin{Rem}\label{comprearrmap}
In \cite[Lemma 2.5]{Talagrand}, this estimate is derived for the monotone rearrangement $\phi$ of $\mu_1$ onto $\mu_{\alpha}$. But since,
\begin{equation} \label{link}\forall x \in \RR, \  \psi(x) = \sg(x) \phi(|x|),\end{equation}
one easily deduces the same estimate for $\psi$, together with the fact that if $x,y$ have opposite signs,
\begin{align*}
 \phi^{-1}(|y|) + \phi^{-1}(|x|) &\geq c \big(\max(|x|,|x|^{\alpha}) + \max( |y|,|y|^{\alpha}) \big)\\
& \geq c' \max(|x - y|, |x-y|^{\alpha}),
\end{align*}
where $c'$ is some constant and where we used the fact that $|x-y| = |x|+|y|$.
\end{Rem}

%
%
To get the exact asymptotic of the tail distribution of $\nu_{\alpha}$ we will need of the following finer estimate on the monotone rearrangement.

\begin{Lem}\label{brenier1}Let $\alpha \geq1$.
Define for any $m\geq 1$,
\begin{equation}\label{defdeltam}\forall s\geq 0, \ \Delta_m(s) = \inf \{ |\psi^{-1}(x) - \psi^{-1}(y)| : |x|\leq m, \ |x-y|=s \}.\end{equation}
There is a constant $\gamma$ depending on $\alpha$, such that for any $\eps \in (0,1)$, and $s \geq m\eps^{-1}$,
$$ \Delta_m(s) \geq 
(1-\gamma \eps) s^{\alpha}.$$

\end{Lem}
\begin{proof}
By definition of $\psi$, we have for any $x\in \RR$,
$$ \psi^{-1}(x) = -\sg(x) \log \int_{|x|}^{+\infty}  e^{-u^{\alpha}} \frac{du}{Z_{\alpha}},$$
where $Z_{\alpha}$ is the normalizing constant of $\mu_{\alpha}$.
Let $s \geq m\eps^{-1}$ and $x,y\in \RR$ such that $0\leq |x| \leq m$, and $|x-y| =s$. If $x$ and $y$ have the same signs, we can assume without loss of generality, that both $x, y\geq 0$. As $x\leq m \leq s$, we have $y= x+s$. Thus,
$$ \psi^{-1}(y) - \psi^{-1}(x) \geq \psi^{-1}(s) - \psi^{-1}(m).$$  
 We have, on one hand, as $s\geq 1$,
$$ \int_{s}^{+\infty} e^{-u^{\alpha}} du \leq \frac{1}{\alpha}  \int_{s}^{+\infty} \alpha u^{\alpha-1} e^{-u^{\alpha}} du = \frac{1}{\alpha}e^{-s^{\alpha}}.$$
And on the other hand,
$$ \int_m^{+\infty} e^{-u^{\alpha}} du \geq e^{-(m+1)^{\alpha}}.$$
Therefore, as $s \geq m\eps^{-1}$,
$$ \psi^{-1}(y) - \psi^{-1}(x) \geq s^{\alpha}  -(m+1)^{\alpha} +\log \alpha \geq s^{\alpha}(1-\gamma \eps^{\alpha}),$$
for some constant $\gamma>0$.  
Now, if $x$ and $y$ have opposite signs, we can assume without loss of generality that $x\leq 0$ and $y\geq 0$. Then, $y\geq s-m$ so that,
$$ |\psi^{-1}(y) - \psi^{-1}(x) | =  \psi^{-1}(y)+\psi^{-1}(-x)  \geq \psi^{-1}(s-m) \geq (s-m)^{\alpha}+\log \alpha.$$
Thus, we can find some constant $\gamma'$ such that $|\psi^{-1}(y) - \psi^{-1}(x) |  \geq (1-\gamma' \eps) s^{\alpha}$. 

\end{proof}

\begin{Rem}
The truncation we performed here is made to ensure we get the best constant (that is $1$) in the estimate of the large increments of the monotone rearrangement. Indeed, defining $\Delta$ as in \eqref{defdelta}, we would get for $s\gg 1$,
$$\Delta(s) \leq \Big|\psi\Big( \frac{s}{2} \Big) - \psi\Big( \frac{-s}{2} \Big) \Big| = 2\psi\Big( \frac{s}{2} \Big)  \simeq 2\Big( \frac{s}{2}\Big)^{\alpha} = 2^{1-\alpha} s^{\alpha},$$
with $2^{1-\alpha}<1$.

\end{Rem}

When $\alpha<1$, we get the following estimate  on the monotone rearrangement of $\nu_1$ onto $\nu_{\alpha}$. Note that as $\nu_{\alpha}$ does not have an exponential tail, the rearrangement map cannot be a Lipschitz function.

\begin{Lem}\label{Brenieralpha}Let $\alpha \in (0,1)$.
Let $\phi$ be the monotone rearrangement of $\mu_1$ onto $\mu_{\alpha}$. There is a constant $K>0$ depending on $\alpha$ such that for any $x,y \in [0,+\infty)$, 
$$ |\phi(x) - \phi(y)|\leq K \max \Big(  |x-y|, x^{\frac{1}{\alpha}-1}|x-y|, |x-y|^{\frac{1}{\alpha}} \Big).$$

\end{Lem}

\begin{proof}
This proof is very much in the spirit of \cite[Lemma 2.5]{Talagrand}. We begin by bounding from above
$$\int_x^{+\infty} e^{-y^{\alpha}} dy,$$
when $x\geq1$. The change of variable $u= y^{\alpha}$ gives,
$$\int_x^{+\infty} e^{-y^{\alpha}} dy = \frac{1}{\alpha} \int_{x^{\alpha}}^{+\infty} u^{\frac{1}{\alpha}-1} e^{-u} du.$$
Let $m = \lceil \frac{1}{\alpha} \rceil$. Integrating by parts $m$ times, we get
\begin{align*}
\int_{x^{\alpha}}^{+\infty} u^{\frac{1}{\alpha}-1} e^{-u} du &= \sum_{k=1}^{m-1} \big(\frac{1}{\alpha}-1\big)...\big(\frac{1}{\alpha}-k+1\big) x^{1-k\alpha} e^{-x^{\alpha}} \\
&+\big(\frac{1}{\alpha}-1\big)...\big(\frac{1}{\alpha}-m+1\big)\int_{x^{\alpha}}^{+\infty} u^{\frac{1}{\alpha}-m} e^{-u} du.
\end{align*}
As $\frac{1}{\alpha} -m \leq 0$, we deduce for any $x\geq 1$,
$$\int_{x^{\alpha}}^{+\infty} u^{\frac{1}{\alpha}-1} e^{-u} du \leq K x^{1-\alpha}e^{-x^{\alpha}},$$
where $K>0$ is some constant depending on $\alpha$ which will vary along the proof. Therefore, for any $x\geq 1$,
 \begin{equation} \label{maj0}\int_x^{+\infty} e^{-y^{\alpha}} dy \leq Kx^{1-\alpha}e^{-x^{\alpha}}.\end{equation}
By definition $\phi$ satisfies for any $x> 0$,
\begin{equation} \label{eqtransport} e^{-x} = \int_{\phi(x)}^{+\infty} e^{-y^{\alpha}} \frac{dy}{Z_{\alpha}}.\end{equation}
This implies that $\phi$ is an increasing homeomorphism of $\RR_+$.
For $\phi(x)\geq 1$, we have
\begin{equation} \label{aprioribound}e^{-x} \leq K\phi(x)^{1-\alpha}e^{-\phi(x)^{\alpha}}.\end{equation}
From \eqref{eqtransport}, we see that $\phi$ is differentiable, and $\phi'$ satisfies for any $x\geq 0$,
$$ e^{-x} = \frac{1}{Z_{\alpha}}\phi'(x) e^{-\phi(x)^{\alpha}}.$$
Thus by \eqref{aprioribound}, we get for $t \geq \phi^{-1}(1)$, 
\begin{equation} \label{derivtransp} \phi'(t) \leq K\phi(t)^{1-\alpha}.\end{equation}
Dividing by $\phi(t)^{1-\alpha}$ and integrating on $[\phi^{-1}(1),x]$ we get
$$ \phi(x)^{\alpha} -1 \leq K(x-\phi^{-1}(1)),$$
for any $x \geq \phi^{-1}(1)$. Hence, 
\begin{equation} \label{bornepsi} \phi(x)\leq Kx^{\frac{1}{\alpha}},\end{equation}
for $x\geq \phi^{-1}(1)$. By \eqref{derivtransp} we deduce
$$\phi'(x) \leq K x^{\frac{1}{\alpha}-1}.$$
Since $\phi'$ is continuous, at the price of taking $K$ larger, we have
$$\forall x \geq 0, \ \phi'(x) \leq K \max(1, x^{\frac{1}{\alpha}-1}).$$
Let $x \geq 0$, and $y \in \RR$ such that $x +y \geq 0$. If $x,x+y\leq 1$,
$$ |\phi(x+y) - \phi(x)| \leq Ky.$$
Whereas if $x,x+y\geq 1$,
$$ |\phi(x+y) - \phi(x)| \leq K \int_x^{x+y} t^{\frac{1}{\alpha}-1} dt = \alpha K\big( (x+y)^{\frac{1}{\alpha}} - x^{\frac{1}{\alpha}}\big).$$
Now, if $ 0\leq x \leq 1\leq x+y$, 
\begin{align*}
 |\phi(x+y) - \phi(x)|& \leq K\int_{x}^{x+y}(1+t^{\frac{1}{\alpha}-1})dt\\
& \leq K\big( y + \alpha\big( (x+y)^{\frac{1}{\alpha}}-x^{\frac{1}{\alpha}}\big).
\end{align*}
In conclusion, for any $x\geq 0$, $x+y\geq 0$,
\begin{equation} \label{controlpsi} |\phi(x+y) - \phi(x)|\leq K\max\Big( y,\big( (x+y)^{\frac{1}{\alpha}}-x^{\frac{1}{\alpha}}\big)\Big).\end{equation} 
The mean value theorem yields
$$ |(x+y)^{\frac{1}{\alpha}} - x^{\frac{1}{\alpha}} | \leq \frac{1}{\alpha} \max\big( x^{\frac{1}{\alpha}-1}, (x+y)^{\frac{1}{\alpha}-1}\big) y.$$
Using the convexity of $x \mapsto |x|^{\frac{1}{\alpha}-1}$, if $1/\alpha\geq1$, or its sub-additivity, when $1/\alpha -1 \in (0,1)$,  we get
$$ |(x+y)^{\frac{1}{\alpha}} - x^{\frac{1}{\alpha}} | \leq \frac{a_{\alpha}}{\alpha} \max\big( x^{\frac{1}{\alpha}-1}, x^{\frac{1}{\alpha}-1} + y^{\frac{1}{\alpha}-1} \big) y,$$
with $a_{\alpha} = \max(1, 2^{\frac{1}{\alpha}-2})$. Together with \eqref{controlpsi}, this gives the claim.
\end{proof}

As in the case $\alpha\geq 1$, we can refine the estimate of Lemma \ref{Brenieralpha} to get the following result.

\begin{Lem}\label{trans0}Let $\alpha \in (0,1)$. Let $\psi$ be the monotone rearrangement of $\nu_1$ onto $\nu_{\alpha}$.
Let $\eps\in(0,1)$. Define the function $\Delta_m$ by,
$$\forall s \geq 0, \ \Delta_m(s) = \inf \big\{  |\psi^{-1}(y) - \psi^{-1}(x) | :  |x| \leq m, \ |x-y|=s \big\}.$$
There is some constant $\gamma>0$, such that 
$$\Delta_m(s) \geq 
\begin{cases}
\gamma^{-1} (m/\eps)^{\alpha-1} s& \text{ if } s <\frac{m}{\eps},\\
 \big(1-\gamma \eps^{\alpha/2}\big)|s|^{\alpha} & \text{ if } s\geq \frac{m}{\eps}.
\end{cases}$$
\end{Lem}

\begin{proof}Since $\phi$ and $\psi$ are linked by the the relation \eqref{link}, the same estimate as in Lemma \ref{Brenieralpha} holds for the Brenier map $\psi$. Therefore, we have for any $|s| \leq \psi^{-1}(m)$, and $t\in \RR$,
$$ | \psi(t) - \psi(s)| \leq K\max\big( \psi^{-1}(m)^{\frac{1}{\alpha}-1} |t-s|, |t-s|^{\frac{1}{\alpha}} \big),$$
with $K\geq 1$. 
Fix $|x| \leq m$, and $y\in\RR$. We have
$$ | \psi^{-1}(y) - \psi^{-1}(x) | \geq K^{-\alpha} \min(|y-x|^{\alpha}, \psi^{-1}(m)^{1-\frac{1}{\alpha}} |y-x| ),$$
But we know from \eqref{bornepsi} that  for $m\geq 1$, $ \psi^{-1}(m)\geq c_0 m^{\alpha}$, with some constant $c_0>0$. Thus, for $m\geq 1$, there is a constant $\gamma>0$, which will vary along the proof without changing name, such that
$$ | \psi^{-1}(y) -  \psi^{-1}(x) | \geq \gamma^{-1} \min(|y-x|^{\alpha}, m^{\alpha-1} |y-x| ).$$
 We deduce that for $|y-x| \leq m/\eps$,
$$ | \psi^{-1}(y) - \psi^{-1}(x) | \geq \gamma^{-1} \Big(\frac{m}{\eps}\Big)^{\alpha- 1} |y-x|.$$
Let $s = |y-x|$. Assume now $s \geq m/\eps$. 
Proceeding as in the proof of Lemma \ref{brenier1} in the case $\alpha\geq 1$, we assume first that $x,y\geq 0$. As $s\geq m \geq x$, we must have $y= x+s$. Then,
$$| \psi^{-1}(y)-\psi^{-1}(x)| \geq \psi^{-1}(s) - \psi^{-1}(m).$$
On one hand, as $\alpha<1$, we have using the sub-additivity of $u\in \RR^+ \mapsto u^{\alpha}$,
$$ \int_m^{+\infty} e^{-u^{\alpha}} du = \int_0^{+\infty} e^{-(u+m)^{\alpha} } du \geq \big( \int_0^{+\infty} e^{-u^{\alpha}} du \big) e^{-m^{\alpha}}= \frac{1}{C}e^{-m^{\alpha}},$$
and on the other hand, by \eqref{maj0},
$$ \int_{s}^{+\infty} e^{-u^{\alpha}} du \leq C s^{1-\alpha} e^{-s^{\alpha}},$$
where $C$ is some constant depending on $\alpha$. Thus, 
$$| \psi^{-1}(y)-\psi^{-1}(x)| \geq s^{\alpha} - m^{\alpha} -(1-\alpha )\log s -2\log C.$$
As $\log s \leq (2/\alpha) s^{\alpha/2}$ for $s\geq 1$, we deduce that
$$| \psi^{-1}(y)-\psi^{-1}(w)| \geq s^{\alpha}(1-\gamma  \eps^{\alpha/2}).$$
 If $x$ and $y$ have opposite signs, we can assume $x\leq 0$ and $y\geq 0$, thus $y= s -m$ and we get,
\begin{align*}
| \psi^{-1}(y)-\psi^{-1}(x)|& \geq \psi^{-1}(y) \geq \psi^{-1}(s-m) \\
&\geq (s-m)^{\alpha} -(1-\alpha) \log (s-m) - \log C.
\end{align*}
As $s\geq m/\eps$, we deduce
$$| \psi^{-1}(y)-\psi^{-1}(x)| \geq s^{\alpha}(1-\gamma  \eps^{\alpha/2}),$$
which ends the proof of the claim.
\end{proof}

\subsection{A family of weights for $\nu_{\alpha}$}
Using transport arguments, we will work in this section at obtaining a family of weights for $\nu_{\alpha}$ which capture its exact tail distribution.

\begin{Pro}\label{IC1}Let $\alpha
> 0$, $\alpha\neq 1$, and $m\geq 1$.
There exist some constants $ \kappa, \eps_0>0$ depending on $\alpha$ such that for any $\eps\in (0,\eps_0)$, $(\nu_{\alpha}, w_{\alpha,\eps }^{(m)})$ satisfies the $[-m,m]$-truncated $\tau$-property where,
$$ w_{\alpha,\eps }^{(m)}(t) = \begin{cases}
\kappa^{-1} e^{-(\frac{m}{\eps})^{\alpha/2}} t^2 & \text{ if } |t| \leq m\eps^{-1},\\
(1- \kappa \eps^{(\alpha/2)\wedge 1}) |t|^{\alpha}  & \text{ if } |t| >  m \eps^{-1}.
\end{cases}$$
\end{Pro}

\begin{proof}
Let $\eps\in (0,1)$ and $m\geq 1$. Let $\delta>0$ such that
\begin{equation*} \frac{1}{2} \Big( \frac{m}{\eps} \Big)^{\alpha} = \frac{2}{\delta^2}.\end{equation*}
With this choice of $\delta$, we will prove that for $s\geq 0$,
$$w_{\delta}(\Delta_m(s)) \geq w_{\eps,\alpha}^{(m)}(s),$$
with the appropriate constants $\kappa$ and  $\eps_0$, $w_{\delta}$  defined in Corollary \ref{corexp}, and where $\Delta_m$ is as in \eqref{defdeltam}. Using the result of Lemma \ref{transp0}, this will yield the claim.

Let $\eps$ be small enough such that $w_{\delta}$ is non-decreasing. This is possible since $\delta^2 \leq 2 \eps^{\alpha}$.
Let $s\geq m/\eps$. If $\eps$ is small enough, we have by Lemma \ref{brenier1} or \ref{trans0},
$$ \Delta_m(s) \geq  \frac{1}{2} \Big(\frac{m}{\eps}\Big)^{\alpha} = \frac{2}{\delta^2}.$$
If $\alpha>1$, then by  Lemma \ref{brenier1} we get, as $\delta^2 \leq 4\eps^{\alpha}$,
$$w_{\delta}(\Delta_m(s)) \geq (1-2\delta)(1-\gamma \eps) s^{\alpha}\geq (1- \kappa \eps^{(\alpha/2)\wedge 1}) s^{\alpha},$$
for some constant $\kappa$ which will vary along the proof.
Similarly, when $\alpha<1$, we get by Lemma \ref{trans0},
$$w_{\delta}(\Delta_m(s)) \geq (1-2\delta)(1-\gamma \eps^{\alpha/2}) s^{\alpha}\geq (1- \kappa \eps^{\alpha/2}) s^{\alpha}.$$
Now let $ s \leq m/\eps$. Assume $\alpha \geq 1$. By Lemma \ref{Taltransp} and the fact that $w_{\delta}$ is non-decreasing, we have
$$w_{\delta}(\Delta_m(s)) \geq w_{\delta}(c s),$$
where $c$ is some positive constant.
Without loss of generality, we can assume $c\leq 1/2$. Then, as $m \eps^{-1}\leq 4\delta^{-2}$, we have $cs \leq 2\delta^{-2}$, so that we get
$$w_{\delta}(\Delta_m(s)) \geq \frac{c^2\delta e^{-\frac{1}{\delta}} s^2}{8}.$$
Using the fact that $\delta e^{ -\frac{1}{\delta}} \geq c_1 e^{-2/\delta}$, for some constant $c_1>0$, we get the claim in the case $\alpha>1$.
Assume now $\alpha<1$. From Lemma \ref{trans0} and the fact that $w_{\delta}$ is non-decreasing, we deduce
$$w_{\delta}(\Delta_m(s)) \geq w_{\delta}(\gamma^{-1} (m/\eps)^{\alpha-1} s).$$
Without loss of generality, we can assume that $\gamma \geq 2$. As $m\eps^{-1} \leq 4\delta^{-2}$ and $s\leq m/\eps$, we have 
$$\gamma^{-1} (m/\eps)^{\alpha-1} s \leq \frac{2}{\delta^2}.$$
Thus,
$$w_{\delta}(\Delta(s)) \geq \frac{1}{8} \delta e^{-\frac{1}{\delta}} \big( \gamma^{-1} (m/\eps)^{\alpha-1} s \big)^2\geq \kappa^{-1}\delta^{a} e^{-1/\delta},$$
with some $a>0$. But, we can find some constant $c_2>0$ such that
$$\delta^{a} e^{- 1/\delta} \geq c_2 e^{-2/\delta},$$
which, recalling that $(m\eps^{-1})^{\alpha} = 4\delta^{-2}$ gives the claim.

\end{proof}

We can now give a proof of Proposition \ref{tau}.

\begin{proof}[Proof of Proposition \ref{tau}]
 As  $(\nu_{\alpha}, w_{\alpha,\eps}^{(m)})$ satisfies the $[-m,m]$-truncated $\tau$-property for $\eps \in (0,\eps_0)$, for some $\eps_0>0$ and any $m\geq 1$ by Proposition \ref{IC1}, we deduce by the tensorization property of the $\tau$-property (see Lemma \ref{tensortau}) that $(\nu_{\alpha}^{ n}, W_{\alpha,\eps}^{(m)})$ satisfies the $mB_{\ell^{\infty}}$-truncated $\tau$-property with $W_{\alpha,\eps}^{(m)}$ defined as in \eqref{weight0}. 

\end{proof}

\section{Large deviations}\label{prooftheorem}
We will prove in this section Theorem \ref{theoremgene}. As sketched in the introduction, the proof will consist in looking for, in a first phase, large deviations inequalities for $\nu_{\alpha}^n$ and lower bounds estimates of the probability of translates. 

As a consequence of the truncated $\tau$-property of Proposition \ref{tau}, satisfied by $\nu_{\alpha}^n$ and the weight functions $W_{\alpha, \eps}^{(m)}$, we deduce an isoperimetric-type bound for $\nu_{\alpha}^n$ with respect to the metric $|| \ ||_{\ell^{\alpha}}$ (or $|| \ ||_{\ell^{\alpha}}^{\alpha}$ in the case $\alpha<1$). This estimate will be of paramount importance to derive the upper bound of Theorem \ref{theoremgene}.
\begin{Pro}\label{grandev}Let $\alpha>0$, $\alpha \neq 2$.  Let $r>0$.
 Let $v(n)$, $t(n)$ be two sequences going to $+\infty$ as $n$ goes to $+\infty$. Let $E$ and $F$ be Borel subsets of $\RR^n$ such that
$$ F +t(n) B_{\ell^2} \subset E, \ \liminf_{n\to+\infty}\nu_{\alpha}^n( F) >0.$$
For $\alpha \neq 1$, we assume that 
$$(\log n)^{\alpha/2} = o(\log \frac{t(n)^2}{v(n)}),$$
whereas for $\alpha=1$, we assume $v(n) = o(t(n)^2)$.  Then,
\begin{equation} \label{ineqdevPro} \limsup_{n\to +\infty} \frac{1}{v(n)} \log \nu_{\alpha}^n\big( x\notin E +  (rv(n))^{1/\alpha} B_{\ell^{\alpha}} \big ) \leq -r.\end{equation}

\end{Pro}
\begin{Rem}

For $\alpha=2$, the Gaussian isoperimetric inequality (see \cite[Theorem 2.5]{Ledouxmono}) entails the same result without any further assumption on the speed $v(n)$ or the set $E$ than $\liminf_n \nu_2^n(E)>0$.
\end{Rem}
\begin{proof} 
Before going into the proof per say, we need to relate the enlargements by the weights $W_{\alpha,\eps}^{(m)}$, for which we know that $(\nu_{\alpha}^n, W_{\alpha,\eps}^{(m)})$ satisfies the $\tau$-property, and therefore a deviation inequality of the type \eqref{devtrunc}, to  the $\ell^{\alpha}$-balls. This is the subject of the following lemma.
\begin{Lem} \label{decoupfnpoids}Let $\alpha>0$. With the notation of Proposition \ref{tau}, for any $r>0$, $m\geq 1$ and $\eps \in (0,\eps_0)$,
$$\big\{ W_{\alpha,\eps}^{(m)} \leq r\big(1-\kappa \eps^{(\alpha/2)\wedge 1}\big) \big\}  \subset k_m(\eps)\sqrt{r} B_{\ell^2} +r^{1/\alpha} B_{\ell^{\alpha}},$$
with $k_m(\eps) =  \sqrt{\kappa} e^{\frac{1}{2} (\frac{m}{\eps})^{\alpha/2} }$.
Moreover, there is a function $l : \RR_+ \to\RR_+$, such that
$$\{ W_{1,\eps} \leq r(1-2\eps) \} \subset l(\eps) \sqrt{r} B_{\ell^2} + rB_{\ell^1}.$$
\end{Lem}
\begin{proof}We will prove only the first statement, the proof for the second one being similar.
Let $y\in\RR^n$. By cutting the entries of $y$, we can find $y_1,y_2\in \RR^n$, such that $y=y_1+y_2$, for any $i\in \{1,...,n\}$, $y_1(i)y_2(i) = 0$, and
$$ |y_1(i)| \leq \frac{m}{\eps}, \quad |y_2(i)|> \frac{m}{\eps}.$$
By the very definition of $W_{\alpha, \eps}^{(m)}$,
$$ \kappa^{-1}  e^{-(m/\eps)^{\alpha/2}}\sum_{i=1}^n |y_1(i)|^2 = W_{\alpha, \eps}^{(m)}(y_1) \leq  W_{\alpha, \eps}^{(m)}(y),$$
and

$$(1-\kappa \eps^{(\alpha/2)\wedge 1}) ||y_2||_{\ell^{\alpha}}^{\alpha} = W_{\alpha,\eps}^{(m)}(y_2) \leq W_{\alpha,\eps}^{(m)}(y).$$
Thus, if we let
$$k_m(\eps)^2 = e^{(m/\eps)^{\alpha/2}} \kappa,$$
and if $W_{\alpha,\eps}^{(m)}(y)\leq r(1-\kappa \eps^{(\alpha/2)\wedge 1})$, then $||y_1 ||_{\ell^2}\leq  k_m(\eps) \sqrt{r}$, and $||y_2||_{\ell^{\alpha}}^{\alpha}\leq r$.

\end{proof}
With this lemma proven, we can now give the proof of Proposition \ref{grandev}. We start with the case $\alpha=1$. As $v(n) = o(t(n)^2)$, for $n$ large enough, we have $l(\eps) \sqrt{rv(n)} \leq t(n)$. Then, by Lemma \ref{decoupfnpoids}, we have
$$ F + \{W_{1,\eps} \leq r(1-2\eps) v(n)\} \subset F + t(n)B_{\ell^2} + r v(n) B_{\ell^{1}}.$$
But by assumption, $ F +t(n) B_{\ell^2} \subset E$. Thus,
$$ F + \{W_{1,\eps} \leq r(1-2 \eps) v(n)\} \subset  E + r v(n) B_{\ell^{1}}.$$
We deduce that,
$$\nu_{1}^n( x \notin E + r v(n) B_{\ell^{1}}) \leq \nu_{1}^n\big( x \notin  F + \{W_{1,\eps} \leq r(1-2\eps) v(n)\}\big).$$
As $(\nu_{1}^{ n}, W_{1,\eps})$ satisfies the $\tau$-property by Corollary \ref{corexp}, we have the following deviation inequality (see \eqref{dev}),
$$\nu_{1}^n( x \notin E + r v(n)B_{\ell^{1}}) \leq \frac{1}{\nu_{1}^n(F)}e^{-r(1-2 \eps)v(n)}.$$
As $\liminf_{n} \nu_{1}^n(F)>0$, we get
$$\limsup_{n\to +\infty} \frac{1}{v(n)} \log \nu_{1}^n( x \notin E +r v(n) B_{\ell^{1}} ) \leq -r(1-2 \eps).$$
Letting $\eps$ going to $0$, we get the claim. 

Let now $\alpha \neq 1$.
 Let $\eps \in (0,\eps_0)$ and set $m = c (\log n)^{1/\alpha}$, with some $c>0$ which is to be chosen later. By Lemma \ref{decoupfnpoids}
$$ F + \{W_{\alpha,\eps}^{(m)} \leq r(1-\kappa \eps^{(\alpha/2)\wedge 1}) v(n)\} \subset F + k_n(\eps) \sqrt{rv(n)}B_{\ell^2} + (r v(n))^{1/\alpha} B_{\ell^{\alpha}}).$$
From the assumption that $(\log n)^{\alpha/2} = o(\log \frac{t(n)}{\sqrt{v(n)}})$ we deduce that
for $n$ large enough,
$$\frac{t(n)}{\sqrt{v(n)}}\geq e^{ (\frac{(\log n)}{\eps})^{\alpha/2}}.$$
In particular for $n$ large enough,
$$\sqrt{\frac{\kappa}{r}}\frac{t(n)}{\sqrt{v(n)}}\geq e^{\frac{1}{2} (\frac{(\log n)}{\eps})^{\alpha/2}}.$$
Put in another way
$$k_m(\eps) \sqrt{rv(n)} \leq t(n).$$
Thus,
$$ F + \{W_{\alpha,\eps}^{(m)} \leq r(1-\kappa \eps^{(\alpha/2)\wedge 1}) v(n)\} \subset F + t(n)B_{\ell^2} + (r v(n))^{1/\alpha} B_{\ell^{\alpha}}.$$
As by assumption $F + t(n)B_{\ell^2} \subset E$, we get
$$\nu_{\alpha}^n ( x \notin E + (r v(n))^{1/\alpha} B_{\ell^{\alpha}}) ) \leq \nu_{\alpha}^n( x \notin  F + \{W_{\alpha,\eps}^{(m)} \leq r(1-\kappa \eps) v(n)\}).$$
As $(\nu_{\alpha}^{ n}, W_{\alpha,\eps}^{(m)})$ satisfies the $mB_{\ell^{\infty}}$-truncated $\tau$-property by Proposition \ref{tau}, we deduce the following deviation inequality (see \eqref{devtrunc}),
$$\nu_{\alpha}^n \big( x \notin  F + \{W_{\alpha,\eps}^{(m)} \leq r(1-\gamma \eps^{(\alpha/2)\wedge 1}) v(n)\}\big) \leq \frac{1}{\nu_{\alpha}^n(F\cap mB_{\ell^{\infty}}) }e^{-r(1-\gamma \eps^{(\alpha/2)\wedge 1 })v(n)}.$$
But, 
$$ \int ||x||_{\infty} d\nu_{\alpha}^n(x)= 
\int ||x||_{\infty} d\mu_{\alpha}^n(x).$$
Let $\Phi = \phi^{\otimes n}$, defined by $\Phi(x) = (\phi(x_i))_{1\leq i \leq n}$, where $\phi$ is the monotone rearrangement map sending $\mu_1$ to $\mu_{\alpha}$. Then $\Phi$  sends $\mu_1^n$ to $\mu_{\alpha}^n$, so that,
$$\int ||x||_{\infty} d\mu_{\alpha}^n(x) 
= \int ||\Phi(x)||_{\infty} d\mu_1^n(x).$$
From \eqref{bornepsi}, we deduce
$$\int ||\Phi(x)||^{\alpha}_{\infty} d\mu_1^n(x)\leq K( 1+ \int ||x||_{\infty}^{1/\alpha} d\mu_1^n(x)),$$
 for some constant $K>0$. But $\int ||x||_{\infty}^{1/\alpha} d\mu_1^n(x) \leq c_0 (\log n)^{1/\alpha}$, for some constant $c_0\geq 1$. Therefore,
\begin{equation} \label{momentalpha} \int ||x||_{\infty} d\mu_{\alpha}^n(x) \leq 2Kc_0(\log n)^{1/\alpha}. \end{equation}
Thus by Markov's inequality, 
$$\nu_{\alpha}^n(x\notin mB_{\ell^{\infty}})\leq \frac{2Kc_0}{c},$$
since we chose $m = c (\log n)^{1/\alpha}$.
As $\liminf_n \nu_{\alpha}^n(F)>0$ by assumption, we deduce that for $c$ large enough,
$$ \liminf_{n\to +\infty} \nu_{\alpha}^n( F\cap mB_{\ell^{\infty}} )>0.$$
Therefore,
$$\limsup_{n\to +\infty} \frac{1}{v(n)} \log \nu_{\alpha}^n (x \notin E + (r v(n))^{1/\alpha} B_{\ell^{\alpha}}\}) \leq -r(1-\kappa\eps^{(\alpha/2)\wedge 1}),$$
which gives the claim by taking $\eps \to 0$.
\end{proof}

We show in the next proposition that we can bound from below the probability of translates under $\nu_{\alpha}^n$.
\begin{Pro}\label{lowb}Let $\alpha\in (0,2]$.
Let $v(n)$ be a sequence going to $+\infty$ as $n$ goes to $+\infty$. Fix some $r>0$. Let $E$ be some Borel subset of $\RR^n$ such that
 $$ \liminf_{n \to +\infty} \nu_{\alpha}^n( E)>0.$$ 
(i). For any sequence $h_n$ of elements of $\RR^n$,
$$\liminf_{n \to +\infty}\frac{1}{v(n)} \log \nu_{\alpha}^n( E + v(n)^{1/\alpha}h_n)\geq - \limsup_{n \to +\infty} || h_n||_{\ell^{\alpha}}^{\alpha}.$$
(ii). If $\alpha \in (0,1]$, then for any sequence $h_n \in \RR_+^n$,
$$\liminf_{n \to +\infty}\frac{1}{v(n)} \log \mu_{\alpha}^n( E + v(n)^{1/\alpha}h_n)\geq - \limsup_{n \to +\infty} || h_n||_{\ell^{\alpha}}^{\alpha}.$$
\end{Pro}

\begin{Rem}
On can obtain the estimate $(ii)$ when $\alpha \in (1,2]$ for the measures $\mu_{\alpha}$ with the additional assumption $n = o(v(n))$ on the speed, which is actually very restrictive in the applications we have in mind. This is one of the reasons of the limitation of Theorem \ref{theoremgenesup} to the case $\alpha\leq 1$, since we do not know how to produce a meaningful lower bound of such translated sets in this case. 
Similarly, when $\alpha>2$, one can see, at least for $\alpha$ integer, that the estimate $(i)$ does not hold unless $n = o(v(n))$.

\end{Rem}

\begin{proof}The proof will essentially follow the lines of \cite[Theorem 5.1]{Ledouxflour}. Indeed, in the Gaussian case $\alpha=2$, this lower bound is derived from the translation formula of the Gaussian measure. The proof for $\alpha<2$ will consist in mimicking the Gaussian case.

If the $\limsup$ in the right-hand side of $(i)$ is infinite, then the statement is trivial. If it is finite, we take some $\tau>0$, such that $|| h_n||_{\ell^{\alpha}}^{\alpha} \leq \tau$, for all $n\in \NN$. Let for any $h \in \RR^n$, $W_{\alpha}(h ) = \sum_{i=1}^n |h_i|^{\alpha}$.
Then, we have,
$$\nu_{\alpha}^n( E + v(n)^{1/\alpha}h_n) =\frac{1}{Z_{n}} \int_{ E} e^{-W_{\alpha}(y+v(n)^{1/\alpha}h)} d{\ell_n}(y),$$
where $\ell_n$ denotes the Lebesgue measure on $\RR^n$, and $Z_{n}$ is the normalizing factor.
If $\alpha \in (0,1]$, then for any $s,t \in \RR$,
$$|s+t|^{\alpha} \leq |s|^{\alpha} +|t|^{\alpha}.$$
Thus,
$$W_{\alpha}( y+v(n)^{1/\alpha} h_n ) \leq W_{\alpha}(y) + v(n)W_{\alpha}( h_n)$$
Therefore,
$$\nu_{\alpha}^n ( E + v(n)^{1/\alpha}h_n) \geq e^{-v(n)W_{\alpha}(h_n) }\nu_{\alpha}^n(E),$$
which gives the claim in the case $\alpha \in (0,1)$. Note that the same argument for $\mu_{\alpha}$ instead of $\nu_{\alpha}$ gives without changes the estimate $(ii)$.

Now, if $\alpha \in (1,2]$, we have for any $s, t \in \RR$, 
\begin{equation} \label{inesub1}| s+t|^{\alpha} \leq |s|^{\alpha} + \alpha \mathrm{sg}(st)|s|^{\alpha-1} |t| + |t|^{\alpha},\end{equation}
where $\mathrm{sg}(st)$ stands for the sign of $st$.
Thus, for any $y,h \in \RR^n$, 
$$ W_{\alpha}(y+v(n)^{1/\alpha}h)\leq W_{\alpha}(y)+\alpha v(n)^{1/\alpha} V(y,h) + v(n)W_{\alpha}( h),$$
where 
\begin{equation} \label{defPhi}V(y,h) = \sum_{i=1}^n v(y_i,h_i),\end{equation}
and  $v(y,h) =  \mathrm{sg}(y h)|y|^{\alpha-1} |h|$. We have,
\begin{align*} 
\frac{1}{Z_{n}} \int_{E} e^{-W_{\alpha}(y+v(n)^{1/\alpha} h_n)} d\ell_n(y)
&\geq \frac{e^{-v(n) W_{\alpha} (h_n)} }{Z_{n} }
 \int_{ E} e^{-W_{\alpha} (y) -\alpha v(n)^{1/\alpha}V(y,h_n) } d{\ell_n(y)}\\
& = e^{-v(n) W_{\alpha} (h_n)} \int_{E}e^{ -\alpha v(n)^{1/\alpha} V(x,h_n)}d \nu_{\alpha}^n(x).
\end{align*}
Jensen's inequality yields,
$$ \int_{E} e^{ -\alpha v(n)^{1/\alpha} V(x,h_n)} d \nu_{\alpha}^n(x) \geq \nu_{\alpha}^n(E)\exp\Big( -\frac{\alpha v(n)^{1/\alpha} }{\nu_{\alpha}^n(E)} \int_{E}  V(x,h_n)d \nu_{\alpha}^n(x)\Big).$$ 
But, by Cauchy-Schwarz inequality, 
$$ \int_{E}  V(x,h_n)d \nu_{\alpha}^n(x) \leq \nu_{\alpha}^n(E)^{1/2} \Big(\int V(x,h_n)^2 d\nu_{\alpha}^n(x) \Big )^{1/2}.$$
But $\int v(x,h)d\nu_{\alpha}(x) =0$ for any $h \in \RR$, since $v(-x,h) = -v(x,h)$ and $\nu_{\alpha}$ is symmetric. Thus,
$$\int V(x,h_n)^2 d\nu_{\alpha}^n(x) = \int |t|^{2(\alpha-1)}d\nu_{\alpha}^n(t) \big(\sum_{i=1}^n |h_n(i)|^2\big).$$
Using the fact that $\alpha \leq 2$, we get,
$$\Big(\int V(x,h_n)^2 d\nu_{\alpha}^n(x) \Big )^{\frac{\alpha}{2} } \leq   c^{\frac{\alpha}{2}} W_{\alpha}(h_n),$$
where $c>0$ is some constant. As $W_{\alpha}(h_n) \leq \tau$, we have
$$ \int_{ E} e^{ -\alpha v(n)^{1/\alpha} V(x,h_n)} d\nu_{\alpha}^n(x)  \geq \nu_{\alpha}^n( E) \exp\Big( -\frac{c^{1/2} \tau \alpha v(n)^{1/\alpha} }{\nu_{\alpha}^n( E) ^{1/2}}\Big).$$ 
Note that is was actually very important that we did not bound $\sg(xy)$ by $1$ in \eqref{inesub1}, so that $v(.,h)$ is of mean $0$ under $\nu_{\alpha}$, and $\int V(x,h_n)^2 d\nu_{\alpha}^n(x)$ is not too big. When one replaces $\nu_{\alpha}$ by $\mu_{\alpha}$, this is exactly where one  needs to make an assumption on the speed to identify the leading term.

By assumption, we know that there is some $\eta>0$ such that for $n$ large enough, $\nu_{\alpha}^n( E) >\eta $. Thus, we get for $n$ large enough,
$$\nu_{\alpha}^n( x \in E + v(n)^{1/\alpha} h_n) \geq \eta\exp\Big( -v(n) W_{\alpha}(h_n) -2 \Big(\frac{c}{\eta}\Big)^{1/2} \tau \alpha v(n)^{1/\alpha}\Big ).$$ 
Taking the $\liminf$ at the exponential scale $v(n)$, we get the claim.


\end{proof}

We can now give a proof of Theorem \ref{theoremgene}. We will essentially follow the proof of the LDP of Wiener chaoses (see \cite{LedouxWiener}), replacing the use of the Cameron-Martin formula by Proposition \ref{lowb}, and the Gaussian isoperimetric inequality with Proposition \ref{grandev}.

\begin{proof}[Proof of Theorem \ref{theoremgene}]
Without loss of generality we can and will assume that $N=\NN$.
\textbf{Property of the rate function: }
By assumption $(iv)$, for any $x \in \mathcal{X}$,
$$I_{\alpha}(x) = \sup_{\delta>0} \limsup_{n \to +\infty} I_{n,\delta}(x).$$
This formulation shows that $I_{\alpha}(x)< +\infty$ if and only if there is a sequence $h_n \in \RR^n$, such that
$$\lim_{n \to +\infty} F_n(h_n)=x,\quad \limsup_{ n \to +\infty} W_{\alpha}(h_n) = I_{\alpha}(x).$$
Thus, $I_{\alpha}(x) \leq \tau$, for some fixed $\tau\geq 0$, if and only if $x$ is a limit point of a sequence $(F_n(h_n))_{n\in N}$ such that $\limsup_n W_{\alpha}(h_n) \leq \tau$. Therefore, $I_{\alpha}$ is lower semi-continuous. 
Moreover,
$$ \{ I_{\alpha} \leq \tau \} \subset \overline{ \cup_{n \in \NN} F_n(2\tau B_{\ell^{\alpha}} ) }.$$
As by assumption $(iv)$ the set on the right-hand side is compact, we conclude that $I_{\alpha}$ is a good rate function.

\textbf{Lower bound: } 
Let $x \in \mathcal{X}$ such that $I_{\alpha}(x) < +\infty$. By assumption $(iv)$, there is a sequence $h_n \in \RR^n$ such that
$$\lim_{n \to +\infty} F_n(h_n)=x,\quad \limsup_{ n \to +\infty} W_{\alpha}(h_n) = I_{\alpha}(x).$$
Let $\delta>0$. For $n$ large enough,
$$\PP\big( f_n(X_n) \in B(x,2\delta) \big) \geq \PP\big( f_n(X_n) \in B(F_n(h_n), \delta) \big).$$
Let 
$$E = \big\{ Y \in \RR^n : d( f_n(Y + v(n)^{1/\alpha} h_n), F_n(h_n)) < \delta \big\}.$$
Note that 
$$\PP\big( f_n(X_n) \in B(F_n(h_n), \delta)\big) = \PP\big( X_n \in E + v(n)^{1/\alpha} h_n \big).$$
By assumption $(i)$, $\PP(X_n \in E)$ goes to $1$ as $n$ goes to $+\infty$. From Proposition \ref{lowb}, we deduce
$$\liminf_{n\to+\infty} \frac{1}{v(n)} \log \PP\big( f_n(X_n) \in B(x,2\delta)\big) \geq -I_{\alpha}(x).$$

\textbf{Upper bound: }
  Let $A$ be a closed subset of $\mathcal{X}$. We can assume without loss of generality that $\inf_A I_{\alpha}>0$. Let $r>0$ such that $\inf_A I_{\alpha}>r$. Put in another way,
$$A\cap \{ I_{\alpha}\leq r\} = \emptyset.$$
As $I_{\alpha}$ is a good rate function, we can find a $\delta>0$ such that 
$$ A \cap V_{\delta}(\{I_{\alpha} \leq r\}) = \emptyset,$$
where $V_{\delta}$ denotes the $\delta$-neighborhood for the distance $d$. Thus,
$$\PP\big( f_n(X_n) \in A\big) \leq \PP\big( f_n(X_n)\notin V_{\delta}(\{ I_{\alpha} \leq r \})\big).$$
Let 
$$U = \big\{ x \in \RR^n  : f_n(x)\in V_{\delta}(\{ I_{\alpha} \leq r \}) \big\}.$$
Define, similarly as for the lower bound, the event
$$ E_{\delta} = \big\{ x \in \RR^n : \sup_{ h \in r^{1/\alpha}B_{\ell^{\alpha}} } d\big( f_n(x+v(n)^{1/\alpha}h_n), F_n(h_n)\big ) < \delta \big\}.$$ 
By assumption $(i)$, we know that $\PP( X_n \in E_{\delta})$ goes to $1$ as $n$ goes to $+\infty$. We claim that 
$$ E_{\delta} +  (v(n)r)^{\frac{1}{\alpha}}B_{\ell^{\alpha}}  \subset  U.$$
Indeed, if $h_n \in r^{1/\alpha} B_{\ell^{\alpha}}$ and $x \in E_{\delta}$, then $I_{\alpha}(F_n(h_n)) \leq r$, from the definition \eqref{deftaux} of $I_{\alpha}$, and
$$ d(f_n(x + v(n)^{1/\alpha} h_n ), F_n(h_n))< \delta,$$
so that $x+v(n)^{1/\alpha} h_n \in U$.
With this observation we get,
$$\PP\big( f_n(X_n) \in A\big) \leq \PP\big( X_n \notin E_{\delta} +  (v(n)r)^{\frac{1}{\alpha}}B_{\ell^{\alpha}}  \big).$$
If $\alpha =2$, we get by the Gaussian isoperimetric inequality (see \cite[Theorem 2.5]{Ledouxmono}) for any $n$ large enough so that $\PP( X_n \in E_{\delta}) \geq 1/2$,
$$ \PP\big( X_n \notin E_{\delta} +\sqrt{v(n)r}B_{\ell^{2}}   \big)\leq e^{-v(n) r},$$
which gives the upper bound.

Let now $\alpha <2$, and $t = t_{\delta/4}$, where $t_{\delta/4}$ is given by assumption $(ii)$. With the notation of Theorem \ref{theoremgene}
define, 
$$F = E_{\delta/2}\cap \big\{ y \in \RR^n : \sup_{|| h_n||_{\ell^2} \leq t } \mathcal{L}_n(h_n) \leq \frac{\delta}{2} \big\}.$$
By Markov's inequality and assumption $(ii)$, we deduce 
$$ \PP\big(  \sup_{|| h_n||_{\ell^2} \leq t } \mathcal{L}_n(h_n)\leq \frac{\delta}{2}\big ) \geq \frac{1}{2}.$$
From assumption $(i)$, we deduce that $\liminf_n \PP(X_n \in F)>0$. Furthermore, we claim that
\begin{equation} \label{claiminclu} F +tB_{\ell^2} \subset E_{\delta}.\end{equation}
Recall that 
$$\mathcal{L}_n(h) = \sup_{X_n + rv(n)^{1/\alpha} B_{\ell^{\alpha} }}d\big( f_n(x+h) , f_n(x)\big).$$
Now, if $X_n \in F$ and $h \in  tB_{\ell^2}$, then by definition of $\mathcal{L}_n$, for all $k \in rB_{\ell^{\alpha}}$
$$d\big(f_n(X_n + v(n)^{1/\alpha}k+h),f_n(X_n+v(n)^{1/\alpha}k)\big) \leq \frac{\delta}{2},$$
which yields \eqref{claiminclu} by triangular inequality. Thus the requirements of Lemma \ref{grandev} are met, and we get
$$ \limsup_{n \to +\infty} \frac{1}{v(n)} \log \PP\big( f_n(X_n) \in A\big) \leq -r.$$
As this inequality is true for any $r < \inf_A I_{\alpha}$, we get the upper bound.

\end{proof}

We will end this section with the proof of Theorem \ref{theoremgenesup}.

\begin{proof}[Proof of Theorem \ref{theoremgenesup}]
We will follow the same steps as for the proof of Theorem \ref{theoremgene}. The compactness assumption $(iii)$, and the assumption $(iv)'$ yield that $I_{\alpha}$ is a good rate function. 
As shown in the proof of Theorem \ref{theoremgene}, a large deviations upper bound holds with speed $v(n)$ and rate function $I_{\alpha}$, under the assumptions $(i)-(ii)-(iii)$. Thus, we only have to prove the lower bound. Let $x \in \mathcal{X}$ such that $I_{\alpha}^+(x)<+\infty$.  We know that there is a sequence $h_n \in \RR_+^n$ such that 
$$\lim_{n\to +\infty} F_n(h_n) = x, \quad \limsup_{n\to +\infty} W_{\alpha}(h_n) = I_{\alpha}^+(x).$$
Proceeding as in the proof of Theorem \ref{theoremgene}, if $\delta>0$, then for $n$ large enough,
$$\PP\big( f_n(X_n) \in B(x,2\delta) \big) \geq \PP\big( f_n(X_n) \in B(F_n(h_n), \delta) \big).$$
Let 
$$E = \big\{ y \in \RR^n : d( f_n(y + v(n)^{1/\alpha} h_n), F_n(h_n)) < \delta \big\}.$$
Note that 
$$\PP\big( f_n(X_n) \in B(F_n(h_n), \delta)\big) = \PP\big( X_n \in E + v(n)^{1/\alpha} h_n \big).$$
By assumption $(i)$, $\PP(X_n \in E)$ goes to $1$ as $n$ goes to $+\infty$. From Lemma \ref{lowb}, we deduce
$$\liminf_{n\to+\infty} \frac{1}{v(n)} \log \PP\big( f_n(X_n) \in B(x,2\delta)\big) \geq  -I_{\alpha}^+(x),$$
which ends the proof of the lower bound. Due to assumption $(iv)'$ the lower bound and upper bound rate functions match so that a full LDP holds. 
\end{proof}

\section{Concentration inequalities}\label{Chapconc}
We will prove in this section the concentration inequalities of Propositions \ref{conclinearstat}, \ref{concspintro} and \ref{concvpintro} for the linear statistics, the empirical spectral measure and largest eigenvalue of Wigner matrices satisfying the concentration property $\mathcal{C}_{\alpha}$ introduced by definition \ref{defCalpha}.

\subsection{Some examples of Wigner matrices satisfying $\mathcal{C}_{\alpha}$}
Before going into the proofs, we will review some workable criterion for a Wigner matrix to satisfy the concentration property $\mathcal{C}_{\alpha}$ when $\alpha\in [1,2]$. The case of $\alpha=2$ of normal concentration has drawn most of the attention, and we refer the reader to \cite[section 8.5]{Ledouxmono}, \cite{GZconc} or also \cite[Part II]{GuionnetFlour} for a presentation of the different examples of classical models of random matrices having normal concentration. 

When $\alpha\in [1,2]$ we introduce the notion of \textit{Poincaré-type inequalities} in the finite-dimensional setting. Let $d_m$ be some distance on $\RR^m$. For a smooth function $f: \RR^m \to \RR$, we define the length of the gradient of $f$ with respect to the distance $d_m$ by,
$$ \forall x \in \RR^m,\  |\nabla f(x) | = \limsup_{y\to x} \frac{|f(y)-f(x)|}{d_m(y,x)}.$$
We say that a probability measure $\mu$ satisfies a \textit{Poincaré-type inequality} on $(\RR^m, d_m)$ if there is some $\lambda>0$, such that for any smooth $f :\RR^m \to \RR$,
$$ \lambda \Var_{\mu} f \leq \int |\nabla f|^2 d\mu,$$ 
where the length of the gradient is taken with respect to $d_m$. 

Following Gozlan \cite[Definition 1.1]{Gozlan}, we will say that a probability measure $\mu$ on $\RR^m$ satisfies $\mathbb{SP}(\omega_{\alpha}, \lambda)$ if it satisfies the Poincaré-type inequality on $(\RR^m,d_{\omega_{\alpha}})$ with spectral gap $\lambda$, where $d_{\omega_{\alpha}}$ is the distance defined in \eqref{defdalpha}.


By the results of Bobkov-Ledoux \cite[Corollary 3.2]{BobLed}, and Gozlan \cite[Proposition 1.2]{Gozlan}, we know that if a Wigner matrix $X$ has entries satisfying $\mathbb{SP}(\omega_{\alpha},\lambda)$, then it satisfies a two-level deviations inequality:  for any Borel subset $A$ of $\mathcal{H}_n^{(\beta)}$ such that $\PP( X\in A)\geq 1/2$, and $r>0$,
\begin{equation}\label{concpropalpha} \PP( X\notin A +  \sqrt{r} B_{\ell^2}+r^{\frac{1}{\alpha}} B_{\ell^{\alpha}} ) \leq e^{-L r},\end{equation}
where $L$ only depends on $\lambda$, and by \cite[Proposition 1.2]{Gozlan}) can be taken as
\begin{equation}\label{defLlambdaintro} L(\lambda) = \frac{w( \frac{\sqrt{\lambda}}{\kappa})}{16},\quad \kappa = \sqrt{18 e^{\sqrt{5}}},\end{equation}
 with $w(t) = \min(|t|^2, |t|)$ for any $t\in \RR$. In particular, such a Wigner matrix has concentration $\mathcal{C}_{\alpha}$.
%
%
\begin{Rem}
We note that when $\alpha>2$, the Poincaré-type inequality  $\mathbb{SG}(\omega_{\alpha},\lambda)$ yields a different deviation inequality (the one above is also true for $\alpha>2$ but not sharp) where the mixed enlargement is replaced by $\sqrt{r} B_{\ell^2}\cap r^{\frac{1}{\alpha}} B_{\ell^{\alpha}} $ (see \cite{Gozlan} for more details).
\end{Rem}


A workable criterion for a probability measure on $\RR$ of the form $\mu = e^{-V} dx$ is given by Gozlan \cite[Proposition 1.2]{Gozlan} in terms of a growth condition of the potential $V$. More precisely, if 
\begin{equation}\label{condpot} \liminf_{x \to \pm \infty} \frac{\mathrm{sg}(x)V'(x)}{x^{\alpha-1}} >0,\end{equation}
then $\mu$ satisfies $\mathbb{SG}(\omega_{\alpha},\lambda)$ on $\RR$. 
We mention also that a criterion is available in higher dimension (although more intricate) in \cite[Proposition 3.5]{Gozlan}, which one may use for the complex entries of Wigner matrices.


In the case $\alpha=1$ of the classical Poincaré inequality, we know by Bobkov \cite{Bobkovlogconc}  (or by Bakry, Barthe, Cattiaux, and Guillin \cite{BBCG}) that  any log-concave law on $\RR^n$ satisfies a Poincaré inequality with a certain spectral gap depending on the dimension.  
Thus, any Wigner matrix with entries whose laws are log-concave will satisfy $\mathcal{C}_1$. 

%
%
%
%
%
%

When $\alpha \in [1,2]$, the concentration property $\mathcal{C}_{\alpha}$ is equivalent (see \cite[Proposition 1.3]{Ledouxmono}) to the following deviation inequality of Lipschitz functions around their medians, which will be useful in the applications. 

\begin{Lem}\label{conclip}Let $\alpha \in [1,2]$.
Let $X$ be a Wigner matrices with entries satisfying $\mathcal{C}_{\alpha}$ for some $\kappa>0$. Let $f : \mathcal{H}_n^{(\beta)}\to \RR$ be a function respectively $L_2$-Lipschitz and $L_{\alpha}$-Lipschitz with respect to $|| \ ||_{\ell^2}$, and $||\ ||_{\ell^{\alpha}}$. Then, for any $t>0$,
$$ \PP( f(X)>m_f +t) \leq 2\exp\Big( - \min \Big( \frac{t^2}{4\kappa^2 L_2^2}, \frac{t^{\alpha}}{2^{\alpha}\kappa^{\alpha} L_{\alpha}^{\alpha}} \Big) \Big),$$
where $m_f$ denotes the median of $f(X)$.
\end{Lem}

\subsection{A deviation inequality for $\nu_{\alpha}^n$, $\alpha \in (0,1)$}\label{sectiondevineq0}
In the case $\alpha\in (0,1)$, we will show that the Wigner matrices in the class $\mathcal{S}_{\alpha}$ satisfy the concentration property $\mathcal{C}_{\alpha}$. This  fact will follow from the study of the concentration property of the product measures $\mu_{\alpha}^n$ and $\nu_{\alpha}^n$.
 It can be shown that the probability measure $\nu_{\alpha}^n$ satisfies a weak Poincaré inequality (see  \cite[Chapter 7 \S 7.5]{BGL}). The derivation of a deviations inequality from the weak Poincaré inequality has been investigated by Barthe, Cattiaux and Roberto \cite{barthe}, and yields a concentration inequality with respect to Euclidean enlargements. We will follow another path which consists, as it was the case for $\alpha\geq 1$, in transporting Talagrand's deviation inequality for the symmetric exponential law \eqref{Taldev} onto $\nu_{\alpha}$ with $\alpha <1$, using the estimate on the monotone rearrangement map proved in Lemma \ref{Brenieralpha}. We start with the one-sided probability measure $\mu_{\alpha}$.

\begin{Pro}\label{devexpo}
Let $n\in \NN$, $n\geq 2$, and $\alpha \in (0,1)$. There is a constant $c>0$ depending on $\alpha$, such that for any $r>0$, $A$ Borel subset of $\RR_+^n$, and $C>0$ such that $\mu_{\alpha}^{n}(A) > 1/C$,
$$\mu_{\alpha}^{n} \Big( x \notin A + C(\log n)^{\frac{1}{\alpha}-1}\big( \sqrt{r} B_{\ell^2} + r B_{\ell^1} \big) + r^{\frac{1}{\alpha}} B_{\ell^{\alpha}} \Big) \leq \frac{e^{-c r}}{\mu_{\alpha}^{n}(A) - 1/C}.$$

\end{Pro}

\begin{Rem}\label{remdev0}
This deviation inequality is not optimal in the sense that it fails to capture the Gaussian fluctuations of empirical means from the central limit theorem. This is due to the $(\log n)^{1/\alpha-1}$ factor in front of the $\ell^2$-ball, which comes from the fact that the increasing rearrangement from $\mu_1$ to $\mu_{\alpha}$ is not a Lipschitz function.

But on the other hand, the $(\log n)^{\frac{1}{\alpha}-1}$ factor seems to be sharp, since it yields a non-trivial deviation inequality for
$$ (\log n)^{\frac{1}{\alpha}-1}\big(\max_{1 \leq i \leq n} x_i - m \big),$$
where $m$ is the median of the maximum function under $\mu_{\alpha}^n$.
But from the extreme value theory (see \cite[Theorem 1.6.2, Corollary 1.6.3]{Leadbetter}), 
$$a_n\big(\max_{1 \leq i \leq n} x_i -b_n\big),$$
converges in law to the Gumbel distribution $G$, where 
$$ a_n \sim c_1 (\log n)^{\frac{1}{\alpha}-1}, \text{ and  } b_n \sim c_2(\log n)^{\frac{1}{\alpha}},$$
for some constant $c_1, c_2$. Moreover, as the Gumbel distribution has a right-tail behaving like $e^{-t}$, we see that the $B_{\ell^1}$ part in the enlargement of the deviations inequality of Proposition \ref{devexpo} is justified.
\end{Rem}

\begin{proof}[Proof of Proposition \ref{devexpo}]
Let $\Phi = \phi^{\otimes n} : \RR^n \to \RR^n$, defined by $\Phi(x) = (\phi(x_i))_{1\leq i \leq n}$, which sends $\mu_1^{n}$ to $\mu_{\alpha}^{ n}$.
Let  $r>0$, and $A$ be a measurable subset of $\RR_+^n$ such that $\mu_1^{n}(A)>0$. In a first step, we will use Lemma \ref{Brenieralpha} to see how the map $\Phi$ transform the set $A + \sqrt{r} B_{\ell^2}  + rB_{\ell^1}$. Actually, to transport the deviation inequality of $\mu_1^{n}$ it is sufficient to understand how $\Phi$ deforms $A' + \sqrt{r} B_{\ell^2}  + rB_{\ell^1}$ for a well-chosen subset $A'$ of $A$ such that $\mu_1^{ n}(A')>0$.
To this end, define
$$ B = \{ x \in \RR^n : || x||_{\infty} \leq  C \log n\}, \ A' = A\cap B,$$
where $C$ is some constant which will be chosen later.
Let $x\in A'$, $y\in B_{\ell^2}$, and $z \in B_{\ell^1}$. By Lemma \ref{Brenieralpha}, we have 
$$| \Phi(x+\sqrt{r} y)-\Phi(x)| \leq K \big( \sqrt{r} |y|+|x|^{\frac{1}{\alpha}-1}\sqrt{r}|y| +|\sqrt{r} y|^{\frac{1}{\alpha}} \big),$$
where the inequality has to be understood coordinate-wise, the functions being applied coordinate by coordinate to the vectors in $\RR^n$, and where $K$ is a  constant depending on $\alpha$ which will vary in the rest of the proof without changing name. Thus,
$$\Phi(x+\sqrt{r} y)-\Phi(x) \in K \Big(\sqrt{r} B_ {\ell^2} + (C \log n)^{\frac{1}{\alpha}-1} \sqrt{r} B_{\ell^2}  + r^{\frac{1}{2\alpha}} B_{\ell^{2\alpha}}\Big).$$
For $C \log n\geq 1$, we have
$$\Phi(x+\sqrt{r} y)-\Phi(x) \in K \Big( (C \log n)^{\frac{1}{\alpha}-1} \sqrt{r} B_{\ell^2}  + r^{\frac{1}{2\alpha}} B_{\ell^{2\alpha}}\Big).$$
Once again by Lemma \ref{Brenieralpha}, we get
$$| \Phi(x+\sqrt{r} y +r z)-\Phi(x+\sqrt{r}y)| \leq K\big( |rz| + |x+\sqrt{r}y|^{\frac{1}{\alpha}-1}|rz| + |rz|^{\frac{1}{\alpha}} \big),$$
where again this inequality is valid coordinate-wise. Using the convexity of the power function $t \mapsto |t|^{\frac{1}{\alpha}-1}$, or its sub-additivity, we get
$$| \Phi(x+\sqrt{r} y +r z)-\Phi(x+\sqrt{r}y)| \leq K\big( |rz|+(|x|^{\frac{1}{\alpha}-1}+ |\sqrt{r}y|^{\frac{1}{\alpha}-1})|rz| + |rz|^{\frac{1}{\alpha}} \big).$$
Note that Hölder's inequality implies 
$$ |y|^{\frac{1}{\alpha}-1} |z| \in B_{\ell^{\gamma}},$$
with $\frac{1}{\gamma} = \frac{1}{2}(\frac{1}{\alpha}+1)$. Thus,
$$ \Phi(x+\sqrt{r} y +r z)-\Phi(x+\sqrt{r}y) \in  K\big((C\log n)^{\frac{1}{\alpha}-1} rB_{\ell^1} +r^{\frac{1}{\gamma}}B_{\ell^{\gamma}} + r^{\frac{1}{\alpha}} B_{\ell^{\alpha}} \big).$$
Therefore,
$$ \Phi(x+\sqrt{r} y +r z) \in A +K\big((C\log n)^{\frac{1}{\alpha}-1}(\sqrt{r} B_{\ell^2}+ rB_{\ell^1}) +r^{\frac{1}{\gamma}}B_{\ell^{\gamma}} + r^{\frac{1}{\alpha}} B_{\ell^{\alpha}} + r^{\frac{1}{2\alpha}} B_{\ell^{2\alpha}} \big).$$
We now simplify the enlargement on the right-hand side. Observe that for any $0< a\leq b \leq c$,
$$ r^{1/b} B_{\ell^b} \subset r^{1/a} B_{\ell^a} + r^{1/c} B_{\ell^c}.$$
Indeed, if $x \in  r^{1/b} B_{\ell^b}$, then
$$\sum_{ |x_i|\geq 1} |x_i|^{a} \leq  \sum_{ |x_i|\geq 1} |x_i|^{b}\leq r,$$
and
$$\sum_{ |x_i|\leq 1} |x_i|^{c} \leq  \sum_{ |x_i|\leq 1} |x_i|^{b}\leq r.$$
Thus, $x = x\Car_{x\geq 1} + x\Car_{x< 1}$, with $x\Car_{|x|\geq 1} \in r^{1/a} B_{\ell^a}$ and $x\Car_{|x|< 1} \in r^{1/c} B_{\ell^c}$.
 Therefore, as $\alpha \leq 2\alpha \leq 2$,  $\alpha \leq \gamma \leq 2\alpha$, and $C \log n\geq 1$,
$$ \Phi(x+\sqrt{r} y +r z) \in A +K\big((C\log n)^{\frac{1}{\alpha}-1}(\sqrt{r} B_{\ell^2}+ rB_{\ell^1}) + r^{\frac{1}{\alpha}} B_{\ell^{\alpha}} \big).$$
Thus,
\begin{equation}
\label{transpenlar}
 \Phi\big(A+\sqrt{r} B_{\ell^2} +r B_{\ell^1}\big) \subset A +K\big((C\log n)^{\frac{1}{\alpha}-1}(\sqrt{r} B_{\ell^2}+ rB_{\ell^1}) + r^{\frac{1}{\alpha}} B_{\ell^{\alpha}} \big).
\end{equation}
Applying the deviation inequality \eqref{Taldev} of $\mu_1^{n}$, we get
$$ \mu_1^{n} \big( x \notin A' + \sqrt{r}B_{\ell^2} + rB_{\ell^1} \big) \leq \frac{e^{-Lr}}{ \mu_1^{n}(A')},$$
where $L>0$ is some constant independent of $n$.
But, since $$\int ||x||_{\infty} d\mu_1^{n}(x) \leq c_0 \log n,$$
for some numerical constant $c_0>0$, we have by Markov's inequality
$$\mu_1^{n}(A') \geq  \mu_1^{n}(A) - \mu_1^{ n}(B^c)\geq  \mu_1^{ n}(A) - \frac{c_0}{C}.$$
Thus,
$$ \mu_1^{n} \big( x \notin A' + \sqrt{r}B_{\ell^2} + rB_{\ell^1} \big) \leq \frac{e^{-cr}}{ \mu_1^{n}(A)-c_0/C}.$$
But, as $\mu_{\alpha}^{n} = \mu_1^{n} \circ \Phi^{-1}$, and $\Phi$ is a bijection,
\begin{align*}
\mu_1^{n} \big( x \notin A' + \sqrt{r}B_{\ell^2} + rB_{\ell^1} \big) & =
\mu_{\alpha}^{n} \big( \Phi(\RR_+^n\setminus ( A' + \sqrt{r}B_{\ell^2} + rB_{\ell^1}) \big)\\
&=  \mu_{\alpha}^{n} \big( \RR_+^n\setminus \Phi( A' + \sqrt{r}B_{\ell^2} + rB_{\ell^1}) \big).
\end{align*}
Using \eqref{transpenlar}, we deduce
$$\mu_{\alpha}^{n} \big( x \notin A +K\big((C\log n)^{\frac{1}{\alpha}-1}(\sqrt{r} B_{\ell^2}+ rB_{\ell^1}) + r^{\frac{1}{\alpha}} B_{\ell^{\alpha}} \big)\big) \leq \frac{e^{-cr}}{ \mu_1^{n}(A)-c_0/C}.$$
Adjusting the constant $c$ we get the claim.
\end{proof}

As observed in remark \ref{comprearrmap}, the 
monotone rearrangement $\psi$ of $\nu_1$ onto $\nu_{\alpha}$,  satisfies the same estimate of Lemma \ref{Brenieralpha} as $\phi$. 
Therefore, the same arguments as for the proof of Proposition \ref{devexpo} can be carried out, and yield a similar deviation inequality for $\nu_{\alpha}^n$ which we stated in Proposition \ref{devnualpha0}. 

In view of this deviation inequality for $\nu_{\alpha}^n$, we see that a Wigner matrix in the class $\mathcal{S}_{\alpha}$ when $\alpha\in (0,1)$ satisfies the concentration property $\mathcal{C}_{\alpha}$.

As for the case where $\alpha \in [1,2]$, the concentration property $\mathcal{C}_{\alpha}$ can be translated into a deviation inequality for Lipschitz or Hölder functions when $\alpha \in (0,1)$, as stated in the following lemma. 

\begin{Lem}\label{conclip0}
Let $\alpha \in (0,1)$. Assume $X$ satisfies the concentration property $\mathcal{C}_{\alpha}$ for some $\kappa>0$. Let $f : \mathcal{H}_n^{(\beta)} \to \RR$ be a function respectively $L_1$-Lipschitz and $L_2$-Lipschitz with respect to $||\ ||_{\ell^{1}}$, and $|| \ ||_{\ell^2}$. There is a constant $c>0$ depending on $\alpha$, such that if $f$ is moreover $L_{\alpha}$-Lipschitz with respect to $||\ ||_{\ell^{\alpha}}^{\alpha}$, then for any $t>0$,
$$ \PP\big( f(X) > m_f +t \big) \leq 4 \exp\Big( -c \min\Big( \frac{t^2}{\kappa^2(\log n)^{2(\frac{1}{\alpha}-1)}L_2}, \frac{t}{\kappa(\log n)^{\frac{1}{\alpha}-1}L_1 +\kappa L_{\alpha}} \Big) \Big),$$
whereas if 
 $$\forall A,B \in \mathcal{H}_n^{(\beta)}, \ f(A)- f(B) \leq L_{\alpha}'||A-B||_{\ell^{\alpha}},$$
for some $L_{\alpha}'>0$, then for any $t>0$,
$$ \PP\big( f(X) > m_f +t \big) \leq 4 \exp\Big( -c \min\Big( \frac{t^2}{\kappa^2(\log n)^{2(\frac{1}{\alpha}-1)}L_2}, \frac{t}{\kappa(\log n)^{\frac{1}{\alpha}-1}L_1 }, \frac{t^{\alpha}}{\kappa^{\alpha} L_{\alpha}'^{\alpha}} \Big) \Big),$$
where $m_f$ is the median of $f(X)$.
\end{Lem}

\subsection{Concentration inequalities for the largest eigenvalue}
We will prove in this section Proposition \ref{concvpintro}. We will see that it will fall easily form Weyl's inequality \cite[Theorem III.2.1]{Bhatia}, as it enables one to compute the Lipschitz constants of the largest eigenvalue function with respect to the distances $|| \ ||_{\ell^p}$ when $p \in [1,2]$ and $|| \ ||_{\ell^p}^p$ when $p\in (0,1)$ on $\mathcal{H}_n^{(\beta)}$.

\begin{proof}[Proof of Proposition \ref{concvpintro}]\label{proofconcvp}Let $\alpha \in (0,2]$. Let $X$ be a Wigner matrix satisfying the concentration property $\mathcal{C}_{\alpha}$ for some $\kappa>0$.
By Weyl's inequality \cite[Theorem III.2.1]{Bhatia}, the function $$f :Y\in \mathcal{H}_n^{(\beta)} \mapsto \lambda_{Y/\sqrt{n}}$$ is $n^{-1/2}$-Lipschitz with respect to the $p$-Schatten (pseudo-)norm $|| \ ||_{p}$ for any $p>0$, which is defined by 
\begin{equation}  \label{normSchatten}  \forall A \in \mathcal{H}_n^{(\beta)}, \ ||A||_p = \big( \tr |A|^p\big)^{1/p}.\end{equation}

Let $m_f$ denote the median of $f(X)$, and $t>0$. As  $\alpha \leq 2$, we have $ || \ ||_{\alpha}\leq || \ ||_{\ell^{\alpha}}$ by \cite[Theorem 3.32]{Zhan}. Thus, $f$ is also $n^{-1/2}$-Lipschitz with respect to $|| \ ||_{\ell^{\alpha}}$. Applying Lemmas \ref{conclip} and \ref{conclip0} successively to $f$ and $-f$, we deduce that for any $t>0$,
\begin{equation} \label{devmedianvp}\PP( |f - m_f|> t) \leq 8\exp\big(-c_{\alpha} h_{\alpha}(t)  \big),\end{equation}
with $h_{\alpha}$ defined in Proposition \ref{concvpintro}, and where $c_{\alpha}$ is some constant depending on $\alpha$.
Integrating the above inequality \eqref{devmedianvp}, we get
\begin{equation} \label{compmed} |\EE f(X) - m_f| =O(\kappa n^{-1/2} (\log n)^{\frac{1}{\alpha}-1}),\end{equation}
if $\alpha\in (0,1)$, and
$$|\EE f(X) - m_f| =O( \kappa n^{-1/2}),$$
if $\alpha \in [1,2]$, which gives the claim.

\end{proof}

\subsection{Two lemmas on spectral variation of Hermitian matrices}\label{spvarsection}
In view of Lemmas \ref{conclip} and \ref{conclip0}, proving the concentration inequalities of Propositions \ref{conclinearstat} and  \ref{concspintro} require to compute the Lipschitz constants of  the empirical spectral measure of Hermitian matrices, with respect to $|| \ ||_{\ell^p}$ when $p\in [1,2]$, and $|| \ ||_{\ell^p}^p$ when $p\in (0,1)$, and a well-chosen distance on $\mathcal{P}(\RR)$.

We will prove and discuss in this subsection Lemmas \ref{spvar1intro} and \ref{spvar0intro}.
For $p>0$, we denote by $\mathcal{W}_p$ the $L^p$-Wasserstein distance, defined for any probability measures $\mu$, $\nu$ on $\RR$ with finite $p^{\text{th}}$-moments by,
$$\mathcal{W}_p(\mu,\nu) = \Big(\inf_{\pi} \int |x-y|^p d\pi(x,y)\Big)^{1/p},$$
 if  $p\geq 1$ and by,
$$\mathcal{W}_p(\mu,\nu) =\inf_{\pi} \int |x-y|^p d\pi(x,y),$$
  if $p\in (0,1)$,
where the infimum is taken on all coupling $\pi$ between $\mu$ and $\nu$.

We begin with the proof of Lemma \ref{spvar1intro}.


\begin{proof}[Proof of Lemma \ref{spvar1intro}]

By Lidskii's theorem (see \cite[Corollary III 4.2]{Bhatia}), we have
$$ \lambda^{\downarrow}(A)-\lambda^{\downarrow}(B)\prec   \lambda^{\downarrow}(A-B),$$
where $ \lambda^{\downarrow}(A)$ denotes the vector of eigenvalues of $A$ in decreasing order, and $\prec$ the majorisation relation between vectors of $\RR^n$ (see \cite[Chapter II]{Bhatia} for a proper definition). Thus, by \cite[Theorem II.3.1]{Bhatia} we get, since $x \mapsto |x|^p$ is convex as $p\geq 1$,
$$ \tr | \lambda^{\downarrow}(A)-\lambda^{\downarrow}(B)|^p \leq \tr|\lambda^{\downarrow}(A-B)|^p.$$
Using the decreasing coupling between the spectra of $A$ and $B$, we get
\begin{equation} \label{Hoeffmanalphagene}\mathcal{W}_{p}( \mu_{A}, \mu_{B}) \leq \frac{1}{n^{1/p }} || A-B||_{p},\end{equation}
where $|| \ ||_p$ denotes the $p$-Schatten norm, defined in \eqref{normSchatten}.
But as $ p\leq 2$, we have by \cite[Theorem 3.32]{Zhan},
\begin{equation} \label{compnrom}  || A-B||_{p} \leq  || A-B||_{\ell^p},\end{equation}
which ends the proof of the first inequality of Lemma \ref{spvar1intro}.

As a consequence of the Kantorovitch-Rubinstein duality (see \cite[Particular case 5.16]{Villani}), we have
$$d \leq \mathcal{W}_1,$$
where $d$ is as in \eqref{defdStiel}. Besides, Jensen's inequality yields for any $p\geq 1$,
$$ \mathcal{W}_{1} \leq \mathcal{W}_p,$$
Therefore,
\begin{equation} \label{compdis} d \leq \mathcal{W}_1 \leq \mathcal{W}_p,\end{equation}
which gives the second claim of the lemma.
\end{proof}

\begin{Rem}\label{spvarrm}
When $p>2$, the inequality for $A,B \in \mathcal{H}_n^{(\beta)}$,
$$\mathcal{W}_p(\mu_A,\mu_B) \leq \frac{1}{n^{1/p}} || A-B||_{\ell^p},$$
 is no longer true, since for $B=0$ it amounts to \eqref{compnrom}, which is false when $p>2$, by taking $A=u u^*$, where $u$ is the constant vector. 

When $p<1$, one may hope for the inequality
\begin{equation} \label{hope} \mathcal{W}_{p}(\mu_A,\mu_B)\leq \frac{1}{n} || A-B||_{\ell^p}^p,\end{equation}
to hold. But taking formally $p \to 0$, would yield
\begin{equation} \label{hope2} | \lambda(A) \Delta \lambda(B) | \leq | (i,j) : A_{i,j} \neq B_{i,j}|,\end{equation}
where $\lambda(A),\lambda(B)$ denote the set of eigenvalues of $A$ and $B$. But one can see that changing $1$ entry to a matrix can change the whole spectrum, which disproves \eqref{hope2}. 


\end{Rem}

The moral of remark \ref{spvarrm} is that one cannot have \eqref{hope} with a constant $1$ on the right-hand side. As the cost function $| \ |^p$ behaves quite badly when $p<1$ as it is not convex (see \cite{McCann} for this transportation problem with concave costs), in particular, the optimal transport map is not necessarily the monotone rearrangement contrary to the case $p\geq 1$, we will not investigate further the question of having a spectral variation inequality involving the $L^p$-Wasserstein distance. We prefer to deal with another distance on $\mathcal{P}_p(\RR)$, the set probability measures on $\RR$ with finite $p^{\text{th}}$ moments, which induces the same topology as $\mathcal{W}_p$ and dominates $d$.  This distance is chosen so that, applied to empirical spectral measures, it will be controlled by $|| \ ||_{\ell^p}^p$ in the case where $p \in (0,1)$.

To this end, let $p\in (0,1)$ and define for any $\mu, \nu  \in \mathcal{P}_p(\RR)$, 
\begin{equation} \label{distdp} d_{p}(\mu,\nu) = \sup_{t\in \RR }\Big| \int (t-x)^{p}_+ d\mu(x) - \int (t-x)^{p}_+ d\nu(x) \Big|.\end{equation}
Taking formally $p$ to $0$, we retrieve the Kolmogorov-Smirnov distance $d_{KS}$. Recall that by integrating by parts, we can write
$$d_{KS}(\mu,\nu) = \sup \big\{ \big| \int f d\mu - \int f d \nu\big| : f \in \text{NBV},\  || f||_{BV} \leq 1\big\} ,$$
where $\text{NBV}$ denotes the set of normalized functions with bounded variations, that is, functions which are the integrals of finite signed measures, and 
$$ || f||_{BV} = | \sigma|(\RR),$$
whenever $f$ is the distribution function of the finite signed measure $\sigma$, and $|\sigma|$ is the total variation of $\sigma$. 

We can actually have a similar formulation for $d_p$, by introducing the fractional integrals of order $p+1$ on $\mathcal{M}_s^{p}$, the set of finite signed measures $\sigma$ such that $|\sigma|$ has a finite $p^{\text{th}}$-moment, which we defined in \eqref{defintfrac}. 
We recall that fractional integrals enjoy the following integration by parts formula (see \cite[(5.16)]{Samko}): for $\mu,\nu \in \mathcal{M}_s^{p}$,
\begin{equation} \label{IPP} \int (\mathcal{I}_+^{p+1} \mu)(t) d\nu(t) = \int (\mathcal{I}_-^{p+1}\nu)(x) d\mu(x).\end{equation}
Thus, we can write
\begin{align}
 d_{p}(\mu,\nu) &= \Gamma(p+1) \sup_{t\in \RR }\big| (\mathcal{I}_+^{p+1} \mu)(t)  - (\mathcal{I}_+^{p+1} \mu)(t) \big| \nonumber\\
&= \Gamma(p+1)\sup_{\sigma} \big| \int (\mathcal{I}_-^{p+1}\sigma) d\mu - \int (\mathcal{I}_-^{p+1}\sigma) d\nu \big|,\label{formulation}
\end{align}
where the supremum is taken on all  $\sigma \in \mathcal{M}_s^p$, such that $|\sigma|(\RR) \leq 1$.
The inequality $d_p \geq$ \eqref{formulation} is the consequence of the integration by parts formula \eqref{IPP}, whereas the equality is given by taking $\sigma = \delta_t$, for $t\in \RR$. 
We investigate now the link between the distances $d$, defined in \eqref{defdStiel}, $\mathcal{W}_p$ and $d_p$ when $p \in (0,1)$.
\begin{Pro}\label{compdist}
Let $p\in (0,1)$. Then, $d_{p}$, defined in \eqref{distdp}, is a distance on $\mathcal{P}_p(\RR)$, and metrizes the weak topology. More precisely, there is a constant $C_p>0$ such that
\begin{equation} \label{compdistdp} d(\mu,\nu) \leq C_{p} d_{p} (\mu,\nu),\end{equation}
for all $\mu,\nu \in \mathcal{P}_p(\RR)$. One can choose 
\begin{equation} \label{defCp} C_p = \sqrt{\pi}(p+1) \frac{\Gamma\big( \frac{p+1}{2}\big)}{\Gamma\big( 1+\frac{p}{2}\big)}.\end{equation}
 Furthermore,
\begin{equation} \label{compdistWp} d_p \leq \mathcal{W}_p.\end{equation}
\end{Pro}

\begin{Rem}
We actually do not know if the distances $d_{p}$ and $\mathcal{W}_{p}$ are comparable, meaning that the reversed inequality $d_p \geq K_p \mathcal{W}_{p}$ is true for some $K_p>0$. We do know however, by the remark \ref{spvarrm}, that such an inequality cannot hold with some constant $K_p$ staying bounded when $p\to 0$.

\end{Rem}

\begin{proof}
In view of the formulation of $d_p$ as \eqref{formulation}, the stake behind \eqref{compdistdp} is to represent the function $t \mapsto (z-t)^{-1}$ as the fractional integral of order $p+1$ of some function. The constant $C_p$ will arise as a bound on the $L^1$ norm of this function as $\Im z \geq 1$, over $\Gamma(p+1)$.

The fractional integral of order $p+1$ of the function $t \mapsto (z-t)^{-1}$ is given in \cite{Samko}, which we state in the next lemma.
\begin{Lem}[{\cite[Chapter 2 (5.25)]{Samko}}]\label{stiel}Let $p\in (0,1)$.
For any $z \in \CC$, $\Im z>0$, we have
$$\forall x \in \RR, \  \frac{1}{z-x} =  \mathcal{I}_-^{p+1}(h)(x),$$
with
\begin{equation} \label{defphi} \forall t \in \RR, \  h(t) = e^{i \pi (p+1)} \Gamma(p+2) \frac{1}{(z-t)^{p+2} },\end{equation}
where $\zeta^{p}$ is the principal branch of the $\alpha^{\text{th}}$-root on $\CC\setminus \RR_-$.
\end{Lem}


Let $\Im z \geq 1$ and $h$ as in \eqref{defphi}. We have
\begin{equation*} \frac{1}{\Gamma(p+1)} ||h ||_1\leq (p+1) \int _{-\infty}^{+\infty}\frac{dt}{(1+ t^2)^{1+p/2}} = 2(p+1) \int_0^{+\infty}\frac{dt}{(1+ t^2)^{1+p/2}} := C_p,\end{equation*}
where we used $\Gamma(p+2) = (p+1)\Gamma(p+1)$.
Therefore,
$$ d  \leq C_{p} d_{p}.$$
But, one can recognize an Euler integral of the first kind in the definition of $C_p$, by making successively the changes of variables $t = \tan u$, and $v = (\cos u)^2$, which yields,
$$C_p = (p+1) \int_0^{1} v^{\frac{p-1}{2}} (1-v)^{-\frac{1}{2}} dv.$$
Therefore by \cite[(2.13)]{Artin}, we deduce the value for $C_p$ claimed in \eqref{defCp}.

Inequality \eqref{compdistWp} is the consequence of the sub-additivity of the function $x \mapsto x^p$ on $\RR^+$. More precisely, for any $x,y,t \in \RR$,
$$(t-x)_+^p - (t-y)_+^p \leq |x-y|^p.$$
Integrating the above inequality under a coupling $P$ of two probability measures with finite $p^{\text{th}}$-moment yields the claim.

From \eqref{compdistdp}, we deduce that the topology induced by $d_p$ on $\mathcal{P}_p(\RR)$ is finer than the weak topology, and by \eqref{compdistWp} that it is coarser than the one induced by $\mathcal{W}_p$. But $\mathcal{W}_p$ induces the weak topology on $\mathcal{P}_p(\RR)$ by \cite[Theorem 6.9]{Villani} (as $| \ |^p$ is a metric on $\RR$ for $p\leq1$), therefore $d_p$ induces the weak topology on this set.
\end{proof}

%

We finally prove that the distance $d_p$ we introduced, when applied to spectral measures of Hermitian matrices, is dominated by $|| \ ||_{\ell^p}^p$ for $p\in (0,1)$, this will directly imply the result of Lemma \ref{spvar0intro}. 

\begin{Lem}\label{spvar0}Let $p\in (0,1)$.
Let $A,B \in \mathcal{H}_n^{(\beta)}$. 
\begin{equation} \label{ineqRot}d_{p}(\mu_{A},\mu_B)  \leq \frac{1}{n} || A-B||_{\ell^p}^p,\end{equation}
where $d_p$ is defined in \eqref{distdp}.
In particular,
\begin{equation} \label{controlStielp} d(\mu_{A},\mu_B) \leq \frac{C_p}{n} ||A-B||_{\ell^p}^p,\end{equation}
where $C_p$ is as in \eqref{defCp}.

\end{Lem}
\begin{Rem}\label{autredis}
Defining the distance 
$$d_p^{-}(\mu,\nu) = \sup_{t \in \RR} \big| \int (t-x)_-^{p} d\mu(x) - \int (t-x)_-^{p} d\nu(x)\big|,$$
for any $\mu,\nu \in \mathcal{P}_p(\RR)$, we see that we have a similar representation as for $d_p$, that is,
$$d_p^-(\mu,\nu) = \sup_{\sigma} \big| \int \mathcal{I}_+^{p+1}(\sigma) d\mu - \int \mathcal{I}_+^{p+1}(\sigma) d\nu \big|,$$
where $\sigma$ run in $\mathcal{M}_s^p$ such that $|\sigma|(\RR)\leq 1$. Moreover, we clearly get the same inequality as \eqref{ineqRot} for $d_p$.
\end{Rem}

\begin{proof}
As $\alpha \leq 2$, the second inequality of \eqref{ineqRot} is due to \cite[Theorem 3.32]{Zhan}. 
To prove the first inequality, we begin by recalling an inequality due to Rotfel'd originally, and then to Thompson \cite{Thompson} (for an extension and a simpler proof). Let $F: \RR_+^{2n} \to \RR$ be a concave symmetric function. Then for any $A,B \in \mathcal{H}_n^{(\beta)}$ positive semi-definite,
$$ F(\lambda(A+B), 0) \leq F(\lambda(A),\lambda(B)),$$
where $\lambda(C)$ denotes the vector of eigenvalues of a Hermitian matrix $C$. Note that since $F$ is symmetric, there is no ambiguity in the writing. 
Let $t\in \RR$. We have,
$$t-A-B \leq (t-A)_+ + |B|.$$
In particular, if we denote $\lambda_1(C) \geq\lambda_2(C)\geq  ... \geq \lambda_n(C)$ the eigenvalues of some Hermitian matrix $C$, then by Weyl's inequality \cite[Theorem III.2.1]{Bhatia}, for any $i\in \{1,...,n\}$, 
$$\lambda_i(t-A-B) \leq \lambda_i\big((t-A)_+ + |B|\big).$$
Therefore,
$$\lambda_i(t-A-B)_+ \leq \lambda_i\big((t-A)_+ + |B|\big).$$
Define
$$\forall x \in \RR_+^{2n}, \ F(x) = \sum_{i=1}^{2n} x_i^{\alpha}.$$
Since $A,B$ are Hermitian,
$$ \lambda (t-A-B)_+ =(t-\lambda(A+B))_+.$$
As $F$ is non-decreasing coordinate-wise,
$$F\big( (t-\lambda(A+B))_+,0\big) \leq F\big( \lambda((t-A)_+ + |B| ),0\big).$$
Rotfel'd inequality gives
$$F\big( \lambda\big((t-A)_+ + |B| \big),0\big) \leq F\big( \big(t-\lambda(A)\big)_+, |\lambda(B)| \big).$$
Thus,
$$\sum_{i=1}^n \big(t-\lambda_i(A+B)\big)_+^{\alpha} \leq \sum_{i=1}^n \big(t-\lambda_i(A)\big)_+^{\alpha} + \sum_{i=1}^n |\lambda_i(B)|^{\alpha}.$$
Applying this inequality with $A+B$, $-B$ instead of $A$ and $B$, we get the first claim.
The inequality \eqref{ineqRot} is a just reformulation of the above inequality and a use of the comparison \eqref{compnrom} between $\ell^p$-(quasi)-norm and $p$-Schatten (quasi)-norm.
Finally, using Proposition \ref{compdist}, we deduce that \eqref{controlStielp} is true.
\end{proof}

With the Lemmas \ref{spvar0} and \ref{spvar1intro}, we can now give a proof of Propositions \ref{conclinearstat} and \ref{concspintro}.

\begin{proof}[Proof of Proposition \ref{conclinearstat}]
Let $\alpha \in (0,2]$ and $X$ to be a Wigner matrix satisfying the concentration property $\mathcal{C}_{\alpha}$ with some $\kappa>0$. Lemma \ref{spvar1intro} and Hölder's inequality allow us to say that if $f : \RR \to \RR$ is $1$-Lipschitz, then the function
\begin{equation} \label{linearstat}Y \in \mathcal{H}_n^{(\beta)} \mapsto \frac{1}{n} \sum_{i=1}^n f(\lambda_i(Y/\sqrt{n})),\end{equation}
where $\lambda_1(Y),...,\lambda_n(Y)$ denote the eigenvalues of $Y$, is $n^{-\frac{1}{2}-\frac{1}{p}}$-Lipschitz with respect to $|| \ ||_{\ell^p}$ for any $p\in [1,2]$. Thus, using Lemma \ref{conclip}, we deduce the concentration inequality for the linear statistics of Lipschitz functions of Proposition \ref{conclinearstat} in the case $\alpha \in [1,2]$.

Assume now that $\alpha \in (0,1)$ and $f$ is $1$-Lipschitz and moreover can be written $f = \mathcal{I}_{\pm}^{1+\alpha}(\sigma)$ for some $\sigma \in \mathcal{M}_s^{\alpha}$ such that $|\sigma|(\RR)\leq m$, then by Lemma \ref{spvar0} (and remark \ref{autredis}), we know that the map \eqref{linearstat} is $\Gamma(\alpha+1)^{-1}n^{-1-\frac{\alpha}{2}}m$-Lipschitz with respect to $|| \ ||_{\ell^{\alpha}}^{\alpha}$. Thus we can deduce from Lemma \ref{conclip0} the second concentration inequality of Proposition \ref{conclinearstat}.

\end{proof}

We prove now Proposition \ref{concspintro}. 
\begin{proof}[Proof of Proposition \ref{concspintro}]
Fix some $z \in \mathcal{K}$. Let $f_z$ denote the function on $\mathcal{H}_n^{(\beta)}$ defined by,
$$ \forall Y \in \mathcal{H}_n^{(\beta)}, \ f_z(Y) = g_{ \mu_{Y/\sqrt{n}}}(z).$$

As $\Im z\geq 1$, we see that the function $x\mapsto (z-x)^{-1}$ is $1$-Lipschitz. Moreover, we know by Lemma \ref{stiel} that when $\alpha \in (0,1)$,
$$ \frac{1}{z-x} =\mathcal{I}_-^{\alpha+1}(h),$$ 
with $||h||_{\ell^1} \leq \Gamma(\alpha+1)C_{\alpha}$, where $C_{\alpha}$ is as in \eqref{defCp}.
Let $m_z$ be the median of $f_z(X)$. Let also $r>0$.  We deduce by Proposition \ref{conclinearstat}, and using remark \ref{concsimpl0} in the case $\alpha \in (0,1)$, that there is a constant $c_{\alpha}$ depending on $\alpha$ such that,
\begin{equation} \label{concspect}\PP( |f_z - m_z| > t) \leq 8\exp ( -c_{\alpha}k_{\alpha}(t)),\end{equation}
where $k_{\alpha}$ is defined in the statement of Proposition \ref{concspintro}.
Integrating this inequality, we get
$$|\EE f_z(X) - m_z| \leq \eps_n $$
with $\eps_n = O( \kappa n^{-1} (\log n)^{(\frac{1}{\alpha}-1)_+})$, uniformly in $z\in \CC$, $\Im z\geq 1$. With this notation, we get for any $t>0$,
$$\PP( | f_z-\EE f_z| > t + \eps_n)\leq 8\exp\big ( -c_{\alpha} k_{\alpha}(t)\big).$$
Let $\mathcal{N}_{t}$ be a $t$-net of $\mathcal{K}$. As $z \mapsto f_z(X)$ is $1$-Lipschitz on $\{ z \in \CC: \Im z\geq 1\}$, we have
$$ \PP\big( \sup_{z\in \mathcal{K}} |f_z - \EE f_z | >2 t +\eps_n \big) \leq 
8 |\mathcal{N}_{t}|\exp\big ( -c_{\alpha} k_{\alpha}(t)\big),$$
As $\mathcal{K}$ is a subset of $\CC$ of diameter inferior to $1$, we can find a $t$-net $\mathcal{N}_t$ such that $|\mathcal{N}_t|\leq t^{-2}$.
Thus,
$$ \PP\big( d( \mu_{X/\sqrt{n}}, \EE \mu_{X/\sqrt{n} }) >2 t +\eps_n \big) \leq 
\frac{8}{t^2}\exp\big ( -c_{\alpha} k_{\alpha}(t)\big),$$
which, adjusting the constant $c_{\alpha}$, gives the claim.

\end{proof}

\section{Deterministic equivalents for  Wigner matrices}\label{secdetermexpoWigner}

 We will prove in this section some uniform deterministic equivalents for the spectral measure and largest eigenvalue of deformed Wigner matrices having concentration $\mathcal{C}_{\alpha}$ for $\alpha \in (0,2)$ (see definition \ref{defCalpha}), using the inequalities proved in the preceding section. We will also prove a deterministic equivalent for traces of polynomials of deformed Wigner matrices, but which will not rely on concentration arguments. In particular, these deterministic equivalents will entail that assumption $(i)$ of Theorem \ref{theoremgene} holds for the spectral measure, the largest eigenvalue and the traces of polynomials of Wigner matrices in $\mathcal{S}_{\alpha}$. More precisely, we will prove the following propositions.

\begin{Pro} \label{chaining0} Let $\alpha \in (0,2)$. Let $X$ be a Wigner matrix such that $\EE|X_{1,2} -\EE X_{1,2}|^2=1$ and satisfying the concentration property $\mathcal{C}_{\alpha}$.
For any $r>0$, 
$$\sup_{H \in r n^{1/\alpha}B_{\ell^{\alpha}}} d\big( \mu_{X/\sqrt{n} +H}, \mu_{sc}\boxplus \mu_{H} \big) \underset{n\to +\infty}{\longrightarrow} 0,$$
in probability, where $d$ is the distance defined in \eqref{defdStiel} .

\end{Pro}
\begin{Rem}
This statement fails when $\alpha=2$ since $X/\sqrt{n}$ is in $r n^{1/2} B_{\ell^2}$ for some $r>0$, with positive probability uniform in $n$. Whereas on one hand, by Wigner's theorem (see \cite{Guionnet})
$$\mu_{2X/\sqrt{n}} \underset{n\to +\infty}{\leadsto } \mu_{sc,2},$$
in probability, where for any $a>0$,
$$\mu_{sc,a} = \frac{1}{2a^2 \pi} \sqrt{4a^2-x^2} \Car_{|x|\leq 2a}dx.$$
On the other hand, by continuity of the free convolution (see \cite[Proposition 4.13]{BV}), 
$$ \mu_{sc}\boxplus \mu_{X/\sqrt{n}} \underset{n\to +\infty}{\leadsto} \mu_{sc}\boxplus \mu_{sc},$$
in probability, and we have $\mu_{sc}\boxplus \mu_{sc}= \mu_{sc,\sqrt{2}}$ by \cite[Example 5.3.26]{Guionnet}.

\end{Rem}

\begin{Pro}
\label{chaining1}
 Let $\alpha \in (0,2)$. Let $X$ be a centered Wigner matrix satisfying the concentration property $\mathcal{C}_{\alpha}$ such that $\EE|X_{1,2}|^2=1$.  Define the function $\rho$ by,
\begin{equation} \label{deff} \forall x \in \RR, \  \rho(x) = \begin{cases} x +\frac{1}{x} & \text{ if } x \geq 1,\\
2 & \text{ otherwise.}
\end{cases}
\end{equation}
For any $r>0$,
$$\sup_{A \in r B_{\ell^{\alpha}}} \big| \lambda_{X/\sqrt{n} + A}-\rho(\lambda_A) \big| \underset{n\to +\infty}{\longrightarrow} 0,$$
in probability.
\end{Pro}

For the traces of polynomials of independent Wigner matrices we will prove the next proposition.
\begin{Pro}\label{convunifpoly}Let $\alpha \in (0,2]$. Let $P\in \CC \langle \textbf{X}\rangle$ be a non-commutative polynomial of total degree $d>\alpha$. Let $\textbf{X} = (X_1,...,X_p)$ be a family of independent centered Wigner matrices with entries having finite $(d+1)^{\text{th}}$-moments, such that $\EE|M_{1,2}|^2 = 1$ for any $M\in \{X_1,...,X_p\}$. For any $r>0$,
$$\sup_{\textbf{H}\in r B_{\ell^{\alpha}}} \big| \tau_n[P(\textbf{X}/\sqrt{n}+n^{1/d}\textbf{H})] - \tau[P(\textbf{s})]-\tr[P_d(\textbf{H})] \big| \underset{n\to+\infty}{\longrightarrow} 0,$$
in probability, where  $P_d$ is the homogeneous part of degree $d$ of $P$,  $\textbf{s} =(s_1,...,s_p)$ is a free family of $p$ semi-circular variables in a non-commutative probability space $(\mathcal{A},\tau)$ and,
$$ B_{\ell^{\alpha}} = \big\{ \textbf{H} \in (\mathcal{H}_n^{(\beta)})^p : \sum_{i=1}^p \tr|H_i|^{\alpha}\leq 1\big\}.$$
\end{Pro}

 It is interesting to note that we are able for polynomials, to make the approximation hold uniformly in $H\in r B_{\ell^2}$, which is why we can consider the Gaussian case in our large deviations principle of Theorem \ref{LDPpoly}.

\subsection{Deterministic equivalents in expectation}

Our approach to prove Propositions \ref{chaining0} and \ref{chaining1} consists is showing in a first step the proposed uniform deterministic equivalents in expectation, and then make use the concentration inequalities of the last section \ref{Chapconc} together with a chaining argument to show that these equivalent hold uniformly in probability.

For the empirical spectral measure, we have such a uniform deterministic equivalents in expectation by the following result of Bordenave and Caputo \cite{Bordenave}.

\begin{The}[{\cite[Theorem 2.6]{Bordenave}}]\label{citeBC}
Let $X$ be a Wigner matrix such that $\EE|X_{1,2}-\EE X_{1,2}|^2 = 1$, $\EE |X_{1,2}|^3<+\infty$, and $\EE X_{1,1}^2<+\infty$. There exists a universal constant $c>0$ such that for any $H\in \mathcal{H}_n^{(\beta)}$,
$$\delta( \EE\mu_{X/\sqrt{n}+H}, \mu_{sc}\boxplus \mu_H) \leq c \frac{\sqrt{\EE X_{1,1}^2} + \EE|X_{1,2}|^3}{\sqrt{n}},$$ 
where $\delta$ is defined for any $\mu,\nu\in \mathcal{P}(\RR)$,
$$ \delta(\mu,\nu) = \sup \big\{\big| g_{\mu}(z) - g_{\nu}(z)\big| : \Im z\geq 2\big\},$$
where $g_{\mu}$ and $g_{\nu}$ denote the Stieltjes transforms of $\mu$ and $\nu$.
\end{The}

For the largest eigenvalue, we will prove the following proposition.

\begin{Pro}\label{compexpvp}Let $\alpha \in (0,2)$.
Let $X$ be a centered Wigner matrix such that $\EE|X_{1,2}|^2=1$ and $\EE| X_{1,1}|^4, \EE|X_{1,2}|^4<+\infty$. For any $r>0$,
$$ \sup_{H \in r B_{\ell^{\alpha}}} | \EE\lambda_{X/\sqrt{n} + H} - \rho(\lambda_H)| \underset{n\to +\infty}{\longrightarrow} 0,$$
where $\rho$ is the function defined in \eqref{deff}.
\end{Pro}

%

\begin{proof}
In a first step, we will perfom a truncation and convolution argument as to the one used in \cite[Proposition 4.1, step 1]{BordCap}, in order to reduce the problem to the case the entries of $X$ satisfies a Poincaré inequality.  Let $\eps>0$ and let $G$ be a GUE matrix, that is, $G = \frac{1}{\sqrt{2}}(B+B^*)$ where $B$ is a matrix with i.i.d complex Gaussian entries with covariance $\frac{1}{2}I_2$, independent from $X$. We set $X^{(\eps)}$ to be the Hermitian matrix with $(i,j)$-entry,
$$ X_{i,j}^{(\eps)} = \frac{X_{i,j}\Car_{|X_{i,j}|\leq \eps^{-1}} - \EE X_{i,j}\Car_{|X_{i,j}|\leq \eps^{-1}} }{(\Var (X_{i,j} \Car_{|X_{i,j}|\leq \eps^{-1}})^{1/2}},$$
and $Y^{(\eps)} = (1+\eps^2)^{-1/2}(X^{(\eps)} + \eps G)$. By \cite[Theorem 1.2]{BGMZ}, $Y^{(\eps)}$  has entries satisfying a Poincaré inequality .

We know by \cite[Theorem 2]{Latala} that there is some constant $C>0$ such that for any centered Wigner matrix $H$,
$$\EE|| H || \leq C \Big( \max_i\big(  \sum_j \EE |H_{i,j}|^2 \big)^{\frac{1}{2}} + \big (\sum_{i,j} \EE |H_{i,j}|^4\big)^{\frac{1}{4}}\Big).$$
This inequality yields as the entries of $X$ have finite fourth moments,
$$ \lim_{\eps \to 0} \limsup_{ n\to\infty} \EE || X - Y^{(\eps)}|| =0.$$
But, using Weyl's inequality \cite[Theorem III.2.1]{Bhatia}, and the fact that $\rho$ is $1$-Lipschitz, we see that $A\mapsto |\lambda_{X/\sqrt{n}+A}-\rho(\lambda_A)|$ is $2$-Lipschitz with respect to $|| \ ||_{\ell^2}$. Thus, we can focus on proving Proposition \ref{compexpvp} when $X$ has entries satisfying a Poincaré inequality.
We make now another reduction of the statement to a convergence in probability and to the case where the supremum is taken on the set of matrices which we denote by $\mathcal{G}$, consisting of $m$-sparse matrices $A$ (meaning at most $m$ entries are non-zero) with spectral radius bounded by $r$, for some fixed $r,m>0$. 

Note that by Weyl's inequality   and  \eqref{compnrom}, we have for any $A\in rB_{\ell^{\alpha}}$,
\begin{equation} \label{inequnif}|\lambda_{X/\sqrt{n}+A}|\leq r + ||X/\sqrt{n}||, \quad 2\leq \rho(\lambda_A)\leq \rho(r).\end{equation}
As $|| X/\sqrt{n}||$ converges in $L^2$ by \cite[Theorem 2.1.22, 27]{Guionnet}, we deduce that, uniformly in $A\in rB_{\ell^{\alpha}} $, $|\lambda_{X/\sqrt{n}+A}-\rho(\lambda_A)|$ is uniformly integrable. Therefore it suffices to prove that for any $t>0$,
$$ \sup_{A\in rB_{\ell^{\alpha}}} \PP( | \lambda_{X/\sqrt{n}+A} - \rho(\lambda_A)|>t)\underset{n\to +\infty}{\longrightarrow} 0.$$

 Let $A\in rB_{\ell^{\alpha}}$, and $M_1\geq ...\geq M_{n^2}$ be the values $|A_{i,j}|$ in non-increasing order. We have,
\begin{equation} \label{estimcoeffball}\forall k \in \{1,...,n^2\}, \  M_k \leq r^{1/\alpha} k^{-1/\alpha}.\end{equation}
Let now $m\in \NN$ and $v_1,...,v_m$ the locations of the $m$ largest values of $|A_v|$. Define $A^{(m)}$ to be the matrix,
$$\forall v \in \{1,...,n\}^2, \ A_v^{(m)} = \sum_{i=1}^m A_{v_i} \delta_{v_i,v}.$$ 
As $\alpha <2$, we deduce,
$$ || A-A^{(m)}||_{\ell^2}^2 \leq r^{\frac{2}{\alpha}} \sum_{k>m} k^{-\frac{2}{\alpha}} =  O(  m^{1- \frac{2}{\alpha}}).$$
Thus, again by Weyl's inequality, it is sufficient to prove for any fixed $m\in \NN$, $r>0$, and $t>0$,
$$ \sup_{A\in \mathcal{G}} \PP( | \lambda_{X/\sqrt{n}+A} - \rho(\lambda_A)|>t)\underset{n\to +\infty}{\longrightarrow} 0.$$
To prove this claim, we will follow a rather classical argument  relying on the Frobenius formula used in the study of finite rank perturbations as in \cite{Benaych} for example, to determine the behavior of the largest eigenvalue of deformed models.

Diagonalize $A=UDU^*$, with $U$ of size $n\times m$ such that $U^* U=I_m$. By Frobenius formula (see \cite[section 4.1]{Benaych}), $\lambda_{X/\sqrt{n}+A}$ is either in the spectrum of $X/\sqrt{n}$, denoted $\sigma(X/\sqrt{n})$, or the largest zero of the function,
\begin{equation} \label{deffn}\forall x\notin \sigma(X/\sqrt{n}), \ f_{n,A}(x) = \det(I _m- U^*R(x)U D).\end{equation}

Our main task consists in proving that this function is uniformly close on any compact subset of $\{ \Re z> \lambda_{X/\sqrt{n}} \}$ to the following deterministic limit function, 
\begin{equation} \label{deffa}\forall \Re z>2, \ f_A(z) = \det(I_m- g_{\mu_{sc}}(z) D).\end{equation}
\begin{Lem}\label{compeq}
Let $\delta>0$ and define,
$$W_{\delta} = \{ \lambda_{X/\sqrt{n}} \leq  2 +\delta\}.$$
For all subset $\Omega$ compactly included in $\{ z \in \CC : \Re z > 2+\delta\}$ and $t>0$,
$$\sup_{A\in \mathcal{G}} \PP( \{\sup_{z \in \Omega} |f_{n,A}(z) - f_A(z)| >t\}\cap W_{\delta} ) \underset{n\to+\infty}{\longrightarrow} 0,$$
where $f_{n,A}$ and $f_A$ are defined in \eqref{deffn}, \eqref{deffa}.
\end{Lem}

Assume  for the moment that this lemma is true. Note that the functions $f_A$, $A\in \mathcal{G}$, form a normal family of holomorphic functions on $\{ z \in \CC : \Im z >2\}$. By \cite[Chapter 5, Theorem 2]{Ahlfors}, it is thus a pre-compact family in the space of holomorphic functions on $\{ z \in \CC : \Im z >2\}$. We deduce by Hurwitz's theorem \cite[Chapter 5, Theorem 10]{Ahlfors} that for any $\delta>0$ and $\Omega$ open subset compactly included in  $\{ z : \Re z>2\}$, there is some $t>0$ such that for any holomorphic function $g$ defined on a neighborhood of $\Omega$, and $A\in rB_{\ell^{\alpha}}$ such that $\sup_{\Omega} || f_A - g|| <t$, then either $f_A$ does not have any zeros in $\Omega$ and therefore $g$ neither, or for any zeros of $f_A$ in $\Omega$, corresponds a zero of $g$ in $\Omega$ which is $\delta$-close.

Let $\delta, r>0$. We set
$$V_{\delta,r} = \big\{ \lambda_{X/\sqrt{n}} \leq 2+\delta, \ 2-\delta \leq \lambda_{X/\sqrt{n}+A} \leq \rho(r),  \}.$$
Let also $\Omega$ be some open subset compactly included in $\{ z : \Re z>2+\delta\}$ such that $[2+2\delta, \rho(r) ] \subset \Omega$. We deduce that for any $\delta>0$ there is a $t>0$, such that,
$$ \PP(  \{ | \lambda_{X/\sqrt{n}+A} - \rho(\lambda_A)|>3\delta\} \cap V_{\delta,r})\leq \PP( \{ \sup_{z \in \Omega} | f_{n,A}(z) - f_A(z)| >t \}\cap W_{\delta}).$$
As this $t$ does not depend on $A\in \mathcal{G}$, we get from Lemma \ref{compeq}
$$ \sup_{A\in \mathcal{G}}  \PP(  \{ | \lambda_{X/\sqrt{n}+A} - \rho(\lambda_A)|>3\delta\} \cap V_{\delta,r}) \underset{n\to+\infty}{\longrightarrow} 0.$$
It remains to show that $\PP(V_{\delta,r})$ goes to $0$ as $n\to +\infty$ uniformly in $A\in rB_{\ell^{\alpha}}$. Note that almost surely (taking an arbitrary coupling of the matrices $X$), we have by Hoeffman-Weilandt inequality \eqref{Hoeffmanalphagene},
$$\sup_{A\in rB_{\ell^2}} \mathcal{W}_2(\mu_{X/\sqrt{n}+A},\mu_{X/\sqrt{n}}) \underset{n\to +\infty}{\longrightarrow} 0.$$
Thus, by Wigner's theorem, almost surely, $\mu_{X/\sqrt{n}+A}$ converges weakly towards $\mu_{sc}$ uniformly in $A\in rB_{\ell^2}$. By lower-semicontinuity of the map 
$$ \mu \in \mathcal{P}(\RR) \mapsto \sup \mathrm{supp}(\mu),$$
we deduce that 
$$ \liminf_{n\to +\infty} \sup_{A\in rB_{\ell^2}}(\lambda_{X/\sqrt{n}+A}-2) \geq 0,$$
almost surely. Using the above convergence, \eqref{inequnif} and  the convergence of the largest eigenvalue of $X/\sqrt{n}$ to $2$  in probability, we can conclude that 
$$ \limsup_{n\to +\infty} \sup_{A\in \mathcal{G}} \PP( V_{\delta,r}) \underset{r\to +\infty}{\longrightarrow} 0,$$
which gives the claim of Proposition \ref{compexpvp}. Thus, we are reduced to show Lemma \ref{compeq}.

\end{proof}

\begin{proof}[Proof of Lemma \ref{compeq}]
Let $\delta>0$ and $\Omega$ as in the statement of Lemma \ref{compeq}. Let $\eta = \inf \{ \Re z : z \in \Omega \} -2$ and $\zeta$ a Lipschitz function such that 
\begin{equation} \label{defzeta} \Car_{(-\infty, 2+\delta]} \leq \zeta \leq \Car_{(-\infty, 2+\eta)}.\end{equation}
Let $u,v$ be some unit vectors and $z\in \Omega$. We set
$$\forall Y \in \mathcal{H}_n^{(\beta)}, \  \mathcal{R}_z(Y)= \langle u,R(z) v\rangle \zeta(\lambda_{Y/\sqrt{n}}),$$
where $R(z) = (z-Y/\sqrt{n})^{-1}$. By Weyl's inequality, this defines a $L_{\Omega}$-Lipschitz function with respect to $|| \ ||_{\ell^2}$, where $L_{\Omega}$ is a constant depending on the set $\Omega$. As the entries of $X$ satisfies a Poincaré inequality, $X$ has concentration $\mathcal{C}_1$. We deduce from Lemma \ref{conclip} that for $n$ large enough,
$$ \PP( |\mathcal{R}_z(X) - \EE\mathcal{R}_z(X) |>t+\delta_n ) \leq 4\exp\big( -c L_{\Omega}^{-1} \sqrt{n} t\big),$$
where $\delta_n =O( n^{-1/2})$. Note that $\mathcal{R}_z$ defines a $1/(\eta -\delta)^2$-Lipschitz function in $z\in \Omega$. As $\Omega$ is  relatively compact, we deduce by an $\eps$-net argument that for any $t>0$,
$$\sup_{||u||=||v||=1} \PP( \sup_{z \in \Omega}|\mathcal{R}_z(X) - \EE\mathcal{R}_z(X) |>t ) \underset{n\to+\infty}{\longrightarrow } 0.$$
In the following lemma, we show an isotropic-like property.

\begin{Lem}
Let $\delta >0$ and $\Omega$ a subset compactly included in $\{ z \in \CC : \Re z>2+\delta\}$. Let $X$ be a Wigner matrix satisfying the assumptions of Proposition \ref{compexpvp}. For any $m\in \NN$,
$$\sup_{z \in \Omega}\sup_{u,v \in \mathcal{V}_m} \big| \langle u, \EE[ \zeta(\lambda_{X/\sqrt{n}})R(z)]v \rangle - \langle u,v \rangle g_{\mu_{sc}}(z) \big| \underset{n\to+\infty}{\longrightarrow} 0,$$
where $\mathcal{V}_m$ denotes the set of unit $m$-sparse vectors, meaning with at most $m$ non-zero entries, $\zeta$ is as in \eqref{defzeta}, and $R(z) = (z-X/\sqrt{n})^{-1}$.
\end{Lem}

\begin{proof}
By polarization, it is sufficient to prove this lemma where the supremum ranges over vectors $v=u$. Moreover, by symmetry, it is enough to show this statement for $\Omega \cap \CC^+$.  Because  
$\mathcal{R}_z$, as a function of $z$, is a Lipschitz function on $\Omega$, we only need to show for any $\eps>0$,
$$\sup_{z \in \Omega_{\eps}}\sup_{u\in \mathcal{V}_m} \big| \langle u, \EE[ \zeta(\lambda_{X/\sqrt{n}})R(z)]u \rangle - \langle u,u \rangle g_{\mu_{sc}}(z) \big| \underset{n\to+\infty}{\longrightarrow} 0,$$
with $\Omega_{\eps} = \{ z \in \Omega : \Im z \geq \eps \}$.
 Let $u \in \mathcal{V}_m$. For any $z\in \CC^+$, we have on one hand,
$$\big|\langle u, \EE[ \zeta(\lambda_{X/\sqrt{n}})R(z)]u \rangle - 
\langle u, \EE R(z)u \rangle\big| \leq \frac{1}{\Im z}\PP( \lambda_{X/\sqrt{n}} > 2+\delta).$$
On the other hand, expanding the scalar product,
$$\big|\langle u, \EE Ru \rangle - g_{\mu_{sc}}(z)\big|\leq m \max_{1\leq i,j \leq n}| \EE R_{ij}-\delta_{i,j}g_{\mu_{sc}}(z)|.$$
As $\lambda_{X/\sqrt{n}}$ converges to $2$ in probability, we are reduced to prove for any $\veps>0$,
$$\sup_{z \in \Omega_{\veps}} \max_{1\leq i \leq n} | \EE R_{i,j} - \delta_{i,j} g_{\mu_{sc}}(z) | \underset{n\to +\infty}{\longrightarrow} 0.$$
Even though this is a classical estimate of random matrix theory, for sake of completeness we give here a proof.  
We start with the case of the off-diagonal entries. We set $H = X/\sqrt{n}$ and we write $R$ as a short-hand for $R(z)$. Let $i\neq j$. We have the following resolvent identity (see \cite[Lemma 3.5]{BGKn}), 
$$R_{i,j} = R_{i,i} \sum_{k}^{(i)} H_{i,k} R^{(i)}_{k,j},$$
where $R^{(i)}$ is the resolvent of the matrix $H$ where we removed the $i^{\text{th}}$-row and $i^{\text{th}}$-column, and $\sum^{(i)}$ means that the summation is over $\{ 1,...,n\}\setminus\{i\}$. By Cauchy-Schwarz inequality we have
$$ | \EE R_{i,j}| \leq \frac{1}{\Im z} \Big( \EE \big| \sum_{k}^{(i)} H_{k,i} R^{(i)}_{k,j}\big|^2\Big)^{1/2}.$$
But, as $R^{(i)}$ is independent of $(H_{k,i})_k$ and $(H_{k,l})_{k\leq l}$ are centered and independent,
$$\EE\big|\sum_{k}^{(i)} H_{k,i} R^{(i)}_{k,j}\big|^2 = \frac{1}{n}\EE \sum_{k}^{(i)}  |R^{(i)}_{k,j}|^2.$$
Recall Ward's identity (see \cite[(3.6)]{BGKn}),
$$ \sum_{k}^{(i)}  |R^{(i)}_{k,j}|^2 =- \frac{1}{\Im z} \Im R^{(i)}_{i,i}.$$
Thus, $$|\EE R_{i,j}|\leq \frac{1}{\sqrt{n} (\Im z)^2}.$$
To deal with the diagonal entries, we start from the Schur complement formula (see \cite[Lemma 2.4.6]{Guionnet}),
\begin{equation}\label{Schr}R_{i,i}^{-1} = z -H_{i,i}+ \langle H^{(i)}, R^{(i)} H^{(i)}\rangle,\end{equation}
where $H^{(i)}$ denotes the $i^{\text{th}}$-column of $H$ where the entry $H_{i,i}$ is removed. Let $\mathcal{F}^{(i)}$ be the $\sigma$-algebra generated by the variables $H_{k,l}$ for $k,l \neq i$. We find,
$$\EE \Big(\big|\langle H^{(i)}, R^{(i)} H^{(i)}\rangle - \frac{1}{n} \tr R^{(i)}\big|^2 | \mathcal{F}^{(i)}\Big) = \frac{|\gamma|^2}{n^2} \sum_{k\neq l}^{(i)} R^{(i)}_{k,l} \overline{R_{l,k}^{(i)}} +  \frac{1}{n^2}\sum_{k\neq  l}^{(i)}| R^{(i)}_{k,l}|^2 + \frac{\gamma'}{n^2} \sum_k^{(i)}|R_{k,k}^{(i)}|^2,$$
where $ \gamma = \EE(X_{1,2}^2)$ and $\gamma' = \EE |X_{1,2}|^4-1$. Introducing the missing diagonal terms, using Ward's identity again and the fact that $|\gamma|\leq 1$, we find,
$$\EE \Big(\big|\langle H^{(i)}, R^{(i)} H^{(i)}\rangle - \frac{1}{n} \tr R^{(i)}\big|^2 | \mathcal{F}^{(i)}\Big) \leq \frac{1}{n^2} |\tr R^{(i)} \overline{R^{(i)}}| + \frac{1}{n^2 \Im z } | \Im R^{(i)}_{i,i}| + \frac{c}{n\Im z},$$
where $c$ is some positive constant depending on $\EE |X_{1,2}|^4$.
This yields,
$$\max_{1\leq i\leq n}\EE\big|\langle H^{(i)}, R^{(i)} H^{(i)}\rangle - \frac{1}{n} \tr R^{(i)}\big|^2  = O( n^{-1/2}(\Im z)^{-2}).$$ 
 From Wigner's theorem, we know that $n^{-1}\tr R^{(i)}$ converges to $g_{\mu_{sc}}$ in probability for any $\Im z>0$. Note that $R^{(i)}$ are identically distributed for $i=1,...,n$.
We deduce from \eqref{Schr} and the fact that $g_{\mu_{sc}}(z)^{-1} = z- g_{\mu_{sc}}(z)$ (see \cite[Example 5.3.2.6]{Guionnet}),
$$\max_{1\leq i\leq n}\EE\big|R_{i,i}^{-1} - g_{\mu_{sc}}(z)^{-1}\big|^2\underset{n\to +\infty}{\longrightarrow} 0,$$
which yields, 
$$\max_{1\leq i\leq n}\EE\big|R_{i,i} - g_{\mu_{sc}}(z)\big|^2\underset{n\to +\infty}{\longrightarrow} 0,$$
for any $z\in \CC^+$. As the functions $R_{i,i}$ and $g_{\mu_{sc}}$ are $\eps^{-2}$-Lipschitz on $\{ z \in \CC^+ : \Im z > \eps\}$, we can extend by an $\eps$-net argument, this convergence uniformly on any bounded subset  of $\{z \in \CC^+ : \Im z >\eps\}$, for any $\eps>0$.

\end{proof}
We come back now to the proof of Lemma \ref{compeq}. The above lemma yields that for any $t>0$,
$$\sup_{u\in \mathcal{V }_m} \PP(\{ \sup_{z \in \Omega}|\langle u, Ru \rangle  - g_{\mu_{sc}}(z) |>t\}\cap W_{\delta} ) \underset{n\to+\infty}{\longrightarrow } 0.$$
Note that $m$-sparse matrices have $m$-sparse eigenvectors. Using the fact that the spectral radius of matrices in $\mathcal{G}$ is bounded and a union bound, we deduce that for any $s>0$,
$$ \sup_{A\in \mathcal{G}}\PP( \{ \sup_{z \in \Omega} || U^*RU D - g_{\mu_{sc}}(z)D||_{\ell^{\infty}}>s\}\cap W_{\delta}) \underset{n\to +\infty}{\longrightarrow } 0,$$
where $|| Y||_{\ell^{\infty}} = \sup_{i,j}|Y_{i,j}|$, for any matrix $Y$.
As the matrices $(I_m - g_{\mu_{sc}}(z)D)$, $z \in \Omega$, $A\in \mathcal{G}$ form a pre-compact subset of $\mathcal{H}_m^{(\beta)}$, the continuity of the determinant on $\mathcal{H}_m^{(\beta)}$, allows us to conclude the proof of Lemma \ref{compeq}.

\end{proof}

\subsection{A chaining argument}
We will  now give a proof of Propositions \ref{chaining0} and \ref{chaining1}. 
As it will rely on a chaining argument, we will need the following  lemma.

\begin{Lem}\label{coveringnb}Let $m\in \NN$ and let $B_{\ell^{p}}$ denote the $\ell^p$-ball of $\CC^m$ for any $p>0$. Fix some $0<p<q< \infty$. 
 We denote by $N(B_{\ell^{p}}, \eps B_{\ell^q})$, the covering number of $B_{\ell^{p}}$ by  $\eps B_{\ell^q}$, that is, the minimal number of translates of $\eps B_{\ell^q}$ needed to cover $B_{\ell^{p}}$. There is a constant $c>0$ depending on $p,q$, such that for  $c (\frac{\log m}{m})^{\frac{1}{p}-\frac{1}{q}}\leq \eps \leq c^{-1} $,
$$\log N(B_{\ell^p}, \eps B_{\ell^q}) \leq c  \eps^{\frac{1}{q}-\frac{1}{p}}\log m.$$

\end{Lem}

\begin{proof}
This estimate is a consequence of the upper bound on entropy numbers of embeddings of $\ell_p^m$  in $\ell_q^m$ given in \cite[Proposition 3.2.2]{ETentropy}. Let $0<p<q< \infty$. Denote by $\ell_p^m$ the space $\RR^m$ equipped with the (quasi)-norm $||\ ||_{\ell_p}$.  We define, for $k\in \NN$,
$$e_k(\ell_p^m \to  \ell_q^m)  =\inf \{ \eps>0 : B_{\ell_p} \text{ can be covered by } 2^{k-1} \text{ balls } \eps B_{\ell_q}\}.$$ 
From \cite[Proposition 3.2.2]{ETentropy}, we know that there is a constant $c>0$ such that for $\log_2(2m) \leq k\leq 2m$,
$$ e_k(\ell_p^m \to  \ell_q^m) \leq c \Big( k^{-1}\log_2\big(1+\frac{2m}{k}\big)\Big)^{\frac{1}{p}-\frac{1}{q}}.$$
Thus, if we set $k= \lambda \log_2(2m)$, for some $\lambda\geq 1$ such that $k\leq 2m$, we deduce the following rough bound,
$$ e_k(\ell_p^m \to  \ell_q^m) \leq c' \lambda^{\frac{1}{q}-\frac{1}{p}},$$
for some constant $c'>0$. Let now $\eps>0$ and set $\lambda$ such that $\eps = c' \lambda^{\frac{1}{q} -\frac{1}{p}}$. The above inequality tells us that if $1 \leq \lambda\leq 2m/\log_2(2m)$, then there are $(2m)^{\lambda}$ balls $\eps B_{\ell^q}$ covering $B_{\ell^p}$, that is,
$$N(B_{\ell^p},\eps B_{\ell^q}) \leq (2m)^{\lambda},$$
which yields the claim.
\end{proof}

We are now ready to give a proof of Proposition \ref{chaining0} and \ref{chaining1}.

\begin{proof}[Proof of Proposition \ref{chaining0}]
Let $H\in \mathcal{H}_n^{(\beta)}$. As $X$ satisfies $\mathcal{C}_{\alpha}$ for some constant $\kappa>0$, we see that $X+\sqrt{n}H$ also satisfies $\mathcal{C}_{\alpha}$ with the same constant $\kappa$. We know from Propositions \ref{concspintro} and \ref{citeBC}, that for any $t>0$,
$$ \PP\left( d\big( \mu_{X/\sqrt{n} + H}, \mu_{sc}\boxplus \mu_H\big) > t +\eps_n\right) \leq \frac{32}{t^2}\exp\big ( -c_{\alpha}k_{\alpha}(t) \big),$$
with $k_{\alpha}$ defined in Proposition \ref{concspintro} and $\eps_n =O\big( n^{-1/2} ( \log n)^{(1/\alpha-1)_+}\big)$, uniformly in $H \in \mathcal{H}_n^{(\beta)}$.
Note that the map
$$ S : H \in \mathcal{H}_n^{(\beta)} \mapsto d\big( \mu_{X/\sqrt{n} +H}, \mu_{sc}\boxplus \mu_{H} \big),$$
is $n^{-1/2}$-Lipschitz with respect to $|| \ ||_{\ell^2}$ by Lemma \ref{spvar1intro}. 
We deduce using an $\eps$-net argument that for $n$ large enough,
\begin{equation} \label{chaining}  \PP\big( \sup_{H \in r n^{1/\alpha}B_{\ell^{\alpha}}} S(H) > 2t\big)\leq \frac{32}{t^2} N( r n^{1/\alpha} B_{\ell^{\alpha}}, tn^{ 1/2}  B_{\ell^2})e^{-\frac{c_{\alpha}}{\kappa^{\alpha}}(t-\eps_n)_+^{\alpha}n^{1+\alpha/2} },\end{equation}
where $N( r n^{1/\alpha} B_{\ell^{\alpha}}, tn^{ 1/2}  B_{\ell^2})$ denotes the covering number of $r n^{1/\alpha} B_{\ell^{\alpha}}$ by $tn^{ 1/2}  B_{\ell^2}$.
But, the homogeneity of the norm gives,
$$N( r n^{1/\alpha} B_{\ell^{\alpha}}, tn^{ 1/2}  B_{\ell^2})  = N(  B_{\ell^{\alpha}}, t'n^{ \frac{1}{2} - \frac{1}{\alpha} }  B_{\ell^2}),$$
with $t' = t/r$.
We get from Lemma \ref{coveringnb} applied with $m=n^2$, 
$$ \log N(  B_{\ell^{\alpha}}, t'n^{ \frac{1}{2} - \frac{1}{\alpha} }  B_{\ell^2}) = O(n\log n),$$
This shows that the covering number is negligible with respect to the speed of the deviations, which concludes the chaining argument.
\end{proof}
We finally give a proof of  Proposition \ref{chaining1}.

\begin{proof}[Proof of Proposition \ref{chaining1}]
Let $r>0$. Similarly as in the proof of Proposition \ref{chaining0}, we deduce from Propositions \ref{concvpintro} and \ref{compexpvp}, that for any $A\in \mathcal{H}_n^{(\beta)}$ and $t>0$, 
$$ \PP\left( \big| \lambda_{X/\sqrt{n} + A}- \rho(\lambda_A)\big| > t +\delta_n\right) \leq 8\exp\big(-c_{\alpha} h_{\alpha}(t)  \big),$$
 where $h_{\alpha}$ is defined in Proposition \ref{concvpintro},  $\delta_n = O(n^{-1/2} (\log n)^{(1/\alpha-1)_+})$ uniformly in $A\in r B_{\ell^2}$, and $\rho$ is as in \eqref{deff}. 

Note that the map $x\mapsto \rho(x)$ is $1$-Lipschitz. From Weyl's inequality \cite[Theorem III.2.1]{Bhatia}, we deduce that
$$ A \mapsto |\lambda_{X/\sqrt{n} + A}-\rho(\lambda_A)|,$$
is $2$-Lipschitz with respect to the Hilbert-Schmidt norm on $\mathcal{H}_n^{(\beta)}$. Using an $\eps$-net argument as in the proof of Proposition \ref{chaining0}, it is sufficient to prove that for any fixed $t>0$, the covering number $N(  B_{\ell^{\alpha}} , tB_{\ell^2})$ is negligible at the exponential scale $n^{\alpha/2}$, that is
 $$ \log N(  B_{\ell^{\alpha}} , tB_{\ell^2}) =o(n^{\alpha/2}).$$
But from Lemma \ref{coveringnb}, we know that,
 $$ \log N(  B_{\ell^{\alpha}} , tB_{\ell^2})= O(\log n),$$
 which ends the proof of the claim.
\end{proof}

\subsection{Traces of polynomials of deformed Wigner matrices}
We will now prove Proposition \ref{convunifpoly}. Contrary to the spectral measure or the largest eigenvalue, the proof will consist in a simple moment computation.
\begin{proof}[Proof of Proposition \ref{convunifpoly}]
By linearity it is sufficient to show the statement when $P$ is a monomial, which we will assume from now on. We can write $P = X_{i_1}...X_{i_q}$, with $q\leq d$. Define the matrix $Q$ with coefficients in $\CC\langle \textbf{X}\rangle$, by
$$ Q =   
  \left(
     \raisebox{0.5\depth}{%
       \xymatrixcolsep{1ex}%
       \xymatrixrowsep{1ex}%
       \xymatrix{
         0\ar @{.}[ddddrrrr]& X_{i_1} \ar @{.}[dddrrr] &  & &     \\
         & & & & \\
         &&&& \\
         &&&& X_ {i_{q-1}}\\
         X_{i_q}& & & & 0
       }%
     }\right).
  $$
Observe that by cyclicity of the trace, for any $\textbf{Y} \in (\mathcal{H}_n^{(\beta)})^p$, $\tr Q(\textbf{Y})^q = q \tr P(\textbf{Y})$. Therefore,
\begin{equation} \label{cycl} \tr P(\textbf{X}/\sqrt{n} + n^{1/d}\textbf{H}) =\frac{1}{q} \tr \big( Q(\textbf{X}/\sqrt{n}) + n^{1/d}Q(\textbf{H})\big)^d.\end{equation}
Write $Z = Q(\textbf{X}/\sqrt{n})$ and $K = Q(\textbf{H})$. We know from the proof of \cite[Lemma 2.1]{LDPtr} that,
$$\big | \tr \big( Z + n^{1/d}K\big)^q -\tr Z^q - n^{\frac{q}{d}} \tr K^q\big| \leq 2^q \max_{1\leq k\leq q-1}n^{\frac{q-k}{d}} ( \tr | Z|^{q+1})^{\frac{k}{q+1}} (\tr |K|^2)^{\frac{q-k}{2}}.$$
Let us define $q$-Schatten (quasi-)norm on $(\mathcal{H}_n^{(\beta)})^p$, for any $q>0$ by,
\begin{equation} \label{defschattenp} \forall \textbf{H} \in (\mathcal{H}_n^{(\beta)})^p, \ || \textbf{H} ||_{q} = \Big( \sum_{i=1}^p \tr|H_i|^q \Big)^{1/q}.\end{equation}
Note that for any $\textbf{Y}\in (\mathcal{H}_n^{(\beta)})^p$,
$$|Q(\textbf{Y})|=
  \left(
     \raisebox{0.5\depth}{%
       \xymatrixcolsep{1ex}%
       \xymatrixrowsep{1ex}%
       \xymatrix{
         |Y_{i_1}| \ar @{.}[dddrrr]& 0\ar @{.}[ddrr]  \ar @{.}[rr]  &   & 0 \ar @{.}[dd]   \\
         0 \ar @{.}[dd] \ar @{.}[rrdd]  &  & & \\
         &&& 0\\
        0 \ar @{.}[rr]  & & 0 & |Y_{i_q}|
       }%
     }\right).
   $$
Thus, for any $m\in \NN$,
$$\tr |Q(\textbf{Y})|^{m} = \sum_{j=1}^q \tr|Y_{i_j}|^{m}  \leq \sum_{i=1}^p \tr|Y_{i}|^{m}  = ||\textbf{Y}||_{m}^{m}.$$
As $\textbf{H}\in r B_{\ell^2}$, $\tr|K|^2 \leq r^2$. Without loss of generality we can assume $r\geq1$. Thus,
\begin{equation} \label{linetr}\big | \tr \big( Z + n^{1/d}K\big)^q -\tr (Z^q) - n^{\frac{q}{d}}\tr K^q\big| \leq r^q 2^q \max_{1\leq k\leq q-1}n^{\frac{q-k}{d}} || \textbf{X}/\sqrt{n}||_{q+1}^{k}.\end{equation}
But we know from Wigner's theorem (see \cite[Lemma 2.1.6]{Guionnet}), that there is a constant $c\geq 1$, such that 
$$ \EE || \textbf{X}/\sqrt{n} ||_{q+1}^{q+1} \leq c n.$$
Besides,
$$\EE \max_{1\leq k \leq q-1} n^{-\frac{k}{d}} || \textbf{X}/\sqrt{n}||_{q+1}^{k} \leq \sum_{k=1}^{q-1} n^{-\frac{k}{q}} \EE  || \textbf{X}/\sqrt{n}||_{q+1}^{k}.$$
By Jensen's inequality, we deduce
$$\EE \max_{1\leq k \leq q-1} n^{-\frac{k}{d}} || \textbf{X}/\sqrt{n}||_{q+1}^{k} \leq \sum_{k=1}^{q-1} n^{-\frac{k}{q}}\big( \EE  || \textbf{X}/\sqrt{n}||_{q+1}^{q+1}\big)^{\frac{k}{q+1}}.$$
Therefore,
$$\EE \max_{1\leq k \leq q-1} n^{-\frac{k}{d}} || \textbf{X}/\sqrt{n}||_{d+1}^{k} \leq qc n^{-( \frac{1}{q}-\frac{1}{q+1})}.$$
We deduce from \eqref{cycl} and \eqref{linetr} that 
$$ \big| \tau_n[P(\textbf{X}/\sqrt{n}+n^{1/d}\textbf{H})] - \EE\tau_n[P(\textbf{X}/\sqrt{n})]-n^{\frac{q}{d}-1}\tr[P(\textbf{H})] \big| \underset{n\to +\infty}{\longrightarrow} 0,$$
uniformly in $\textbf{H}\in r B_{\ell^2}$ and where $\tau_n = \frac{1}{n} \tr$. It is now sufficient to prove that $n^{q/d-1}\tr P(\textbf{H})$ converges to $0$ uniformly in $\textbf{H} \in r B_{\ell^{\alpha}}$, as soon as $q<d$.
Assume first $q\geq \alpha$. Using the non-commutative Hölder's inequality (see \cite[Corollary IV.2.6]{Bhatia}), we get
\begin{equation*} \tr [P(\textbf{H})] \leq \prod_{j=1}^q || H_{i_j} ||_q. \end{equation*}
The arithmetic-geometric mean inequality yields,
\begin{equation} \label{ineqSchatten} \tr [P(\textbf{H})] \leq \frac{1}{q} \sum_{j=1}^q  \tr |H_{i_j}|^q.\end{equation}
As $q\geq \alpha$, we deduce
\begin{equation}\label{comptrnorm} \tr [P(\textbf{H})] \leq   || \textbf{H}||_{q}^q  \leq  || \textbf{H}||_{\alpha}^q. \end{equation}
We conclude that when $\alpha\leq q<d$, 
$$\sup_{\textbf{H}\in r B_{\ell^{\alpha}}} n^{\frac{q}{d}-1} \tr [P(\textbf{H})] \underset{n\to+\infty}{\longrightarrow} 0.$$
If $q<\alpha$, then $q=1$ and $\alpha>1$.  By Jensen's inequality, 
$$|\tr H_{i_1}| \leq n^{1-1/\alpha}(\tr |H_{i_1}|^{\alpha})^{1/\alpha}.$$ Thus, as $d>\alpha$,
$$\sup_{\textbf{H}\in r B_{\ell^{\alpha}}} n^{\frac{1}{d}-1} \tr [P(\textbf{H})] \underset{n\to+\infty}{\longrightarrow} 0.$$
Besides, we know by \cite[Theorem 5.4.2]{Guionnet}, that 
$$ \EE \tau_n[P(\textbf{X}/\sqrt{n})] \underset{n\to+\infty}{\longrightarrow} \tau[P(\textbf{s})],$$
where $\textbf{s}$ are a family of $p$ free semi-circular variables defined on a non-commutative probability space $(\mathcal{A},\tau)$. This ends the proof of the proposition.
\end{proof}

\section{Deterministic equivalent for the last-passage time}\label{secdetermequivLPT}
We will prove in this section the analogue of the results for Wigner matrices of the preceding section, for the last-passage time. More precisely, we will provide a deterministic equivalent for the last-passage time when the matrix of weights is deformed by some matrix $nH$, where $||H||_{\ell^{\alpha}}$ is bounded for some $\alpha \in (0,1)$.

Let $\mathcal{A}$ denote the set of finite vectors $(v_1,...,v_m)$, which we will call \textit{admissible}, such that $v_i\in \{0,...,n\}^d$, $v_0=(0,...,0)$, $v_m=(n,...,n)$, and for any $i \in \{0,...,m-1\}$, $v_i< v_{i+1}$, where $<$ denotes the lexicographic order. With this definition we set, for any $H\in \RR^I$, where $I = \{0,...,n\}^d$, 
\begin{equation} \label{defequidetLPP}\mathcal{T}_n(H) = \sup_{V \in \mathcal{A} }\Big\{ \sum_{i=0}^{m} H_{v_i}^+ + \sum_{i=0}^{m-1} g\Big( \frac{v_{i+1}-v_i}{n}\Big)\Big\},\end{equation}
where $V=(v_0,...,v_m)$ for some $m\in \NN$, where $g$ is as in \eqref{defg}, and where we denote here, for better lisibility, $x^+$ the positive part of $x\in \RR$ ($x^+ = x_+$, our former notation).
With this notation, we will prove the following proposition.

\begin{Pro}\label{convunifLPP}Let $\alpha \in (0,1)$. Let $X = (X_{v})_{v \in \ZZ_+^d}$ be a family of i.i.d random variables following the law $\mu_{\alpha}$. For any $r>0$, 
$$\sup_{||H||_{\ell^{\alpha}}\leq r} \Big|\frac{1}{n}T(X+nH)^+ - \mathcal{T}_n(H)\Big| \underset{n\to +\infty}{\longrightarrow} 0,$$
in probability, where $Y^+$ denotes the multi-matrix $(Y_v^+)_v$.
\end{Pro}

We will follow the same arguments as for  the proof of the uniform deterministic equivalent of the empirical spectral measure and the largest eigenvalue of Wigner matrices. We will begin by showing that the deterministic equivalent \eqref{defequidetLPP} we propose, holds uniformly in expectation. This is the object of the following lemma.

\begin{Lem}
Let $\alpha\in (0,1)$. Let $X = (X_{v})_{v \in \ZZ_+^d}$ be a family of i.i.d non-negative random variables with common distribution function satisfying \eqref{condfuncdistr}. For any $r>0$,
$$\sup_{||H||_{\ell^{\alpha}} \leq r }\big|\frac{1}{n}\EE T(X+nH)^+ - \mathcal{T}_n(H)\big|\underset{n\to+\infty}{\longrightarrow} 0,$$
where $\mathcal{T}_n(H)$ is as in \eqref{defequidetLPP}.
\end{Lem}

\begin{proof}
Let $\mathcal{A}_m$ denote the subset of vectors of $\mathcal{A}$ of size less or equal than $m$, and define $\hat{\mathcal{T}}_n^{(m)}$ by,
$$\hat{\mathcal{T}}_n^{(m)}(H) = \sup_{V \in \mathcal{A}_m }\Big\{ \sum_{i=0}^{p} H_{v_i}^+  +\sum_{i=0}^{p-1}\frac{1}{n} \EE T_{v_i,v_{i+1}}(X) \Big\},$$
and $\mathcal{T}_n^{(m)}$,
$$\mathcal{T}_n^{(m)}(H) = \sup_{V \in \mathcal{A}_m }\Big\{ \sum_{i=0}^{p} H_{v_i}^+ +\sum_{i=0}^{p-1} g\Big( \frac{v_{i+1}-v_i}{n}\Big)\Big\},$$
where $V=(v_0,...,v_p)$ for some $p\leq m$, and $g$ is as in \eqref{defg}.
 We begin by proving that there is some constant $C>0$ depending on $\alpha$, such that for any $|| H ||_{\ell^{\alpha}}\leq r$,
\begin{equation} \label{encadrelasttime} -Cr(\log n)^{\frac{1}{\alpha}} n^{\alpha-1}\leq  \frac{1}{n}\EE T(X+nH)^+ - \hat{\mathcal{T}}_n^{(m)}(H) \leq C rm^{1-\frac{1}{\alpha}}.\end{equation}
In the following $C$ will denote a constant which will depend only on $\alpha$ and which will vary along the lines of the proof.
Let $\pi$ be an optimal path for the last-passage time $T(X+nH)^+$, and denote by $v_1,..,v_{m-1}$ be the $m-1$ largest values of $H^+$ on the path $\pi$, sorted in lexicographic order. Add $v_0 = (0,...,0)$ and $v_m=(n,...,n)$, to get $V=\{v_0,...,v_m\} \in \mathcal{A}_m$. We have
$$ \frac{1}{n}T(X+nH)^+ - \sum_{i=0}^{m} H_{v_i}^+  -\sum_{i=0}^{m-1}\frac{1}{n}  T_{v_i,v_{i+1}}(X)\leq \frac{1}{n}\sum_{v\in \pi } (X +nH)^+_{v} - \sum_{i=0}^{m} H_{v_i}^+ -\frac{1}{n} \sum_{v\in \pi} X_{v}.$$
As $(x+y)^+\leq x^+ + y^+$, we deduce
$$ \frac{1}{n}T(X+nH)^+ - \sum_{i=0}^{m} H_{v_i}^+  -\sum_{i=0}^{m-1}\frac{1}{n}  T_{v_i,v_{i+1}}(X)\leq \sum_{v \in \pi\cap V^c }H_{v}^+.$$ 
Now observe that if $M_1\geq ... \geq M_{d(n+1)}$ are the values of $H^+$ (or of $H^-$) along $\pi$ in decreasing order, we have since $\sum_{i} M_i^{\alpha}\leq r^{\alpha}$, for any $k\in \{1,...,d(n+1)\}$,
\begin{equation} \label{boundweight}M_k \leq r k^{-1/\alpha}.\end{equation} 
Therefore,
$$ \sum_{v \in \pi\cap V^c}H_{v}^+\leq r \sum_{k=m-1}^{+\infty} k^{-1/\alpha} \leq C r m^{1-\frac{1}{\alpha}},$$
for some constant $C>0$. This proves the upper bound of \eqref{encadrelasttime}. On the other hand, let $V=\{v_0,....,v_p\}\in \mathcal{A}_m$. Considering the optimal paths from $v_i$ to $v_{i+1}$ in the last-passage time $T_{v_i,v_{i+1}}(X)$, for $i=0,...,p-1$ and their concatenation $\pi$, we get,
\begin{equation} \label{bornesupLPP} \sum_{i=0}^{p} H_{v_i}^+ + \frac{1}{n}\sum_{i=0}^{p-1}T_{v_i,v_{i+1}}(X) - T(X+nH)^+\leq \sum_{X_{v} \geq -n H_{v}} H_{v}^- + \sum_{ X_{v } \leq -n  H_{v}} \frac{X_{v}}{n}.\end{equation}
Indeed, if $v\in \pi $, then
$$H_v^+ + X_v - (X+nH)_v^+ \leq \Car_{\{X_v \geq -nH_v\}}H_v^-+\Car_{\{X_v \leq -nH_v\}} X_v,$$
by considering the cases whether $H_v\geq 0$ or ($H_v\leq 0$ and $X+nH_v\geq 0$) or ($H_v\leq 0 $ and $X+nH_v\leq 0$).
Turning our attention to the first sum in \eqref{bornesupLPP}, we deduce by bounding the first $n^{\alpha}$ largest weights of $H_{v}^-$ by $X_{v}/n$, and using the bound \eqref{boundweight} for the rest of the terms, 
$$ \EE\Big(\sum_{X_{v} \geq -n H_{v}} H_{v}^-\Big) \leq \frac{n^{\alpha}}{n} \EE\sup_v X_v + r \sum_{k>n^{\alpha}} k^{-\frac{1}{\alpha}}.$$
By \eqref{momentalpha}  we have,
$$ \EE \sup_v X_v \leq c(\log n)^{\frac{1}{\alpha}},$$
for some constant $c>0$. We thus proved,
$$ \EE\Big(\sum_{X_{v} \geq -n H_{v}} H_{v}^-\Big) \leq C r (\log n)^{\frac{1}{\alpha}} n^{\alpha-1}.$$
On the other hand, focusing now on the second term of \eqref{bornesupLPP}, 
 $$\EE\Big(\sum_{ X_{v} \leq -n  H_{v}} \frac{X_{v}}{n} \Big) =\frac{1}{n} \EE \big(X_0 |\{ v : X_0\leq -n H_{v}\}| \big).$$
But $||H||_{\ell^{\alpha}} \leq r$, thus
$$ |\{ v : X_0\leq -n H_{v}\}| \Big( \frac{X_0}{n} \Big)^{\alpha}\leq r.$$
Therefore,
$$\EE\Big(\sum_{ X_{v } \leq -n  H_{v}} \frac{X_{v}}{n} \Big) \leq n^{\alpha-1} r \EE X_0^{1-\alpha}.$$
which concludes the proof of the lower bound of \eqref{encadrelasttime}. Comparing $\mathcal{T}^{(m)}_n$ and $\hat{\mathcal{T}}^{(m)}_n$, we get using the translation invariance in law (by vectors of $\ZZ^d_+$) of $(X_v)_{v\in\ZZ^d_+}$,
$$|\mathcal{T}^{(m)}_n(H) - \hat{\mathcal{T}}_n^{(m)}|\leq m \max_{v\in \{0,...,n\}^d} \Big|\frac{1}{n} \EE T_{0,v}(X)- g\Big( \frac{v}{n}\Big)\Big|.$$

As $\EE T_{0, \lfloor n w\rfloor}(X)$ is coordinate-wise non-decreasing as a function of $w\in \RR_+^2$, and converges to $g(w)$ which is continuous by \cite[Theorem 2.3]{Martin}, we deduce that $w \mapsto \EE T_{0, \lfloor n w\rfloor}(X)$ converges uniformly to $g$ on $[0,1]^2$ by Dini's Theorem. Thus, 
\begin{equation}\label{boundLPP1}|\mathcal{T}^{(m)}_n(H) - \hat{\mathcal{T}}^{(m)}_n|\leq m \eps(n),\end{equation}
where $\eps(n) \to +\infty$ when $n\to +\infty$.

Now, using the same argument as for the upper bound of \eqref{encadrelasttime}, we see that 
\begin{equation}\label{boundLPP2}| \mathcal{T}^{(m)}_n(H) - \mathcal{T}_n(H)|\leq C r m^{1-\frac{1}{\alpha}},\end{equation} for any $||H||_{\ell^{\alpha}}\leq r$. Indeed, if $V$ achieves the supremum in $\mathcal{T}_n(H)$, then taking $V'$ the $m$ largest values of $H^+$ on $V$, we get
$$0\leq \mathcal{T}_n(H) - \mathcal{T}_{n}^{(m)}(H) \leq \sum_{v\notin V'} H_v^+.$$ Thus, using  \eqref{boundweight}, we get the claim. To summarize, we got by \eqref{encadrelasttime}, \eqref{boundLPP1}, and \eqref{boundLPP2},
$$\big|\frac{1}{n}\EE T(X+nH)^+ - \mathcal{T}_n(H)\big| \leq C r m^{1-\frac{1}{\alpha}} + m\eps(n) +C r(\log n )^{\frac{1}{\alpha}} n^{\alpha-1},$$
for some constant $C>0$ and for any $||H||_{\ell^{\alpha}}\leq r$, which gives finally the claim by taking the $\limsup$ as $n\to +\infty$, and then as $m\to +\infty$.
\end{proof}

We can now give a proof of Proposition \ref{convunifLPP}.

\begin{proof}[Proof of Proposition \ref{convunifLPP}]Let $H\in \RR^I$.
Note that $X\mapsto T(X+nH)$ is $1$-Lipschitz with respect to $|| \ ||_{\ell^{1}}$ on $\RR^I$. As $|| \ ||_{\ell^1}\leq || \ ||_{\ell^{\alpha}}$ since $\alpha<1$, we deduce that $X\mapsto T(X+nH)$ is also $1$-Lipschitz with respect to $|| \ ||_{\ell^{\alpha}}$. Moreover by Hölder's inequality, $X\mapsto T(X+nH)$ is $\sqrt{n}$-Lipschitz with respect to $|| \ ||_{\ell^2}$. We get by Lemma \ref{conclip0}, for any $t>0$,
$$ \PP( |T(X+nH) - m |> tn)\leq8 \exp(-c p_{\alpha}(t)),$$
where $m$ is the median of $T(X+nH)$, $c$ is some strictly positive constant, and 
$$p_{\alpha}(t)=\min\Big(\frac{t^2n}{(\log n)^{2(\frac{1}{\alpha}-1)}}, \frac{tn}{(\log n)^{\frac{1}{\alpha}-1}}, n^{\alpha}t^{\alpha}\Big).$$
 Integrating this inequality we get,
$$|\EE T(X+nH) - m| = O( (\log n)^{\frac{1}{\alpha}-1} \sqrt{n}),$$
uniformly in $H$. Using the result of Proposition \ref{convunifLPP}, we deduce that for $n$ large enough, 
\begin{equation}\label{concLPP} \PP\big( |T(X+nH) - \mathcal{T}_n(H) |> (t+\delta_n)n\big)\leq8 e^{-cn^{\alpha}t^{\alpha}},\end{equation}
where $\delta_n = O( (\log n)^{\frac{1}{\alpha}-1} n^{-\frac{1}{2}})$.
Let now $r>0$. 
Note that 
$$H \mapsto n^{-1}  |T(X+nH) - \mathcal{T}_n(H) |,$$ is $2$-Lipschitz with respect to $||\ ||_{\ell^1}$ on $\RR^I$. Besides, by Lemma \ref{coveringnb} for any $\eps>0$, the covering number of $r B_{\ell^{\alpha}}$ by $\ell^2$-balls of radii $\eps$ satisfies,
$$\log N(r B_{\ell^{\alpha}}, \eps B_{\ell^2}) = O(\log n).$$ Since this estimate is negligible with respect to the concentration bound \eqref{concLPP}, we deduce using an $\eps$-net arguments as in the proofs of Propositions \ref{chaining0} and \ref{chaining1}, that 
$$ \PP\Big( \sup_{H\in rB_{\ell^{\alpha}}} \big| \frac{1}{n} T(X+nH) - \mathcal{T}_n(H)\big| >t\Big) \underset{n\to +\infty}{\longrightarrow} 0,$$
which ends the proof of the claim.
\end{proof}

\section{Applications to Wigner matrices}\label{Wigner}
We apply in this section Theorem \ref{theoremgene} in the setting of Wigner matrices, and we derive the LDP of Theorems \ref{LDPmsp}, \ref{LDPvp} and \ref{LDPpoly}. In all this section, $X$ will designate a Wigner matrix with the class $\mathcal{S}_{\alpha}$ for some $\alpha \in (0,2]$. It is clear that Theorem \ref{theoremgene} remains valid in the context of Wigner matrices in the class $\mathcal{S}_{\alpha}$, making the according change in the rate function $I_{\alpha}$, by replacing the weight function $|| \ ||_{\ell^{\alpha}}^{\alpha}$ by $W_{\alpha}$, which defines the law of a Wigner matrix in $\mathcal{S}_{\alpha}$ (see \eqref{defW}).

\subsection{Large deviations of the empirical spectral measure}

\begin{proof}[Proof of Theorem \ref{LDPmsp}]
From Proposition \ref{chaining0}, we know that assumption $(i)$ of Theorem \ref{theoremgene} is satisfied with 
$$\forall H \in \mathcal{H}_n^{(\beta)}, \ F_m(H) = \mu_{sc}\boxplus \mu_{n^{1/\alpha} H},$$
and
$$\forall H \in \mathcal{H}_n^{(\beta)}, \ f_m(X) = \mu_{X/\sqrt{n}+n^{1/\alpha} H}.$$
where $m$ is the (real) dimension of $\mathcal{H}_n^{(\beta)}$, with the metric $d$ on $\mathcal{P}(\RR)$ defined in \eqref{defdStiel},  and $v(m) = n^{1 + \frac{\alpha}{2}}$. 

By Lemma \ref{spvar1intro}, we see that $f_m$ is $n^{-1}$ -Lipschitz with respect to $|| \ ||_{\ell^2}$ on $\mathcal{H}_n^{(\beta)}$ and $d$ on $\mathcal{P}(\RR)$. By the remark \ref{remtheogene}
\hyperref[remlipconst]{(c)}, and from the fact that $ \alpha <2$, we deduce that  the assumption $(ii)$ of Theorem \ref{theoremgene} holds. 
Besides, as $\alpha \leq 2$, we have by \cite[Theorem 3.32]{Zhan}
$$ \forall H \in \mathcal{H}_n^{(\beta)}, \  (\tr |H|^{\alpha} )^{1/\alpha}  \leq   || H ||_{\ell^{\alpha}} .$$
Thus for any $r>0$,
$$ F_m(r B_ {\ell^{\alpha}} ) \subset \{ \mu \in \mathcal{P}(\RR) : \mu |x|^{\alpha} \leq r^{\alpha} \},$$
which shows that $\cup_{m} F_m(rB_{\ell^{\alpha}})$ is relatively compact by Prokhorov's theorem, and that $(iii)$ is verified.

To prove $(iv)$ it is sufficient to show that for a fixed $H\in \mathcal{H}_p^{(\beta)}$, there is a sequence $H_n \in \mathcal{H}_n^{(\beta)}$, $n\geq p$, such that
\begin{equation} \label{proptaux} \lim_{n\to +\infty}\mu_{n^{1/\alpha}H_n}=  \mu_{p^{1/\alpha}H}, \quad \text{ and }\quad \lim_{n\to +\infty} W_{\alpha}(H_n) = W_{\alpha}(H).\end{equation}
Let for any $k\in \NN$, $H_{kp} = \oplus_{i=1}^k k^{-1/\alpha}H\in \mathcal{H}_{kp}^{(\beta)}$. We have $W_{\alpha}(H_{kp}) = W_{\alpha}(H)$, as $W_{\alpha}(\lambda Y) =\lambda^{\alpha} W_{\alpha}(Y)$ for any $\lambda>0$, and
$$ \mu_{(kp)^{1/\alpha} H_{kp}} = \mu_{p^{1/\alpha} H}.$$
Now, if $n = kp + l$, with $k\in \NN$ and $1\leq l \leq p$, we define
$$ H_{n} = \Big(\frac{kp}{kp+l}\Big)^{1/\alpha} \left(\begin{array}{cc}
H_{kp} & 0\\
0 & 0
\end{array}\right) \in \mathcal{H}_n^{(\beta)}.$$
We have,
$$ \mu_{n^{1/\alpha} H_n} = \frac{kp}{kp+l} \mu_{(kp)^{1/\alpha} H_{ kp}} + \frac{l}{kp+l} \delta_0.$$
Thus,
$$d( \mu_{n^{1/\alpha} H_n}, \mu_{(kp)^{1/\alpha} H_{kp}}) \leq \frac{2 l }{kp+l} \leq \frac{2p}{n}.$$
 Besides,
$$ W_{\alpha}(H_{kp}) \geq W_{\alpha}(H_n)\geq \big(1+\frac{1}{k}\big)^{-\frac{1}{\alpha}} W_{\alpha}(H_{kp}).$$
As $W_{\alpha}(H_{kp}) = W_{\alpha}(H)$, and $\mu_{(kp)^{1/\alpha} H_{kp}} = \mu_{p^{1/\alpha}H}$, we get the claim \eqref{proptaux}.

\end{proof}

\subsection{Large deviations of the largest eigenvalue}

\begin{proof}[Proof of Theorem \ref{LDPvp}]

We begin by giving back to $J_{\alpha}$ its variational form. We claim that for any $x\in \RR$,
\begin{equation} \label{Jvar1} J_{\alpha}(x) =  \sup_{\delta>0}\inf\big \{ W_{\alpha}(A)  : A \in \cup_{n \in \NN} \mathcal{H}_n^{(\beta)}, \ | x - \rho(\lambda_A)|<\delta \big\},\end{equation}
where $\rho$ is the function
\begin{equation*} \forall x \in \RR, \  \rho(x) = \begin{cases} x +\frac{1}{x} & \text{ if } x \geq 1,\\
2 & \text{ otherwise.}
\end{cases}
\end{equation*}
Let us prove first that
\begin{equation} \label{Jvar}\forall x \in \RR, \ J_{\alpha}(x) = \inf\big \{ W_{\alpha}(A)  : A \in \cup_{n \in \NN} \mathcal{H}_n^{(\beta)}, \ x = \rho(\lambda_A) \big\}.\end{equation}
When $x<2$, both sides of \eqref{Jvar} are infinite. If $x\geq2$, we denote by $\mathcal{J}_{\alpha}$ the right-hand side of \eqref{Jvar}. The function $x \in (0,1] \mapsto \rho(1/x)$ is the inverse of the Stieltjes transform of $\mu_{sc}$ on $[2,+\infty)$ (see  \cite[Example 5.3.2.6]{Guionnet}). Thus, we can write
 $$\mathcal{J}_{\alpha}(x) = \inf\big \{ W_{\alpha}(A)  : A \in \cup_{n \in \NN} \mathcal{H}_n^{(\beta)}, \  1/\lambda_A = g_{\mu_{sc}}(x) \big\}.$$
 As $W_{\alpha}$ is $\alpha$-homogeneous, and $\lambda_{tA} = t \lambda_A$, for any $t\geq 0$, we get
$$ \mathcal{J}_{\alpha}(x) = \mathcal{J}_{\alpha}(1) g_{\mu_{sc}}(x)^{-\alpha}.$$
Thus, $J_{\alpha} = \mathcal{J}_{\alpha}$. As $J_{\alpha}$ is clearly lower semi-continuous, the equality \eqref{Jvar1} holds by the remark  \ref{remtheogene} \hyperref[remlsi]{(e)}.

%
%

We check now the assumptions of Theorem \ref{theoremgene}. Assumption $(i)$ of Theorem \ref{theoremgene} is met by the result of Proposition \ref{chaining1}, with
$$\forall H \in \mathcal{H}_n, \ f_m(H) = \lambda_{X/\sqrt{n}}, \ F_m(H) = \rho(\lambda_H),$$
where as before $m$ is the dimension of $\mathcal{H}_n^{(\beta)}$, and $v(m) = n^{\alpha/2}$.
 Weyl's inequality \cite[Theorem III.2.1]{Bhatia} shows that $f_m$ is $n^{-1/2}$-Lipschitz with respect to $|| \ ||_{\ell^2}$, and thus assumption $(ii)$ is satisfied as $\alpha <2$ by the remark \ref{remtheogene} 
\hyperref[remlipconst]{(c)}.
Besides, note that for any $H\in \mathcal{H}_n^{(\beta)}$,
$$|\lambda_H| \leq (\tr |H|^{\alpha})^{1/\alpha} \leq || H ||_{\ell^{\alpha}},$$
where we used in the second inequality the fact that $\alpha \leq 2$ and \cite[Theorem 3.32]{Zhan}. As $\rho$ is non-decreasing, we deduce for any $r>0$ that, 
$$ \{ F_m(H) : H\in r B_{\ell^{\alpha}} \} \subset [ 2, \rho(r) ],$$
which proves that $(iii)$ is satisfied. To show that $(iv)$ holds, it suffices to observe that if $H \in \mathcal{H}_n^{(\beta)}$, and if we set for any $m\geq n$, \begin{equation} \label{defH} H_m = \left(\begin{array}{cc}
H_{n} & 0\\
0 & 0
\end{array}\right) \in \mathcal{H}_m^{(\beta)},\end{equation}
then $W_{\alpha}(H_m) = W_{\alpha}(H)$, and  provided $\lambda_H\geq 0$, we have $\lambda_H = \lambda_{H_m}$, so that in particular $\rho(\lambda_H) = \rho(\lambda_{H_m})$.
\end{proof}

\subsection{Large deviations of non-commutative polynomials}

Finally, we give a proof of Theorem \ref{LDPpoly}.

\begin{proof}[Proof of Theorem \ref{LDPpoly}]

By a homogeneity argument similar as for the proof of Theorem \ref{LDPvp}, we get for any $x \in \RR$,
$$K_{\alpha}(x) = \inf\big\{ W_{\alpha}( \textbf{H}):  \textbf{H} \in \cup_{n\in \NN} (\mathcal{H}_n^{(\beta)})^p, x = \tr P_d(\textbf{H}) + \tau(P(\textbf{s}))\big\},$$
where $P_d$ denotes the homogeneous part of degree $d$ of $P$. From the remark  \ref{remtheogene} \hyperref[remlsi]{(e)}, we get as $K_{\alpha}$ is lower semi-continuous, that
$$K_{\alpha}(x) =\sup_{\delta>0} \inf\big\{ W_{\alpha}( \textbf{H}):  \textbf{H} \in \cup_{n\in \NN} (\mathcal{H}_n^{(\beta)})^p, |x - \tr P_d(\textbf{H}) - \tau(P(\textbf{s}))|< \delta \big\}.$$ 
Assumption $(i)$ of Theorem \ref{theoremgene} is a consequence of Lemma \ref{convunifpoly} with the speed $v(m) = n^{\alpha(\frac{1}{2}+\frac{1}{d})}$ and 
$$F_m(\textbf{H}) = \tr P_d(\textbf{H}) + \tau(P(\textbf{s})), \ f_m(\textbf{H}) = \tau_n(P(\textbf{X}/\sqrt{n})),$$
where $m$ is the real dimension of $(\mathcal{H}_n^{(\beta)})^p$.

Let us now prove assumption $(ii)$.  Note that by linearity, it suffices to prove assumption $(ii)$ when $P$ is a monomial of  total degree $k\geq 1$ less or equal than $d$, which we will assume from now on. If $k=1$, then there are two cases to consider. First we see by Hölder's inequality that $f_m$ is $n^{-1}$-Lipschitz with respect to $|| \ ||_{\ell^2}$.  If $d=1$ then $\alpha \in (0,1)$, so that as $v(n) = n^{3\alpha/2}$ in this case. We conclude by remark \ref{remtheogene} \hyperref[remlipconst]{(c)} that assumption $(ii)$ holds. If $d\geq 2$ and $k=1$, then we deduce again by remark \ref{remtheogene} \hyperref[remlipconst]{(c)} that assumption $(ii)$ is fulfilled as $v(n) = n^{\alpha(\frac{1}{2}+\frac{1}{d})}$.

In the case $k\geq 2$, we will need to understand the stability of the function $f_m$ with respect to the Euclidean norm. This is the object of the following lemma.

\begin{Lem}\label{meanvalue}  There is a constant $C_{d,p}>0$ depending on $d$ and $p$, such that for any monomial $q \in \CC\langle \textbf{X} \rangle$ of total degree $d\geq 2$, and $\textbf{Y},\textbf{H} \in (\mathcal{H}_n^{(\beta)})^p$, 
$$ | \tr q(\textbf{Y} +\textbf{H}) - \tr q( \textbf{Y})| \leq C_{d,p} \big( || \textbf{Y} ||_{2(d-1)}^{d-1} + || \textbf{H}||_{2}^{d-1}\big)||\textbf{H}||_{2},$$
 where for any $q>0$, $|| \ ||_q$ denotes the $q$-Schatten norm on $(\mathcal{H}_n^{(\beta)})^p$, defined in \eqref{defschattenp}.
\end{Lem}
\begin{proof}Let 
$$ \forall \textbf{H} \in (\mathcal{H}_n^{(\beta)})^p, \ f(\textbf{H}) = \tr q(\textbf{H}).$$
By the mean value theorem, we have
\begin{equation} \label{meanv} | f(\textbf{Y} +\textbf{H}) - f( \textbf{Y})| \leq \max_{0\leq t\leq 1}|| \nabla f(\textbf{Y} + t\textbf{H})||_{2} ||\textbf{H}||_{2}.\end{equation}
Note that if $R\in \CC\langle \textbf{X}\rangle$ is a monomial of degree $d-1$ in $\textbf{X}$, then by \eqref{comptrnorm}, we have
$$ \tr |R(\textbf{Z})|^{2} \leq || \textbf{Z} || _{2(d-1)}^{2(d-1)}.$$
As $ \nabla_{X_i} f$ is the sum of at most $d$ monomials of degree $d-1$ in $\textbf{X}$, we get by triangular inequality and the above observation, 
$$|| \nabla_{X_i} f(\textbf{Z})||_{2} \leq d  || \textbf{Z}||_{2(d-1)}^{d-1}.$$  
Thus, 
$$ || \nabla f(\textbf{Y} + t\textbf{H})||_{2} \leq p d|| \textbf{Y}+t\textbf{H}||_{2(d-1)}^{d-1}.$$  
As $\textbf{Z} \mapsto || \textbf{Z}||_{2(d-1)}^{d-1}$ is convex, we get 
$$ || \nabla f(\textbf{Y} + t\textbf{H})||_{2} \leq dp (1+t)^{d-2}\big( ||\textbf{Y}||_{2(d-1)}^{d-1} + t ||\textbf{H}||_{2(d-1)}^{d-1}\big).$$ 
As $2(d-1) \geq 2$, we have
$$ || \nabla f(\textbf{Y} + t\textbf{H})||_{2} \leq dp 2^{d-2}\big( ||\textbf{Y}||_{2(d-1)}^{d-1} + ||\textbf{H}||_{2}^{d-1}\big).$$
This inequality together with \eqref{meanv} which yields the claim \eqref{meanvalue}.
\end{proof}

We come back now at the proof of assumption $(ii)$ of Theorem \ref{theoremgene}. Let $r\geq 1$.
Let $\textbf{K} \in r B_{\ell^{\alpha}}$, and set $\textbf{Y} = \textbf{X} + n^{\frac{1}{2}+\frac{1}{d}} \textbf{K}$.
As we assumed $P$ is a monomial of total degree $k$, from the preceding Lemma \ref{meanvalue}, we have for any $\textbf{H} \in (\mathcal{H}_n^{(\beta)})^p$,
\begin{equation*} |f_m(\textbf{Y} +\textbf{H}) -f_m(\textbf{Y})| \leq \frac{c}{n}\Big( ||\textbf{Y}/\sqrt{n}||_{2(k-1)}^{k-1} +||\textbf{H}/\sqrt{n}||_2^{k-1}\Big) || \textbf{H}/\sqrt{n}||_2.
\end{equation*}
where $c$ is some constant depending $p$ and $d$. Using the fact that $x^{k-1} \leq 1+ x^{d-1}$ for any $1\leq  k \leq d$ and $x\geq 0$, we get,
\begin{align*}
|f_m(\textbf{Y} +\textbf{H}) -f_m(\textbf{Y})| &\leq \frac{c}{n}( ||\textbf{Y}/\sqrt{n}||_{2(k-1)}^{k-1}  +1 ) ||\textbf{H}/\sqrt{n}||_2\\
&+ \frac{c}{n}|| \textbf{H}/\sqrt{n}||_2^d.
\end{align*}
Let $\delta\in (0,1)$ and $t_{\delta} =  \delta n^{\frac{1}{2}+\frac{1}{d}}$. For $\textbf{H} \in  t_{\delta} B_{\ell^2}$,
$$ |f_m(\textbf{Y} +\textbf{H}) -f_m(\textbf{Y})| \leq 2c\delta (n^{\frac{1}{d}-1}||\textbf{Y}/\sqrt{n}||_{2(k-1)}^{k-1}  +1 ).$$
With the notation of Theorem \ref{theoremgene}, we have
$$\EE \sup_{\textbf{H} \in t_{\delta} B_{\ell^2}} \mathcal{L}_m(\textbf{H}) \leq 2c\delta (n^{\frac{1}{d}-1} \EE||\textbf{Y}/\sqrt{n}||_{2(k-1)}^{k-1}  +1 ),$$
where $m$ is the dimension of $(\mathcal{H}_n^{(\beta)})^p$. By convexity, we deduce 
$$ \EE ||\textbf{Y}/\sqrt{n}||_{2(k-1)}^{k-1}  \leq2^{k-2} \EE ||\textbf{X}/\sqrt{n}||_{2(k-1)}^{k-1}  +2^{k-2}n^{\frac{k-1}{d}} ||\textbf{K}||_{2(k-1)}^{k-1}.$$
But by Wigner's theorem (see \cite[Lemma 2.1.6]{Guionnet}),
$$\EE ||\textbf{X}/\sqrt{n}||_{2(k-1)}^{k-1}   \leq c_0n^{1/2},$$
for some constant $c_0>0$. As $\textbf{K} \in r B_{\ell^{\alpha}}$ with $\alpha \leq 2$, we deduce as $k\geq 2$,
$$ || \textbf{K}||_{2(k-1)}\leq || \textbf{K}||_2\leq r.$$
Thus, 
$$\EE \sup_{\textbf{H} \in  t_{\delta} B_{\ell^2}} \mathcal{L}_m(\textbf{H}) \leq C\delta (n^{\frac{1}{d}+\frac{1}{2}-1}  +r^{d-1}).$$
where  $C$ is some positive constant depending on $p$ and $d$. This shows that assumption $(ii)$ is satisfied.

We show now that assumption $(iii)$ holds.
Using \eqref{comptrnorm} for $q=d$, we get
$$|\tr P_d(\textbf{H})| \leq C' || \textbf{H}||_{\alpha}^{d/\alpha},$$
where $C'$ is some constant depending on $P$.
This proves condition $(iii)$ of Theorem \ref{theoremgene}.
To show  that the last assumption $(iv)$ is met, it suffices to observe that for any fixed $\textbf{H} \in (\mathcal{H}_n^{(\beta)})^p$, with the same construction as in \eqref{defH}, there is a sequence $\textbf{H}_m \in (\mathcal{H}_m^{(\beta)})^p$, for $m\geq n$, such that 
$$ \tr P_d(\textbf{H}_m) = \tr P_d(\textbf{H}),$$
and $W_{\alpha}(\textbf{H}) = W_{\alpha}(\textbf{H}_m)$.

\end{proof}

\section{Application to last-passage time}\label{LPPLDP}
We prove in this last section Theorem \ref{LDPLPP}. 
\begin{proof}[Proof of Theorem \ref{LDPLPP}] We will verify the assumptions of Theorem \ref{theoremgenesup}. Assumption $(i)$ holds due to Proposition \ref{convunifLPP} with $v(n) = n^{\alpha}$, and
$$\forall X \in \RR^I, \ f_m(X) = \frac{1}{n}T(X^+),\  F_m(X) = \mathcal{T}_n(X),$$
where $\mathcal{T}_n$ is defined in \eqref{defequidetLPP}, $X^+$ denotes the matrix with coefficients $(X^+_v)_v$, and $m$ is the dimension of $\RR^I$. 
 As $$X\mapsto T(X^+)/n,$$
 is $n^{-1/2}$-Lipschitz with respect to $|| \ ||_{\ell^{2}}$, assumption $(ii)$ is satisfied by the 
remark \ref{remtheogene} \hyperref[remlipconst]{(c)}. 

Using the fact that $|| \ ||_{\ell^1} \leq ||\ ||_{\ell^{\alpha}}$ when $\alpha \leq 1$, on $\RR^I$, we see that the condition $(iii)$ of Theorem \ref{theoremgenesup} is met. To prove $(iv)'$, we first observe that 
\begin{equation} \label{claimL} L_{\alpha}(x) = \inf \{ || H ||_{\ell^{\alpha}}^{\alpha} : \mathcal{T}_n(H) = x, H\in \RR^I\}.\end{equation}
Indeed, since the function $g$ is superadditive by \cite[Proposition 2.1]{Martin}, we deduce that 
$$\mathcal{T}_n(H) \geq g(1,...,1),$$
for any $H\in \RR^I$. Therefore, both sides of \eqref{claimL} are infinite if $x<g(1,...,1)$. Now if $x\geq g(1,1)$, and $ H\in \RR^	I$ is such that $\mathcal{T}_n(H) = x$, then denoting $\{v_0,...,v_p\}$ the element of $\mathcal{A}_m$ achieving the supremum in  \eqref{defequidetLPP}, we get,
$$
||H||_{\ell^{\alpha}}^{\alpha} \geq  \Big( \sum_{i=0}^{p-1} H_{v_i}^+\Big)^{\alpha} = \Big( x- \sum_{i=0}^{p-1} g \Big( \frac{v_{i+1}-v_i}{n}\Big)\Big)^{\alpha}.$$
Using the superadditivity of $g$, it yields
$$
||H||_{\ell^{\alpha}}^{\alpha} \geq (x-g(1,...,1))^{\alpha},$$
with equality for the matrix $H$ whose entries are all zero except $H_{(n,...,n)} = x-g(1,1)$. This proves the equality \eqref{claimL}. In particular, $L_{\alpha}$ is lower semi-continuous and therefore by the remark \ref{remtheogene} \hyperref[remlsi]{(e)}, we deduce,
$$ L_{\alpha}(x) = \sup_{\delta>0}\inf \{ || H ||_{\ell^{\alpha}}^{\alpha} : |\mathcal{T}_n(H) - x|<\delta, H\in \RR^I\}.$$
As the matrices $H \in \RR^I$ with  $H_{v}= (x-g(1,...,1))_+\Car_{v = (n,...,n)}$, achieves \eqref{claimL} for any $n$, we deduce, 
$$ L_{\alpha}(x) = \sup_{\delta>0} \limsup_{n\to +\infty} \inf\{ || H ||_{\ell^{\alpha}}^{\alpha} : |\mathcal{T}_n(H) - x|<\delta, H\in \RR^I\}.$$
Finally, as $\mathcal{T}_n(H) = \mathcal{T}_n(H^+)$, where $H^+$ is the matrix $(H_v^+)_{v\in \{0,...,n\}^d}$, we get
$$ L_{\alpha}(x) = \sup_{\delta>0} \limsup_{n\to +\infty} \inf\{ || H ||_{\ell^{\alpha}}^{\alpha} : |\mathcal{T}_n(H) - x|<\delta, H\in \RR_+^I\}.$$
 This proves the last assumption $(iv)'$ of Theorem \ref{theoremgenesup}.

\end{proof}

\newpage
\bibliographystyle{plain}
\bibliography{main.bib}{}

\end{document}